\documentclass[11pt,reqno,a4paper]{amsart}

\usepackage[colorlinks=true,linkcolor=blue]{hyperref}
\usepackage{cleveref}

\usepackage{mathrsfs}
\usepackage{amsmath, amssymb, amsthm}
\usepackage{amsfonts}
\usepackage{graphicx}
\makeatother
\usepackage{hyperref}
\usepackage{geometry}
\usepackage{enumerate}
\usepackage{etoolbox}
\usepackage{amscd}
\usepackage{bold-extra}
\usepackage[all,cmtip]{xy}
\usepackage[latin1]{inputenc}

\patchcmd{\subsection}{-.5em}{.5em}{}{}
\patchcmd{\subsubsection}{-.5em}{.5em}{}{}

\numberwithin{equation}{section}

\newcommand{\SL}{\operatorname{SL}}

\newcommand{\PGL}{\operatorname{PGL}}
\newcommand{\GL}{\operatorname{GL}}

\newcommand{\Aut}{\operatorname{Aut}}

\newcommand{\im}{\operatorname{Im}}

\newcommand{\cF}{\mathcal{F}}

\newcommand{\cH}{\mathcal{H}}

\newcommand{\cO}{\mathcal{O}}

\newcommand{\cU}{\mathcal{U}}

\newcommand{\cW}{\mathcal{W}}

\newcommand{\bF}{\mathbb{F}}

\newcommand{\bK}{\mathbb{K}}
\newcommand{\bL}{\mathbb{L}}

\newcommand{\bP}{\mathbb{P}}
\newcommand{\bQ}{\mathbb{Q}}
\newcommand{\bR}{\mathbb{R}}
\newcommand{\bS}{\mathbb{S}}
\newcommand{\bT}{\mathbb{T}}

\newcommand{\bZ}{\mathbb{Z}}

\newcommand{\ra}{\rightarrow}

\newcommand{\qand}{\quad \textrm{and} \quad}

\def\acts{\curvearrowright}
\newcommand\subsetsim{\mathrel{%
\ooalign{\raise0.2ex\hbox{$\subset$}\cr\hidewidth\raise-0.8ex\hbox{\scalebox{0.9}{$\sim$}}\hidewidth\cr}}}
\newcommand{\eps}{\varepsilon}

\DeclareMathOperator{\Hom}{Hom}

\DeclareMathOperator{\pr}{pr}

\DeclareMathOperator{\supp}{supp}

\DeclareMathOperator{\Prob}{Prob}

\DeclareMathOperator{\Char}{char}
\DeclareMathOperator{\Prod}{Prod}

\theoremstyle{theorem}
\newtheorem{theorem}{Theorem}[section]

\newtheorem{corollary}[theorem]{Corollary}

\newtheorem{lemma}[theorem]{Lemma}

\theoremstyle{definition}
\newtheorem{definition}[theorem]{Definition}
\newtheorem{remark}[theorem]{Remark}

\usepackage{changepage}
\usepackage{mathtools}
\usepackage{graphicx}
\usepackage{calc}

\usepackage[toc,page]{appendix}

\usepackage{scalerel}

\usepackage{extarrows}

\DeclareMathSymbol{\shortminus}{\mathbin}{AMSa}{"39}
\usepackage{mathtools}
\begin{document}
\bibliographystyle{plain} 

\title[Dynamics of Multiplicative Groups over Fields]{Dynamics of Multiplicative Groups over Fields and F\o lner-Kloosterman Sums}

\author{Michael Bj\"orklund}
\address{Department of Mathematics, Chalmers, Gothenburg, Sweden}
\email{micbjo@chalmers.se}
\thanks{}

\author{Alexander Fish}
\address{School of Mathematics and Statistics F07, University of Sydney, NSW 2006,
Australia}
\curraddr{}
\email{alexander.fish@sydney.edu.au}
\thanks{}


\begin{abstract}
For two countably infinite fields whose multiplicative groups are isomorphic, we examine invariant couplings between the actions that these groups induce on the additive Pontryagin duals of the fields. We show that the actions are disjoint unless the fields themselves are isomorphic and the group isomorphism extends (possibly after a finite twist) to a field isomorphism. As an application, we establish equidistribution of F\o lner--Kloosterman sums - an extension of classical Kloosterman sums to infinite fields. Unlike the classical case over algebraic closures of finite fields, these averages exhibit an inherent multiplicative asymmetry, revealing new and fundamentally different behavior. Finally, we derive several combinatorial consequences, including results on sum-product phenomena and a Furstenberg--S\'ark\"ozy-type theorem for Laurent polynomials over general fields.
\end{abstract}

\maketitle

\section{Introduction}

\subsection{Reconstructing fields from their multiplicative groups}

The multiplicative group of a field can be remarkably uninformative: even highly non-isomorphic fields may share the same group structure. For algebraically closed fields, whose multiplicative groups are clearly divisible, this is already visible, since
\[
\bK^{*} \cong \mu(\bK)\oplus \mathbb{Q}^{(\kappa)}
\]
with the torsion part determined by the characteristic and the rank by the cardinality; hence any two algebraically closed fields of the same characteristic and cardinality automatically have isomorphic multiplicative groups.

More surprisingly, this occurs for very concrete fields as well. The fields $\mathbb{Q}(x)$ and $\mathbb{Q}(x,y)$ differ drastically in transcendence degree, yet both are countable-dimensional over $\mathbb{Q}$ and admit unique factorizations over countably many irreducibles, making both their additive and multiplicative groups isomorphic. Even fields of different characteristics can exhibit the same multiplicative behavior:
\[
\mathbb{Q}^{*} \cong \{\pm 1\} \oplus \bigoplus_{p\in P} p^{\mathbb{Z}}
\qquad\text{and}\qquad
\mathbb{F}_3(x)^{*} \cong \{\pm 1\} \oplus \bigoplus_{f\in P(x)} f^{\mathbb{Z}},
\]
and since both indexing sets are countable, one obtains the counterintuitive isomorphism
\[
\mathbb{Q}^{*} \cong \mathbb{F}_3(x)^{*}.
\]
This raises a natural and long-standing question: when does a group isomorphism between the multiplicative groups of two fields reflect a deeper algebraic 
connection - specifically, does it extend (perhaps after a mild twist) to a field isomorphism? A central goal is to identify invariants capable of distinguishing non-isomorphic fields that nonetheless have isomorphic multiplicative groups.

A foundational result in this direction is the Neukirch-Uchida-Pop theorem \cite{N,P,U}, which shows that the absolute Galois group determines a number field (and more generally, fields finitely generated over a base). By Kummer theory, an isomorphism of absolute Galois groups induces one of multiplicative groups, and the main challenge lies in recovering the full field structure.

In a related development, Cornelissen, de Smit, Li, Marcolli and Smit \cite{CLMS, CdSLMS,CM} proved that an isomorphism between the abelianized Galois groups of two number fields that also preserves their Dirichlet 
L-functions must arise from a field isomorphism. Further contributions by Bogomolov and Tschinkel \cite{BT} have shown that for various classes of fields, small Galois-theoretic or K-theoretic invariants can rigidly determine the field.

\subsection{Disjointness versus twisted field isomorphisms}

In this paper we develop a dynamical counterpart of this program. We begin by introducing some notation. Let $\bK$ be a discrete countably infinite field. The Pontryagin dual of the additive group $(\bK,+)$,
\[
\widehat{\bK}=\Hom(\bK,\bS^1),
\]
is a compact second countable group, and there is a natural injective homomorphism
\[
\bK^* \longrightarrow \Aut(\widehat{\bK}), \qquad 
a \mapsto a^*,
\]
where $(a^*\xi)(b)=\xi(ab)$ for $\xi\in\widehat{\bK}$ and $b\in\bK$.  
The associated $\bK^*$-action on $\widehat{\bK}$ has exactly two ergodic invariant probability measures: the point mass $\delta_1$ at the trivial character, and the Haar probability measure $m_{\widehat{\bK}}$. Since $\bK$ is infinite, $m_{\widehat{\bK}}$ is non-atomic.

Now let $\bL$ be another countably infinite field such that $\bK^*$ and $\bL^*$ are isomorphic, and fix an isomorphism $\rho:\bK^*\to\bL^*$. Via~$\rho$, the group $\bK^*$ acts on $\smash{\widehat{\bL}}$, and again the only ergodic invariant probability measures are $\delta_1$ and $m_{\widehat{\bL}}$. Our first objective is to compare the two ergodic $\bK^*$-actions
\[
\bK^* \acts (\widehat{\bK},m_{\widehat{\bK}})
\qquad\text{and}\qquad
\bK^* \acts_\rho (\widehat{\bL},m_{\widehat{\bL}}).
\]

A probability measure on $\widehat{\bK}\times\widehat{\bL}$ whose marginals are $m_{\widehat{\bK}}$ and $m_{\widehat{\bL}}$ is called a \textbf{coupling}. Such a coupling is \textbf{invariant} if it is invariant under the diagonal $\bK^*$-action. The product measure $m_{\widehat{\bK}}\otimes m_{\widehat{\bL}}$ is always an invariant coupling, and we refer to it as the \textbf{trivial coupling}. We say that the two $\bK^*$-actions are \textbf{disjoint} if the trivial coupling is the only invariant coupling. The actions are said to be \textbf{measurably isomorphic} if there exists a measurable map $\pi:\smash{\widehat{\bK}}\to \smash{\widehat{\bL}}$, defined almost everywhere, such that $\pi_* m_{\widehat{\bK}} = m_{\widehat{\bL}}$ and
\[
\pi(a^*\xi)=\rho(a)^*\pi(\xi)
\qquad\text{for $m_{\widehat{\bK}}$-almost every $\xi\in\widehat{\bK}$.}
\]
Whenever such a map exists, the measure $\mu_\pi\in\Prob(\widehat{\bK}\times\widehat{\bL})$ defined by
\[
\mu_\pi(f)
=\int_{\widehat{\bK}} f(\xi,\pi(\xi))\, dm_{\widehat{\bK}}(\xi),
\qquad f\in C(\widehat{\bK}\times\widehat{\bL}),
\]
is a non-trivial $\bK^*$-invariant coupling. Thus, disjointness is a particularly strong form of non-isomorphism. 

As discussed in the previous subsection, a central theme is to understand whether an isomorphism between the multiplicative groups of two fields reflects a deeper algebraic relationship, and which invariants can distinguish fields whose multiplicative groups nonetheless coincide. In this spirit, the main theorem below shows that the actions of the multiplicative groups of two countably infinite non-isomorphic fields on their respective additive duals cannot be measurably isomorphic, even when the multiplicative groups themselves are isomorphic.

\begin{theorem}
\label{Thm_main1}
Let \(\bK\) and \(\bL\) be countably infinite fields such that \(\bK^*\) and \(\bL^*\) are isomorphic as groups. Let \(\rho : \bK^* \to \bL^*\) be an isomorphism, and suppose that the ergodic \(\bK^*\)-actions
\[
\bK^* \acts (\widehat{\bK}, m_{\widehat{\bK}}) \quad \text{and} \quad \bK^* \acts_\rho (\widehat{\bL}, m_{\widehat{\bL}}),
\]
where \(\bK^*\) acts on \(\widehat{\bL}\) via \(\rho\), are not disjoint. Then \(\bK\) and \(\bL\) are isomorphic as fields. More precisely, there exists a homomorphism \(w : \bK^* \to \bL^*\) with finite image such that the map \(\kappa : \bK \to \bL\) defined by
\[
\kappa(0) = 0, \quad \text{and} \quad \kappa(x) = w(x)\rho(x) \quad \text{for all } x \in \bK \setminus \{0\},
\]
is a field isomorphism. Furthermore, any non-trivial $\bK^*$-invariant and $\bK^*$-ergodic coupling
of the $\bK^*$-actions on $(\smash{\widehat{\bK}},m_{\widehat{\bK}})$ and 
$(\smash{\widehat{\bL}},m_{\widehat{\bL}})$ is of the form $(e,b^*)_*\mu$ for some $b \in \bL^*$, 
where 
\[
\mu(f)
= \frac{1}{|\im(w)|}
\sum_{c \in \im(w)}
\int_{\widehat{\bL}} f(\kappa^*(\eta),c^*\eta)\, dm_{\widehat{\bL}}(\eta),
\qquad f \in C(\widehat{\bK}\times\widehat{\bL}).
\]
and $\kappa^* : \widehat{\bL} \ra \widehat{\bK}, \enskip \eta \mapsto \eta \circ \kappa$ denotes the $\kappa$-induced isomorphism  between $\smash{\widehat{\bL}}$ and $\smash{\widehat{\bK}}$.
\end{theorem}

\begin{remark}
We prove this theorem in Section \ref{Sec:ThmMain1}.
\end{remark}

\begin{remark}
Let us now give an example of a pair $(\kappa,w)$ with $w$ non-trivial in the case $\mathbb K=\mathbb Q$. Every $x \in \bQ^*$ can be written uniquely as 
\[
x = (-1)^m p_1^{n_1} p_2^{n_2}\cdots, 
\]
where $m\in\{0,1\}$, where $p_1,p_2,\ldots$ is a fixed 
enumeration of the primes, and where the sequence $i \mapsto n_i$ has finite support. Define 
$\kappa(x)=x$ and let $w(x)=(-1)^{n_1}$. Then $w$ is a non-trivial multiplicative 
homomorphism with finite image, $\kappa$ is the (identity) field automorphism of $\mathbb Q$, and the 
map $\rho(x)=w(x)^{-1}\kappa(x)$ is an isomorphism of $\mathbb Q^{*}$.
\end{remark}

\begin{remark}
It follows from our discussion in Subsection~\ref{subsec:measure-rigid}, and is not difficult to prove, that the (mixing) \( \mathbb{K}^* \)-actions
\[
\mathbb{K}^* \curvearrowright (\widehat{\mathbb{K}}, m_{\widehat{\mathbb{K}}}) 
\quad \text{and} \quad 
\mathbb{K}^* \curvearrowright_\rho (\widehat{\mathbb{L}}, m_{\widehat{\mathbb{L}}})
\]
are both \emph{2-simple}. By a classical result of del Junco, Rudolph, and Veech - surveyed in \cite[Chapter 12]{G} - this implies that the only non-trivial \( \mathbb{K}^* \)-factors of these actions arise as quotients by \emph{finite} subgroups of \( \mathbb{K}^* \). Moreover, any joining between these actions must be relatively independent over a common factor.

As a result, rigidity reduces (to a large extent) to showing that any measurable \( \mathbb{K}^* \)-equivariant isomorphism between \( \smash{\widehat{\mathbb{K}}} \) and \( \smash{\widehat{\mathbb{L}}} \) must arise from a field isomorphism.
\end{remark}

\begin{remark}
When \(\mathbb{K}\) and \(\mathbb{L}\) are number fields, a strong 
classification of such joinings - under significantly weaker assumptions - was 
obtained by Einsiedler and Lindenstrauss \cite{EL}. Their approach relies on 
entropy methods. In contrast, our argument is primarily Fourier-analytic in 
nature. This enables us to work in much greater generality (in particular, 
beyond number fields) and to compare actions of multiplicative groups of 
fields that may differ dramatically, even in characteristic. However, this 
comes at the cost of requiring stronger assumptions on the size of the 
acting groups.
\end{remark}

The key mechanism behind our proof (see Theorem~\ref{Thm_CriterionCR}) is the failure of \emph{3-mixing} for the action 
\[
\mathbb{K}^* \curvearrowright (\widehat{\mathbb{K}}, m_{\widehat{\mathbb{K}}}).
\]
To understand this, observe that for any \( a \in \mathbb{K}^* \), the identity
\[
1 + (a - 1) + (-a) = 0
\]
implies the functional identity
\[
\int_{\widehat{\mathbb{K}}} f(\xi) f((a - 1).\xi) f((-a).\xi) \, dm_{\widehat{\mathbb{K}}}(\xi) = 1, 
\quad \text{for all } a \in \mathbb{K}^* \setminus \{1\},
\]
where \( f(\xi) := \xi(1) \) satisfies \( m_{\widehat{\mathbb{K}}}(f) = 0 \). This shows that the triple correlations remain nontrivial for all such \( a \), and in particular, no sequence \( (a_n) \subset \mathbb{K}^* \) diverging to infinity can decorrelate the triple product. Hence, the action fails to be 3-mixing. This lack of 3-mixing enables us to extract many linear relations among the values 
\[
\rho(a_1), \rho(a_2), \rho(a_3) 
\quad \textrm{for triples \((a_1, a_2, a_3) \in (\mathbb{K}^*)^3\) satisfying \( a_1 + a_2 + a_3 = 0 \)}. 
\]
These relations ultimately force \( \rho \) to behave like an "almost" additive homomorphism. \\

In Appendix A - the main technical component of the paper - we show that any such 
"almost homomorphism" must in fact arise from a genuine field isomorphism. Let us now deduce a simple corollary of Theorem \ref{Thm_main1}. 

\begin{corollary}
\label{Cor_main1}
Let $\bK$ be a countably infinite field and consider the $\bK^*$-action on 
$\widehat{\bK} \times \widehat{\bK}$ given by
\[
a.(\xi_1,\xi_2) = (a^*\xi_1,(a^{-1})^*\xi_2), \quad a \in \bK^*, \enskip (\xi_1,\xi_2) \in \widehat{\bK} \times \widehat{\bK}.
\]
Then this action has exactly four ergodic invariant probability measures:
\[
\delta_1 \otimes \delta_1, \quad \delta_1 \otimes m_{\widehat{\bK}}, \quad 
m_{\widehat{\bK}} \otimes \delta_1, \quad m_{\widehat{\bK}}\otimes m_{\widehat{\bK}}.
\]
\end{corollary}

\begin{proof}
Let $\theta$ be an ergodic, invariant probability measure for the $\bK^*$-action on $\widehat{\bK} \times \widehat{\bK}$. Since the projections of $\theta$ onto each factor are themselves ergodic, they must be either $\delta_1$ or the Haar measure $m_{\widehat{\bK}}$. If at least one of these projections is equal to $\delta_1$, then it is clear that $\theta$ must fall into one of the first three classified cases. Thus, we may assume that both projections are equal to the Haar measure, and in particular, $\theta$ is an invariant coupling of Haar measures.

By Theorem \ref{Thm_main1}, applied to the (multiplicative) isomorphism $\rho(a) = a^{-1}$, the product measure is the unique invariant coupling, provided there does \emph{not} exist a homomorphism $w: \bK^* \to \bK^*$ with finite image and a field isomorphism $\kappa: \bK \to \bK$ such that
\[
\frac{1}{a} = w(a) \kappa(a), \quad \text{for all } a \in \bK^*.
\]
Suppose, for the sake of contradiction, that such $w$ and $\kappa$ do exist. Substituting $1 - a$ in place of $a$ in the above identity yields:
\[
\frac{1}{1 - a} = w(1 - a) (1 - \kappa(a)) = w(1 - a) - \frac{w(1 - a)}{w(a)} \cdot \frac{1}{a}, \quad \text{for all } a \in \bK^* \setminus \{1\}.
\]
Since $w$ has finite image, there exists $c \in \bK^*$ such that the set
\[
S = \left\{ a \in \bK^* \setminus \{1\} : \frac{w(1 - a)}{w(a)} = c \right\}
\]
is infinite. For each $a \in S$, the identity above implies
\[
\frac{1}{1 - a} + \frac{c}{a} \in \mathrm{im}(w),
\]
which is impossible: the map $a \mapsto \frac{1}{1 - a} + \frac{c}{a}$ is at most two-to-one, so its image over the infinite set $S$ cannot be contained in the finite set $\mathrm{im}(w)$.
\end{proof}

\subsection{F\o lner-Kloosterman sums}

Let us now discuss an application of Corollary \ref{Cor_main1}. We begin by recalling the definition of the (classical) Kloosterman sums, which are invaluable tools in analytic number theory. Let $\bF_q$ denote the finite field with $q$ elements of characteristic $p$, and fix an algebraic closure $\smash{\overline{\bF}_p}$. For each $n$, we have a natural field inclusion $\smash{\bF_{p^n} \hookrightarrow \overline{\bF}_p}$. Given $(\xi_1,\xi_2) \in \smash{\widehat{\overline{{\bF}}}_p \times \widehat{\overline{{\bF}}}_p}$ and an integer $n \geq 1$, we define the 
(classical) \textbf{Kloosterman sum}
\[
K_{p^n}(\xi_1,\xi_2) := \sum_{a \in \bF_{p^n}^*} \xi_1(a) \xi_2(a^{-1}).
\]
In applications, it is often important to understand the asymptotics of $K_{p^n}(\xi_1,\xi_2)$ as either $p$ or $n$ tends to infinity (or both), and how it depends on the (additive) characters $\xi_1$ and $\xi_2$. In Deligne's classical work \cite{D}, very strong bounds are established for the sums $K_{p^n}(\xi_1,\xi_2)$. As a very special 
case, it follows from these bounds that for every $p$,
\begin{equation}
\label{equiKlooster}
\lim_{n \ra \infty} \frac{1}{|\bF_{p^n}|} \sum_{a \in \bF_{p^n}^*} \xi_1(a) \xi_2(a^{-1}) = 0, \quad \textrm{for all $(\xi_1,\xi_2) \in \smash{\widehat{\overline{{\bF}}}_p \times \widehat{\overline{{\bF}}}_p}$}.
\end{equation}
The subsets $(\bF_{p^n})$ of $\overline{{\bF}}_p$ form an additive F\o lner sequence,
which is also asymptotically invariant under multiplication, as well as symmetric under multiplicative inversion. \\

It turns out that for any countably infinite field \( \bK \), one can construct
sequences of finite subsets that are asymptotically invariant under both addition
and multiplication - so-called \textbf{double F\o lner sequences} \cite[Proposition 2.4]{BM}.
Given such a sequence \( (F_n) \) in \( \bK \), one readily checks that
\begin{equation}
\label{easydouble}
\lim_{n \to \infty} \frac{1}{|F_n|} \sum_{a \in F_n} \xi(a) = 0,
\qquad \text{for every } \xi \in \widehat{\bK} \setminus \{1\}.
\end{equation}
Corollary~\ref{Cor_main1} yields the following generalization of
\eqref{equiKlooster}.

\begin{theorem}[F\o lner-Kloosterman sums]
\label{Thm_main2}
Let $\bK$ be a countably infinite field and let $(F_n)$ be a double F\o lner sequence in $\bK$. Then, for all $(\xi_1,\xi_2) \in \smash{\widehat{\bK}} \times \smash{\widehat{\bK}}$
with $\xi_1 \neq 1$, we have
\[
\lim_{n \ra \infty} \frac{1}{|F_n|} \sum_{a \in F_n \setminus \{0\}} \xi_1(a) \xi_2(a^{-1}) = 0.
\]
\end{theorem}

\begin{remark}
We refer to sums as in the theorem as \textbf{F\o lner-Kloosterman sums}.
\end{remark}

\begin{proof}
Fix \((\xi_1, \xi_2) \in \widehat{\bK} \times \widehat{\bK}\), and let \((F_n)\) be a double F\o lner sequence. Consider the sequence \((\mu_n)\) of probability measures on \(\widehat{\bK} \times \widehat{\bK}\) defined by
\[
\mu_n = \frac{1}{|F_n|} \sum_{a \in F_n \setminus \{0\}} \delta_{a^*\xi_1} \otimes \delta_{(a^{-1})^*\xi_2}, \quad n \geq 1.
\]
Our goal is to show that if \(\xi_1 \neq 1\), then
\[
\lim_{n \to \infty} \widehat{\mu}_n(1,1) = 0,
\]
where \(\widehat{\mu}_n\) denotes the Fourier transform of \(\mu_n\). To prove this, consider a weak$^*$-convergent subsequence \((\mu_{n_k})\) with limit \(\mu\). Since \((F_n)\) is asymptotically invariant under multiplication, \(\mu\) is invariant under the \(\bK^*\)-action. Consequently, by Corollary~\ref{Cor_main1}, \(\mu\) must be of the form
\[
\mu = A \cdot \delta_1 \otimes \delta_1 + B \cdot \delta_1 \otimes m_{\widehat{\bK}} + C \cdot m_{\widehat{\bK}} \otimes \delta_1 + D \cdot m_{\widehat{\bK}} \otimes m_{\widehat{\bK}},
\]
for some non-negative real numbers \(A, B, C, D\) satisfying \(A + B + C + D = 1\). Now observe that \(\widehat{\mu}(1,1) = A\), since \(\widehat{m}_{\bK}(1) = 0\). To show that \(A = 0\), we use the identity from \eqref{easydouble}, which gives
\[
\lim_{k \to \infty} \widehat{\mu}_{n_k}(1,0) = \lim_{k \to \infty} \frac{1}{|F_{n_k}|} \sum_{a \in F_{n_k} \setminus \{0\}} \xi_1(a) = 0.
\]
But since \( \lim_k \widehat{\mu}_{n_k}(1,0) = A + B\), it follows that \(A + B = 0\). As both \(A\) and \(B\) are non-negative, this forces \(A = B = 0\).

Since the choice of weak$^*$-convergent subsequence was arbitrary, it follows that the full sequence \((\widehat{\mu}_n(1,1))\) converges to \(0\), as desired.
\end{proof}

\subsection{On the inherent multiplicative asymmetry of double F\o lner sequences}

The last theorem and its proof naturally raise the question: is the assumption \(\xi_1 \neq 1\) truly necessary? In the proof, it was used to show that \(\lim_k \widehat{\mu}_{n_k}(1,0) = 0\), but could we not instead consider \(\lim_k \widehat{\mu}_{n_k}(0,1)\)? Would this limit also vanish? \\

Explicitly, this amounts to showing that for any double F\o lner sequence \((F_n)\),
\begin{equation}
\label{fail}
\lim_{n \to \infty} \frac{1}{|F_n|} \sum_{a \in F_n \setminus \{0\}} \xi(a^{-1}) = 0.
\end{equation}
This holds when \((F_n)\) is symmetric under inversion - as in the case \(\bK = \overline{\mathbb{F}}_p\) with \(F_n = \mathbb{F}_{p^n}\) discussed above. But what about the general case? \\

Our next result tells us that multiplicative symmetry of double F\o lner sequences is a very rare phenomenon.

\begin{theorem}
\label{Thm_MultAssym_Main}
Let \(\bK\) be a countable field and suppose that there exists a F\o lner sequence \((F_n)\) in \((\bK,+)\) such that \(((F_n \setminus \{0\})^{-1})\) is also a F\o lner sequence. Then \(\bK\) is a union of finite subfields. 
\end{theorem}

\begin{remark}
We prove Theorem \ref{Thm_MultAssym_Main} in Section \ref{Sec:MultAssym}.
\end{remark}

In particular, this theorem implies that multiplicatively symmetric double F\o lner sequences cannot exist in fields of characteristic zero. Consequently, one should not expect the limit \eqref{fail} to hold for purely automatic reasons in such fields. In fact, in the case of real subfields, we can demonstrate that \eqref{fail} fails in a particularly strong sense:

\begin{theorem}
\label{Thm_WeirdInverseMain}
Let $\bK$ be a countable sub-field of $\bR$. Then, for every $\xi \in \widehat{\bR}$
and every double F\o lner sequence $(F_n)$, the sequence
\[
u_n := \frac{1}{|F_n|} \sum_{a \in F_n \setminus \{0\}} \xi(a^{-1}), \quad n \geq 1.
\]
does not converge to zero as $n \ra \infty$. 
\end{theorem}

\begin{remark}
We prove Theorem \ref{Thm_WeirdInverseMain} in Section \ref{Sec:Failure}.
\end{remark}

\begin{remark}
It is unclear whether the limit of $(u_n)$ exists for any real character $\xi$; however, if it does, the theorem above tells us that it cannot be zero. Computing this limit, if it exists, for a given $\xi$, would certainly be of interest.
\end{remark}

\begin{remark}
We emphasize that Theorem~\ref{Thm_WeirdInverseMain} only concerns convergence to zero for those $\xi \in \widehat{\bK}$ arising as restrictions of continuous characters on $\bR$. In fact, the $\bK^*$-action on $\widehat{\bK}$ defined by $a.\xi = (a^{-1})^*\xi$ is ergodic with respect to $m_{\widehat{\bK}}$, so the pointwise ergodic theorem gives
\[
\frac{1}{|F_n|} \sum_{a \in F_n \setminus \{0\}} \xi(a^{-1}) = 0
\quad \text{for $m_{\widehat{\bK}}$-almost every $\xi$}.
\]
Since the set of $\xi$ that are restrictions of continuous characters on $\bR$ is $m_{\widehat{\bK}}$-null, this does not contradict Theorem~\ref{Thm_WeirdInverseMain}.

To see that this set is null, suppose otherwise. Let $U$ be its image under the (injective) dual of the restriction map $\smash{\widehat{\bR}} \to \smash{\widehat{\bK}}$. Then $U$ is a Borel subgroup of positive Haar measure in the compact group $\widehat{\bK}$. By \cite[Chapter II, Theorem 14]{W}, any Borel subgroup of positive measure in a compact group is open. Since $\smash{\widehat{\bK}}$ is connected (\cite[Prop.~3, p.~63]{Morris}), this forces $U = \smash{\widehat{\bK}}$, implying the dual of the restriction map is a (topological) isomorphism $\bR \to \smash{\widehat{\bK}}$, which is impossible since the group $(\bR,+)$ is not compact.
\end{remark}

\subsection{Measure-classification for diagonalizable representations}

Our aim now is to situate Corollary \ref{Cor_main1} within a broader framework and to present a substantial generalization, along with several applications to Ramsey theory. \\

Let $\bK$ be a discrete countably infinite field. Let $U$ be a countable and discrete (possibly infinite-dimensional) vector space over $\bK$, and $\rho_U : \bK^* \ra \GL(U)$ a linear representation. Let $\smash{\widehat{U}}$ denote the Pontryagin dual of $U$, and 
$\rho_U^* : \bK^* \ra \Aut(\widehat{U})$ the \textbf{dual representation} given by
\[
(\rho_U(a)^*.\xi)(u) = \xi(\rho_U(a)u), \quad \textrm{for $a \in \bK^*, \xi \in \widehat{U}$ and $u \in U$}.
\]
Our goal in this subsection is to classify $\rho^*_U$-invariant and $\rho_U^*$-ergodic
probability measures on the compact and metrizable group $\smash{\widehat{U}}$, under some mild technical assumptions. \\

Let us now explicate what these assumptions are. We say that the representation $(U,\rho_U)$ is \textbf{diagonalizable} if there are linear sub-spaces $U_1,\ldots,U_r \leq U$ and distinct integers $m_1,\ldots,m_r$ such that 
\[
U = U_1 \oplus \cdots \oplus U_r \qand 
\rho_U(a)u_i = a^{m_i} u_i, \quad \textrm{for all $i=1,\ldots,r$}.
\]
We refer to $U_1,\ldots,U_r$ as the \textbf{weight spaces} of $(U,\rho_U)$, with 
\textbf{weights} $m_1,\ldots,m_r$ respectively, and we say that a diagonalizable $\bK^*$-representation $(U,\rho_U)$ is \textbf{reduced} if $m_1,\ldots,m_r$ are not divisible by the characteristics of $\bK$. If $\Char(\bK) = 0$, every diagonalizable
representation is reduced. \\

If $(U,\rho_U)$ is a diagonalizable $\bK^*$-representation, with weight spaces $U_1,\ldots,U_r$ and weights $m_1,\ldots,m_r$, then  
\[
\widehat{U} \cong \widehat{U}_1 \times \cdots \times \widehat{U}_r,
\]
and $\rho_U(a)^*\xi_i = (a^{m_i})^*.\xi_i$ for all $\xi_i \in \widehat{U_i}$ and $i=1,\ldots,r$. If $W_i \leq U_i$ is a linear-subspace, we define the 
\textbf{annihilator} $W_i^{\perp} < \widehat{U}_i$ of $W_i$ by
\[
W_i^{\perp} = \{ \xi_i \in \widehat{U}_i \, \mid \, \xi_i \mid_{W_i} = 1\}
\] 
and we note that every probability measure on $\widehat{U}$ of the form
\[
\theta = m_{W_1^{\perp}} \otimes \cdots \otimes m_{W_r^{\perp}},
\]
is $\rho_U^*$-invariant. Moreover, if all of the weights are non-zero, then every such measure is also $\rho_U^*$-ergodic. Our main result in this sub-section essentially states that when all weights are non-zero, then these measures exhaust the space of $\rho_U^*$-invariant and $\rho_U^*$-ergodic probability measures on 
$\smash{\widehat{U}}$. However, there is a small technical subtlety that appears for fields with positive characteristics. One way to  circumvent this issue is by introducing the following notion of a reduced model.

\begin{definition}
Let $(U,\smash{\rho_U})$ and $(V,\smash{\rho_V})$ be diagonalizable $\bK^*$-representations with weight spaces $U_1,\ldots,U_r$ and $V_1,\ldots,V_s$ respectively. 
If $(V,\rho_V)$ is reduced, and there is an isomorphism $\Phi : V \ra U$ between the 
additive groups $(V,+)$ and $(U,+)$ satisfying $\Phi \circ \rho_V(a) = \rho_U(a) 
\circ \Phi$ for all $a \in \bK^*$, we say that $(V,\rho_V,\Phi)$ is a \textbf{reduced model} of $(U,\rho_U)$.
\end{definition}

\begin{remark}
We emphasize that we do not assume that $\Phi$ is \emph{linear}, but only additive. 
\end{remark}

If $(V,\rho_V,\Phi)$ is a reduced model of $(U,\rho_U)$, then the push-forward of 
the inverse dual map 
\[
(\Phi^{-1})^* : \widehat{V} \ra \widehat{U}, \enskip \eta \mapsto \eta \circ \Phi^{-1}
\]
induces an affine isomorphism between the set of $\rho_V^*$-invariant probability measures on $\smash{\widehat{V}}$ and the set of $\rho_U^*$-invariant probability measures on $\smash{\widehat{U}}$, mapping ergodic measures to ergodic measures. In particular, if a reduced model $(V,\rho_V,\Phi)$ of $(U,\rho_U)$ can be found, then to classify the 
$\rho_U^*$-invariant probability measures it is enough to classify the $\rho_V^*$-invariant ones. The following two theorems now provide a complete measure classification for a general diagonalizable representation.

\begin{theorem}
\label{Thm_ReducedModels}
Every diagonalizable $\bK^*$-representation admits a reduced model. 
\end{theorem}

\begin{remark}
Theorem \ref{Thm_ReducedModels} is proved in Appendix B.
\end{remark}

\begin{theorem}
\label{Thm_DiagonalMain} 
Suppose \( (V, \rho_V) \) is a reduced diagonalizable $\bK^*$-representation with weight spaces $V_1,\ldots,V_s$. Then, for every \( \smash{\rho_V^*} \)-invariant and \( \smash{\rho_V^*} \)-ergodic probability measure \( \theta \) on \( \smash{\widehat{V}} \), there are linear subspaces 
\[
W_1 \leq V_1, \dots, W_s \leq V_s  \quad \textrm{such that}   \quad 
\theta =m_{W_1^{\perp}} \otimes \cdots \otimes m_{W_s^{\perp}}.
\]
\end{theorem}

\begin{remark}
Theorem \ref{Thm_DiagonalMain} is proved in Section \ref{Sec:ThmDiagonal}.
\end{remark}

\begin{remark}
When all of the weights are positive (or negative), the result follows from the van der Corput lemma, though the argument is far from straightforward when \( \bK \) has positive characteristic. For a proof in this special case, we refer the reader to \cite[Theorem 2.2.1, p.60]{L}. The novelty of Theorem \ref{Thm_Equin1nr} lies in its applicability to cases where 
the weights have \emph{mixed} signs. In this more general setting, the standard van der Corput approach is no longer effective, prompting the new techniques developed in this paper.
\end{remark}

We can now formulate the following generalization of Theorem \ref{Thm_main2}, which 
corresponds to the case $r = 2, n_1 = 1, n_2 = -1$ and $\eta = 1$.

\begin{theorem}
\label{Thm_Equin1nrMain}
Let $n_1,\ldots,n_r$ be non-zero integers, not divisible by the characteristics of $\bK$. Suppose that the set $I = \{ i=1,\ldots,r \, \mid \, n_i > 0\}$ is non-empty. Then, for every double F\o lner sequence $(F_k)$ in $(\bK,+)$ and 
$\xi_1,\ldots,\xi_r \in \smash{\widehat{\bK}}$ satisfying $\xi_i \neq 1$ for all $i \in I$ and for every multiplicative unitary character $\eta$ on $\bK^*$, we have
\[
\lim_{k \ra \infty} \frac{1}{|F_k|} \sum_{a \in F_k \setminus \{0\}} \eta(a) \, \xi_1(a^{n_1}) \cdots \xi_r(a^{n_r}) = 0.
\]
\end{theorem}

\begin{remark}
Theorem \ref{Thm_WeirdInverseMain} demonstrates that, in general, the condition 
$I \neq \emptyset$ is necessary. We prove Theorem \ref{Thm_Equin1nrMain} in 
Subsection \ref{subsec:Equi1nrMain}.
\end{remark}

\subsection{Pinned twisted patterns in sets of positive upper Banach density}

We now turn to some combinatorial consequences of the theorems above. The following
result can be viewed as an analogue of \cite[Theorem 1.10]{BB} by the first author 
and K. Bulinski, which in turn is a generalization of the investigations in \cite{BF1} by the authors.
Recall that if $A$ is a countable abelian group, $E \subseteq A$, and $(F_n)$ is a F\o lner sequence in $A$, then the \textbf{upper asymptotic density} of $E$ along $(F_n)$ is
\[
\overline{d}_{(F_n)}(E) = \varlimsup_{n\to\infty} \frac{|E \cap F_n|}{|F_n|}.
\]
The \textbf{upper Banach density} of $E$ is defined by
\[
d^*(E) = \sup\{\, \overline{d}_{(F_n)}(E) \mid (F_n)\ \text{is a F\o lner sequence in } A \,\}.
\]
In all the theorems below, when $\mathbb{K}$ is a countable infinite field, any reference to upper Banach density is taken with respect to the additive group $(\mathbb{K},+)$. In particular, we do not require the F\o lner sequences involved to be asymptotically invariant under multiplication.

\begin{theorem}[Pinned images of invariant functions]
\label{Thm_CombMain}
Let \( U \) be a countable vector space over \( \mathbb{K} \), and let \( \rho_U : \mathbb{K}^* \to \operatorname{GL}(U) \) be a diagonalizable representation with only non-zero integer weights. Let $Y$ be a set and suppose $\psi : U \ra Y$ is a $\rho_U$-invariant function. Then, for every $E \subset U$ with positive upper Banach density, 
the following holds: for every finite set $F \subset U$, there exists $u_F \in E$ 
such that
\[
\psi(F) \subseteq \psi(E-u_F). 
\]
In particular, $\psi(E-E) = \psi(U)$.
\end{theorem}

\begin{remark}
Theorem \ref{Thm_CombMain} is proved in Subsection \ref{subsec:FSthm}.
\end{remark}

We now state a simple corollary of Theorem~\ref{Thm_CombMain}. Given a set \( E \subset \mathbb{K}^2 \), define
\[
\Prod(E) = \{ xy \mid (x, y) \in E \} \subset \mathbb{K}.
\]
In particular, if \( E = A \times B \) for subsets \( A, B \subset \mathbb{K} \), then \( \Prod(E) = A \cdot B \). \\

A key theme in the sum-product phenomenon is that when \( E = A \times A \), for some set \( A \subset \mathbb{K} \), the set \( \Prod(E - E) \) is expected to be large, provided \( A \) is sufficiently dense. See, for instance, \cite[Theorem 1.4]{HIS} in the setting of finite fields, and \cite[Corollary 1.1]{F} for the case of integers. \\

To see how this set arises in our context, consider the \( \mathbb{K}^* \)-representation \( (\rho_U, U = \mathbb{K}^2) \) defined by
\begin{equation}
\label{rhoUintro}
\rho_U(a)(x, y) = (a x, a^{-1} y), \quad a \in \mathbb{K}^*, \; (x, y) \in U.
\end{equation}
The function \( \psi : \mathbb{K}^2 \to \mathbb{K} \) given by \( \psi(x, y) = x y \) is \( \rho_U \)-invariant and satisfies \( \Prod(E) = \psi(E) \). Theorem \ref{Thm_CombMain} then yields the following:

\begin{corollary}
\label{Cor_Products}
Let \(\mathbb{K}\) be a countably infinite field, and let \(E \subset \mathbb{K}^2\) be a set with positive upper Banach density. Then \(\Prod(E - E) = \mathbb{K}\). Moreover, for every finite set \(Q \subset \mathbb{K}\), there exists an element \(u_Q \in E\) such that \(Q \subset \Prod(E - u_Q)\).
\end{corollary}

Specializing further to product sets in \( \mathbb{K}^2 \), we deduce:

\begin{corollary}
Let \( \mathbb{K} \) be a countably infinite field, and let \( A, B \subset \mathbb{K} \) be sets with positive upper Banach density. Then, for every finite set \( Q \subset \mathbb{K} \), there exists \( (a_Q, b_Q) \in A \times B \) such that
\[
Q \subset (A - a_Q) \cdot (B - b_Q).
\]
In particular, \( \mathbb{K} = (A - A) \cdot (B - B) \).
\end{corollary}

This result may be viewed as an analogue of the sum-product estimates in \cite{HIS} and \cite{F}, now extended to arbitrary countably infinite fields. 

\begin{remark}[Hyperbolas do not contain infinite difference sets]
The proof of Theorem \ref{Thm_CombMain} for the representation 
\(\rho_U\) in \eqref{rhoUintro} (i.e., Corollary \ref{Cor_Products}) 
proceeds by showing that for every \(t \in \bK^*\), the \emph{hyperbola}
\[
\mathcal{H}_t = \{ (a, t/a) : a \in \bK^* \} \subset \bK^2
\]
is a set of measurable recurrence (i.e., a Poincar\'e set) in \((\bK^2,+)\).
It is well known that any set containing the difference set of an infinite
subset is automatically a set of measurable recurrence. Thus, the corollary 
would follow quickly if \(\mathcal{H}_t\) happened to contain an infinite 
difference set. However, as we show in Subsection \ref{subsec:hyperbola}, 
\(\mathcal{H}_t\) contains no difference set of size larger than four, at least 
when \(\Char(\bK) \neq 2\).
\end{remark}

\begin{remark}[Quotients versus products]
The case of fractions (as opposed to products) is significantly simpler and does not require any disjointness arguments. Consider the representation \( U = \mathbb{K}^2 \) with \( \mathbb{K}^* \)-action
\[
\rho_U(a)(x, y) = (a x, a y), \quad a \in \mathbb{K}^*, \; (x, y) \in \mathbb{K}^2,
\]
which consists of a single weight space. Define the map \( \psi : \mathbb{K}^2 \to \mathbb{K} \) by
\[
\psi(x, y) = 
\begin{cases}
\frac{x}{y}, & y \neq 0, \\
0, & y = 0,
\end{cases}
\]
which is clearly \( \rho_U \)-invariant. Then for any set \( E \subset \mathbb{K}^2 \) with positive upper Banach density, we have
\[
\psi(E - E) = \left\{ \frac{x}{y} \mid (x, y) \in E - E \right\} = \mathbb{K}.
\]
This conclusion does not require the full force of Theorem~\ref{Thm_CombMain}, and can instead be derived directly from Poincar\'e's recurrence theorem.
\end{remark}

Finally, we highlight that Theorem~\ref{Thm_CombMain} is flexible enough to handle more exotic versions of the sum-product phenomenon. As an illustration, consider the \( \mathbb{K}^* \)-representation $\rho$ on $\bK^3$ given by
\[
\rho(a)(x, y, z) = (a x, a^{-2} y, a^{-1} z), \quad a \in \mathbb{K}^*, \; (x, y, z) \in \mathbb{K}^3,
\]
and define the $\rho$-invariant map \( \psi : \mathbb{K}^3 \to \mathbb{K}^2 \) by
\[
\psi(x, y, z) = (x^3 y z, x z).
\]
A direct computation shows that
\[
\psi(\mathbb{K}^3) = \{ (s, t) \in \mathbb{K}^2 \mid t \neq 0 \} \cup \{(0, 0)\},
\]
so Theorem~\ref{Thm_CombMain} yields the following:

\begin{corollary}
Let \( A, B, C \subset \mathbb{K} \) be sets of positive upper Banach density. Then
\[
\left\{ (x^3 y z, x z) \mid x \in A - A,\; y \in B - B,\; z \in C - C \right\}
= \{ (s, t) \in \mathbb{K}^2 \mid t \neq 0 \} \cup \{(0, 0)\}.
\]
\end{corollary}

\subsection{Chromatic number of spacetime}

In a recent paper, Davies proves the following result \cite[Theorem~3]{Davies}, thereby resolving a problem on the chromatic number of space time by Kosheleva and Kreinovich \cite{KoKr}.

\begin{theorem}
For every partition $\bQ^2 = E_1 \sqcup E_2 \sqcup \cdots \sqcup E_r$, there
exist an index $j \in \{1,\dots,r\}$ and two points $(x_1,x_2), (y_1,y_2) \in E_j$ such that
\[
(x_1 - x_2)^2 - (y_1 - y_2)^2 = 1.
\]
\end{theorem}

Observe that in any such partition, at least one cell must have positive upper Banach density. Consequently, Corollary~\ref{Cor_Products} immediately yields the following density strengthening of Davies's theorem, valid over every countably infinite field of characteristic different from~$2$.

\begin{theorem}
Let $\bK$ be a countably infinite field with $\Char(\bK)\neq 2$, and let $E\subset\bK$ have positive upper Banach density. Then for every $z\in\bK$ there exist $(x_1,x_2),(y_1,y_2)\in E$ such that
\[
(x_1 - x_2)^2 - (y_1 - y_2)^2 = z.
\]
\end{theorem}

\begin{remark}
The restriction $\Char(\bK)\neq 2$ is essential. Indeed, let $\bK$ be a countably infinite field with $\Char(\bK)=2$, and let $E_o\subset\bK$ be a set of positive upper Banach density satisfying
\[
1\notin E_o+E_o.
\]
(We leave the construction of such sets to the reader).  
Define
\[
E := \{(x_1,x_2)\in\bK^2 : x_1+x_2\in E_o\}.
\]
Then $E$ has positive upper Banach density in $\bK^2$.  
Suppose, for contradiction, that the theorem above holds for the set $E$. Then there exist $(x_1,x_2),(y_1,y_2)\in E$ such that
\[
1=(x_1-x_2)^2-(y_1-y_2)^2.
\]
Since $\Char(\bK)=2$, we have $(x_1-x_2)^2=(x_1+x_2)^2$ for all $x_1,x_2$, and hence
\[
1=(x_1+x_2+y_1+y_2)^2.
\]
Set $u:=x_1+x_2$ and $v:=y_1+y_2$. By definition of $E$, we have $u,v\in E_o$. Moreover, in characteristic~$2$ we have $-1=1$, and therefore
\[
1 = x_1+x_2+y_1+y_2 = u+v \in E_o+E_o,
\]
contradicting the defining property of $E_o$. This shows that the assumption $\Char(\bK)\neq 2$ cannot be removed.
\end{remark}

\begin{proof}
When $\Char(\bK)\neq 2$, the linear map
\[
(x_1,x_2)\mapsto(x_1+x_2,\; x_1-x_2)
\]
is invertible and hence surjective. Therefore the image
\[
F := \{(x_1+x_2,\; x_1-x_2)\in\bK^2 : (x_1,x_2)\in E\}
\]
has positive upper Banach density in $\bK^2$. Fix $z\in\bK$. By Corollary~\ref{Cor_Products}, we may choose $(u_1,u_2),(v_1,v_2)\in F$ with
\[
(u_1 - v_1)(u_2 - v_2) = z.
\]
By the definition of $F$, there exist $(x_1,x_2),(y_1,y_2)\in E$ such that
\[
(u_1,u_2)=(x_1+x_2,\; x_1-x_2),\qquad
(v_1,v_2)=(y_1+y_2,\; y_1-y_2).
\]
A direct computation then gives
\begin{align*}
(u_1 - v_1)(u_2 - v_2)
&= (x_1 + x_2 - (y_1 + y_2))(x_1 - x_2 - (y_1 - y_2)) \\[0.2cm]
&= ((x_1-y_1)+(x_2-y_2))\,((x_1-y_1)-(x_2-y_2)) \\[0.2cm]
&= (x_1 - y_1)^2 - (x_2 - y_2)^2 = z,
\end{align*}
as required.
\end{proof}

\subsection{A Furstenberg-S\'ark\"ozy theorem for Laurent polynomials}

A classical and celebrated theorem---proved independently by Furstenberg and Sark\H{o}zy---states that if \(E \subset \mathbb{Z}\) has positive upper Banach density and \(p : \mathbb{Z} \to \mathbb{Z}\) is a polynomial with \(p(0) = 0\), then there exists a non-zero integer \(m\) such that \(p(m) \in E - E\). A minor variation of their argument yields the following extension to arbitrary countably infinite fields: if \(\mathbb{K}\) is such a field and \(E \subset \mathbb{K}\) has positive upper Banach density, then for every polynomial \(p : \mathbb{K} \to \mathbb{K}\) with \(p(0) = 0\), there exists \(a \in \mathbb{K}^\ast\) for which \(p(a) \in E - E\).

This naturally raises the question of whether the Furstenberg--Sark\H{o}zy theorem continues to hold when the polynomial \(p\) is replaced by a more general \emph{rational function}. The techniques developed in this paper allow us to answer this question affirmatively for the class of Laurent polynomials.

\begin{theorem}
\label{Thm_FSthm}
Let \( \mathbb{K} \) be a countably infinite field, and let \( E \subset \mathbb{K} \) be a set with positive upper Banach density. Then, for every Laurent polynomial of the form
\[
p(X) = \sum_{k=-N}^N p_k X^k \quad \text{with } p_0 = 0,
\]
there exists \( a \in \mathbb{K}^* \) such that \( p(a) \in E - E \).
\end{theorem}

\begin{remark}
The proof of this theorem is given in Subsection~\ref{subsec:FSthm}.
\end{remark}

\subsection{Acknowledgements}
The authors are grateful to SMRI, University of Sydney and Chalmers University in Gothenburg for their hospitality during our respective visits. M.B. was supported
by VR-grant 2023-03803. A.F. would like to thank the support of the Australian Research council via the grant DP240100472.

\setcounter{tocdepth}{1}
\tableofcontents

\section{Coupling non-resistance and linear relations}

\subsection{Pontryagin duality, measures and Fourier transforms}
\label{subsec:FT}
Let \( A \) be a discrete abelian group. Its \textbf{Pontryagin dual} is the compact metrizable group  
\[
\widehat{A} = \Hom(A, \mathbb{S}^1),
\]  
equipped with the topology of pointwise convergence. We denote its Haar probability measure by \( m_{\widehat{A}} \). Given $a \in A$, we define the character $a^{\vee}$ on $\widehat{A}$ by
\[
a^{\vee}(\xi) = \xi(a), \quad \xi \in \widehat{A}.
\]
Suppose \( \Gamma < \Aut(A) \) is an abelian subgroup. Then \( \Gamma \) acts on \( \widehat{A} \) via  
\[
\Gamma \times \widehat{A} \to \widehat{A}, \quad (\gamma, \xi) \mapsto \gamma.\xi,
\]  
where the action is defined by  
\[
(\gamma.\xi)(a) = \xi(\gamma(a)) \quad \text{for all } a \in A.
\]  
This defines a left action of \( \Gamma \) on \( \widehat{A} \) by continuous group automorphisms (if \( \Gamma \) is non-abelian, this instead defines a right action), preserving the Haar probability measure \( m_{\widehat{A}} \). Moreover, the $\Gamma$-action is \textbf{ergodic} with respect to \( m_{\widehat{A}} \) provided that the \( \Gamma \)-orbit of every nonzero element of \( A \) is infinite. Given any probability measure \( \mu \) on \( \smash{\widehat{A}} \), we define its \textbf{Fourier transform} by  
\[
\widehat{\mu}(a) = \int_{\widehat{A}} \xi(a) \, d\mu(\xi), \quad a \in A.
\]  
The Fourier transform \( \widehat{\mu} \) uniquely determines the measure \( \mu \). In particular, we have  
\[
\widehat{\delta}_1(a) = 1 \quad \text{for all } a \in A,
\qand
\widehat{m}_{\widehat{A}}(a) = 
\begin{cases}
1 & \text{if } a = 0, \\
0 & \text{if } a \neq 0
\end{cases}.
\]
If $\varphi \in L^1(\widehat{A},m_{\widehat{A}})$, we define the Fourier transform by
\[
\widehat{\varphi}(a) = \int_{\widehat{A}} \varphi(\xi) \xi(a) \, dm_{\widehat{A}}(\xi), \quad a \in A.
\]
If $\varphi$ is continuous on $\widehat{A}$ and $\widehat{\varphi} \in \ell^1(A)$, we can write
\[
\varphi(\xi) = \sum_{a \in A} \widehat{\varphi}(a) (-a)^{\vee}(\xi), \quad \textrm{for all $\xi \in \widehat{A}$}.
\]

\subsection{Couplings}

Let \( (X, \mathscr{B}_X) \) be a measurable space. We denote by \( \Prob(X) \) the space of probability measures on \( \mathscr{B}_X \), suppressing the dependence on the sigma-algebra when no confusion arises. \\

Let \( \Gamma \) be a countable discrete group, and suppose \( a : \Gamma \times X \to X \) is a measurable action of \( \Gamma \) on \( X \). We typically use the shorthand notation \( \gamma.x := a(\gamma, x) \) to denote the action, and note that it extends to $\Prob(X)$ by defining
\[
(\gamma.\mu)(A) = \mu(\gamma^{-1}.A), \quad \gamma \in \Gamma, \enskip \mu \in \Prob(X), \enskip A \in \mathscr{B}_X.
\]
The set of \( \Gamma \)-invariant probability measures on \( X \) is denoted by \( \Prob_\Gamma(X) \), which may be empty in general. When \( \mu \in \Prob_\Gamma(X) \), we refer to the measure space \( (X, \mu) \) as a \textbf{Borel \( \Gamma \)-space}, again suppressing explicit reference to the underlying sigma-algebra. \\

Given two Borel \( \Gamma \)-spaces \( (X, \mu) \) and \( (Y, \nu) \), a \textbf{coupling} of \( \mu \) and \( \nu \) is a probability measure \( \eta \) on the product sigma-algebra \( \mathscr{B}_X \otimes \mathscr{B}_Y \) satisfying
\[
\eta(A \times Y) = \mu(A) \quad \text{and} \quad \eta(X \times B) = \nu(B), \quad \text{for all } A \in \mathscr{B}_X,\, B \in \mathscr{B}_Y.
\]
We refer to \( \mu \) and \( \nu \) as the \textbf{marginals} of \( \eta \). A canonical example of a coupling is the \textbf{trivial coupling} \( \mu \otimes \nu \), the product measure on \( X \times Y \). \\

We denote by \( \mathscr{C}(\mu, \nu) \) the space of all couplings of \( \mu \) and \( \nu \), and by \( \mathscr{C}_\Gamma(\mu, \nu) \) the subset consisting of those couplings that are invariant under the diagonal \( \Gamma \)-action on \( X \times Y \). Clearly, if both \( \mu \) and \( \nu \) are \( \Gamma \)-invariant, then the trivial coupling \( \mu \otimes \nu \) lies in \( \mathscr{C}_\Gamma(\mu, \nu) \). \\

More generally, if we are given Borel \( \Gamma \)-spaces \( (X_1, \mu_1), \ldots, (X_r, \mu_r) \), we write \( \mathscr{C}_\Gamma(\mu_1, \ldots, \mu_r) \) for the space of all \( \Gamma \)-invariant probability measures on the product $X_1 \times \cdots \times X_r$ whose marginals on each \( X_i \) is \( \mu_i \), for \( i = 1, \ldots, r \).

\subsection{Markov intertwiners}

Let \( \Gamma \) be a countable discrete group, and let \( (X, \mu) \) and \( (Y, \nu) \) be Borel \( \Gamma \)-spaces. Suppose \( \eta \in \mathscr{C}_\Gamma(\mu, \nu) \) is a \( \Gamma \)-invariant coupling of \( \mu \) and \( \nu \). Then, by disintegrating \( \eta \) over \( (X, \mu) \) and \( (Y, \nu) \), respectively, we obtain measurable maps
\[
x \mapsto \eta_x \in \Prob(Y), \quad y \mapsto \eta^y \in \Prob(X),
\]
defined for \( \mu \)-almost every \( x \in X \) and \( \nu \)-almost every \( y \in Y \), such that
\[
\eta = \int_X \delta_x \otimes \eta_x \, d\mu(x) = \int_Y \eta^y \otimes \delta_y \, d\nu(y).
\]
These disintegrations are \( \Gamma \)-equivariant and give rise to associated \textbf{Markov intertwiners}
\[
S : L^2(X, \mu) \to L^2(Y, \nu), \quad T : L^2(Y, \nu) \to L^2(X, \mu),
\]
defined by
\[
Sf(y) = \int_X f(x) \, d\eta^y(x), \quad Tg(x) = \int_Y g(y) \, d\eta_x(y),
\]
for all \( f \in L^2(X, \mu) \) and \( g \in L^2(Y, \nu) \). These linear maps are bounded, \( \Gamma \)-equivariant, and satisfy the identity
\[
\eta(f \otimes g) = \langle Sf, g \rangle_{L^2(Y, \nu)} = \langle f, Tg \rangle_{L^2(X, \mu)},
\]
for all \( f \in L^2(X, \mu) \), \( g \in L^2(Y, \nu) \). Moreover, both \( S \) and \( T \) preserve non-negativity and the constant function: \( S1 = 1 \) and \( T1 = 1 \). Moreover, $S$ and $T$ clearly maps $L^\infty$-spaces into $L^\infty$-spaces. \\

The following lemma provides a simple characterization of the trivial coupling. We omit the proof.

\begin{lemma}
\label{Lemma_CouplingMarkov}
The coupling \( \eta \) is trivial (i.e., \( \eta = \mu \otimes \nu \)) if and only if either
\[
Sf = \int_X f \, d\mu \quad \text{for all } f \in L^2(X, \mu),
\quad \text{or} \quad
Tg = \int_Y g \, d\nu \quad \text{for all } g \in L^2(Y, \nu),
\]
where the equalities are understood in the \( L^2 \)-sense.
\end{lemma}

\subsection{Measure-rigid groups}
\label{subsec:measure-rigid}

Let $\bK$ be a countably infinite discrete field and let $\bK^*$ denote its multiplicative group. Note that the $\bK^*$ acts on $(\bK,+)$ by automorphism via standard multiplication. We denote by $\widehat{\bK}$ the Pontryagin dual of $(\bK,+)$.

\begin{definition}
A subgroup \( \Delta < \bK^* \) is called \textbf{measure-rigid} if the only \(\Delta\)-invariant and \(\Delta\)-ergodic probability measures on \( \widehat{\bK} \) are \( \delta_1 \) and \( m_{\widehat{\bK}} \).
\end{definition}

\begin{remark}
We emphasize that measure-rigid subgroups do not exist when \(\bK\) is a finite field. Indeed, in this case, the dual group satisfies \(\smash{\widehat{\bK}} \cong \bK\), implying that \(\widehat{\bK}\) consists of exactly two \(\bK^*\)-orbits: the fixed point \(1\) and its complement. The normalized counting measure on the complement of $\{1\}$ is clearly \(\bK^*\)-invariant and \(\bK^*\)-ergodic, yet it clearly differs from the Haar measure \(m_{\widehat{\bK}}\), demonstrating the absence of measure rigidity.
\end{remark}

Our main examples of measure-rigid subgroups are constructed from the following two lemmas.

\begin{lemma}
\label{Lemma_FiniteIndexMeasureRigid}
Every finite-index subgroup $\Delta < \bK^*$ is measure-rigid. 
\end{lemma}

\begin{proof}
Let $\mu$ be a $\Delta$-invariant and $\Delta$-ergodic probability measure on $\widehat{\bK}$. First, suppose there exists some $\xi_o \in \widehat{\bK}$ such that $\mu(\{\xi_o\}) > 0$. Then, by $\Delta$-invariance, we have 
\[
\mu(\{\delta^*.\xi_o\}) = \mu(\{\xi_o\}) > 0, \quad \text{for all } \delta \in \Delta.
\]
This immediately implies that the orbit $\Delta^*.\xi_o \subset \widehat{\bK}$ is finite. However, since $\bK^*$ acts freely on $\widehat{\bK} \setminus \{1\}$ and $\Delta$ has finite index in the infinite group $\bK^*$, we conclude that $\xi_o = 1$. Since $\delta_1$ is clearly $\Delta$-invariant and $\Delta$-ergodic, we may assume from now on that $\mu$ is non-atomic. Our goal is to show that $\mu = m_{\widehat{\bK}}$. Since $\mu$ is non-atomic, Wiener's Lemma gives us the asymptotic formula
\[
\lim_{n \to \infty} \frac{1}{|F_n|} \sum_{b \in F_n} |\widehat{\mu}(b)|^2 = 0
\]
for any F\o lner sequence $(F_n)$ in $(\bK,+)$. Now, using the decomposition
\begin{equation}
\label{fidecomp}
\bK^* = \bigsqcup_{q \in Q} \Delta q, \quad \text{for some finite set } Q \subset \bK^*,
\end{equation}
we obtain
\[
\lim_{n \to \infty} \frac{1}{|F_n|} \sum_{b \in F_n} |\widehat{\mu}(b)|^2
\, = \, \lim_{n \to \infty} \sum_{q \in Q} \frac{|F_nq^{-1} \cap \Delta|}{|F_n|} |\widehat{\mu}(q)|^2.
\]
Choosing $(F_n)$ to be a double F\o lner sequence in $\bK$, we observe that $(F_n q^{-1} \setminus \{0\})$ forms a F\o lner sequence in $(\bK^*,\cdot)$ for every $q \in Q$, leading to
\[
\lim_{n \to \infty} \frac{|F_nq^{-1} \cap \Delta|}{|F_n|} = \frac{1}{|Q|}, \quad \text{for all } q \in Q.
\]
Thus, we conclude that
\[
\lim_{n \to \infty} \frac{1}{|F_n|} \sum_{b \in F_n} |\widehat{\mu}(b)|^2  = \frac{1}{|Q|} \sum_{q \in Q} |\widehat{\mu}(q)|^2 = 0.
\]
This forces $\widehat{\mu}(q) = 0$ for all $q \in Q$. Since $\mu$ is $\Delta$-invariant, it follows from \eqref{fidecomp} that $\widehat{\mu}(b) = 0$ for all $b \neq 0$, implying that $\mu = m_{\widehat{\bK}}$. This completes the proof.
\end{proof}

\begin{lemma}
\label{Lemma_KnIsMeasureRigid}
Let \( n \) be a positive integer that is not divisible by the characteristic of $\bK$.
Then, for every F\o lner sequence $(F_k)$ in $(\bK,+)$ and $\xi \in \widehat{\bK} \setminus \{1\}$, 
\[
\lim_{k \ra \infty} \frac{1}{|F_k|} \sum_{a \in F_k} \xi(a^n) = 0.
\]
In particular, the subgroup \( \Delta_n = \{ a^n \mid a \in \bK^* \} \) is measure-rigid.
\end{lemma}

\begin{remark}
Observe that $\Delta_{-n} = \Delta_n$, ensuring that the measure-rigidity aspect of the lemma remains valid for negative $n$. However, as we will see in Section~\ref{Sec:Failure}, there exist countably infinite fields where the equidistribution conclusion fails when $n$ is negative.
\end{remark}

\begin{proof}
When $\bK$ has zero characteristic, the proof follows from van der Corput's lemma. The result in full generality is established in \cite[Theorem 2.2.1, p. 60]{L}.
\end{proof}

\begin{remark}
If \(\Char(\bK) = 0\), then the condition on the integer \( n \) is automatically satisfied. However, when \(\Char(\bK) = p > 0\), the requirement that \( n \) is not divisible by \( p \) becomes essential. For instance, if \( n = p \) and \(\bL := \{ x^p \mid x \in \bK \}\) is a proper (additive) subgroup of $\bK$, then $
\smash{m_{\bL^{\perp}}}$
is a \(\smash{\Delta}\)-invariant and \(\smash{\Delta}\)-ergodic probability measure on 
$\smash{\widehat{\bK}}$ that is distinct from both \(\smash \delta_1\) and \( \smash m_{\widehat{\bK}} \).
\end{remark}

Let us recall some notation: If $V$ is a vector space over $\bK$, we denote by $\mathscr{L}(V)$ the set of linear subspaces of $V$. If $V$ is finite-dimensional over $\bK$, then $\mathscr{L}(V)$ is countable. 

\begin{lemma}
\label{Lemma_MeasureRigidEquivalence}
Let $\Delta < \bK^*$ be a subgroup. Then the following conditions are equivalent:
\begin{itemize}
\item[$(i)$] $\Delta$ is measure-rigid. \vspace{0.2cm}
\item[$(ii)$] For every countable $\bK$-vector space $V$, the only $\Delta$-invariant and $\Delta$-ergodic probability measures on $\smash{\widehat{V}}$ are of the form $m_{W^{\perp}}$
for some linear subspace $W \leq V$. \vspace{0.2cm}
\item[$(iii)$] For every finite-dimensional $\bK$-vector space $V$, the map 
\[
\Prob(\mathscr{L}(V)) \ra \Prob_\Delta(\widehat{V}), \quad p \mapsto \mu_p = \sum_{W \leq V} p(W) m_{W^{\perp}}
\]
is an affine bijection. Furthermore, $\mu_p$ is $\Delta$-ergodic if and only if $p = \delta_{\{W\}}$ for some linear subspace $W \leq V$. 
\end{itemize}
\end{lemma}

\begin{remark}
The assumption in (ii) that the vector space $V$ is countable is not strictly necessary but suffices for our purposes. If $V$ is uncountable (yet discrete), its dual $\smash{\widehat{V}}$ becomes non-metrizable, introducing technical complications in measure analysis that would divert from the main focus of the paper. Likewise, the assumption in (iii) that $V$ is finite-dimensional can be relaxed, but doing so would introduce additional difficulties that are not essential to our discussion.
\end{remark}

\begin{proof}
Since the only linear subspaces of the $\bK$-vector space $\bK$ are $\{0\}$ and $\bK$, condition (iii) clearly implies (i). It remains to show that (i) implies (ii) and (ii) implies (iii). \\

To prove that (i) implies (ii), assume $\Delta$ is measure-rigid. Fix a countable $\bK$-vector space $V$ and a $\Delta$-invariant, $\Delta$-ergodic probability measure $\mu$ on $\smash{\widehat{V}}$. Consider the pairing  
\[
\widehat{V} \times V \to \widehat{\bK}, \quad (\xi,v) \mapsto \xi_v, \quad \text{where } \xi_v(b) = \xi(bv) \text{ for } b \in \bK.
\]
Since $(a^*\xi)_v = a^*\xi_v$ for all $a \in \bK^*$, the push-forward $\mu_v$ of $\mu$ under $\xi \mapsto \xi_v$ is a $\Delta$-invariant, $\Delta$-ergodic probability measure on $\smash{\widehat{\bK}}$, and thus must be either $\delta_1$ or $m_{\widehat{\bK}}$. Define  
\[
W_\mu = \{ v \in V \mid \mu_v = \delta_1 \}.
\]
Clearly, $0 \in W_\mu$, and for $v \notin W_\mu$, we have $\mu_v = m_{\widehat{\bK}}$. To show that $W_\mu$ is a linear subspace, note that for $v \in W_\mu$,  
\[
\widehat{\mu}_v(b) = \int_{\widehat{V}} \xi(bv) \, d\mu(\xi) = 1 \quad \text{for all } b \in \bK.
\]
Thus, the set  
\[
X_v = \{ \xi \in \widehat{V} \mid \xi(bv) = 1 \text{ for all } b \in \bK \}
\]
is $\mu$-conull in $\widehat{V}$. If $v_1, v_2 \in W_\mu$, then for all $b_1, b_2, b \in \bK$,  
\[
\xi(bb_1v_1) = \xi(bb_2v_2) = 1 \quad \text{for $\mu$-almost every } \xi.
\]
Hence, $\xi(b(b_1v_1 + b_2v_2)) = 1$ for $\mu$-almost every $\xi$, implying $\widehat{\mu}_{b_1v_1 + b_2 v_2}(b) = 1$ for all $b \in \bK$, so $b_1 v_1 + b_2 v_2 \in W_\mu$. Since this holds for all choices of $b_1, b_2, v_1, v_2$, we conclude that $W_\mu$ is a subspace. Since $\mu_v = m_{\widehat{\bK}}$ for $v \notin W_\mu$, we get $\smash{\widehat{\mu}_v} = \chi_{W_\mu}$, and since Fourier transforms uniquely determine measures, it follows that $\mu = m_{W^\perp_\mu}$.

Conversely, for any subspace $W \leq V$, the measure $m_{W^\perp}$ is $\Delta$-invariant and $\Delta$-mixing, hence $\Delta$-ergodic, proving that (i) implies (ii). \\

To show that (ii) implies (iii), consider a finite-dimensional $\bK$-vector space $V$ and a $\Delta$-invariant probability measure $\mu$ on $\smash{\widehat{V}}$. By (ii), every ergodic component of $\mu$ is of the form $m_{W^\perp}$ for some $W \leq V$. Since $V$ has only countably many subspaces, the ergodic decomposition of $\mu$ must match the form of $\mu_p$ for some unique probability measure $p$ on $\mathscr{L}(V)$.
\end{proof}

\subsection{Coupling resistance}

Let $\bK_1,\ldots,\bK_r$ be countably infinite fields. 

\begin{definition}
A subgroup $\Delta < \bK_1^* \times \cdots \bK_r^*$ is said to be \textbf{coupling-resistant}
if 
\[
\mathscr{C}_\Delta
\big(
(\widehat{\bK}_1,m_{\widehat{\bK}_1}),\ldots,(\widehat{\bK}_r,m_{\widehat{\bK}_r})
\big)
= \{ m_{\widehat{\bK}_1} \otimes \cdots \otimes m_{\widehat{\bK}_r} \}.
\]
This is, the product measure is the only $\Delta$-invariant coupling of 
$(\widehat{\bK}_1,m_{\widehat{\bK}_1}),\ldots,(\widehat{\bK}_r,m_{\widehat{\bK}_r})$.
\end{definition}

\begin{remark}
\label{Rmk_CouplingInherits}
If \( \Delta < \bK_1^* \times \cdots \times \bK_r^* \) is coupling-resistant, then each projection \( \Delta_I < \prod_{i \in I} \bK_i^* \) for any subset \( I \subseteq \{1, \ldots, r\} \) is also coupling-resistant. Moreover, if $r=1$, then any subgroup of $\bK^*$ is trivially coupling-resistant. 
\end{remark}

The following lemma extends Lemma \ref{Lemma_MeasureRigidEquivalence}. 

\begin{lemma}
\label{Lemma_CouplingResistantEquivalence}
Let $\Delta < \bK_1^* \times \cdots \bK_r^*$ be a subgroup such that $\Delta_i = \pr_i(\Delta) < \bK_i^*$ is measure-rigid for every $i=1,\ldots,r$. Then the 
following three conditions are equivalent:
\begin{itemize}
\item[$(i)$] $\Delta$ is coupling-resistant. \vspace{0.2cm}
\item[$(ii)$] For any countable vector spaces \( V_1, \dots, V_r \) over \( \bK_1, \dots, \bK_r \) respectively, every $\Delta$-invariant and $\Delta$-ergodic probability measure $\theta$ on $\widehat{V}_1 \times \cdots \widehat{V}_r$ is of the form
\[
\theta = m_{W_1^{\perp}} \otimes \cdots \otimes m_{W_r^{\perp}},
\]
for some linear sub-spaces $W_1 \leq V_1,\ldots,W_r \leq V_r$. \vspace{0.1cm}
\item[$(iii)$] For any finite-dimensional vector spaces \( V_1, \dots, V_r \) over \( \bK_1, \dots, \bK_r \), respectively, the map
\[
\Prob(\mathscr{L}(V_1) \times \cdots \mathscr{L}(V_r)) \ra \Prob_\Delta(\widehat{V}_1 \times \cdots \times \widehat{V}_r), \quad p \mapsto \theta_p
\]
where 
\[
\theta_p = \sum_{W_1 \leq V_1} \cdots \sum_{W_r \leq V_r} p(W_1,\ldots,W_r)  m_{W_1^{\perp}} \otimes \cdots \otimes m_{W_r^{\perp}},
\]
is an affine bijection. Furthermore, $\theta_p$ is $\Delta$-ergodic if and only if 
$p = \delta_{\{(W_1,\ldots,W_r)\}}$ for some linear subspaces $W_1 \leq V_1,\ldots,W_r \leq V_r$.  
\end{itemize}
\end{lemma}

\begin{proof}
Clearly, both (ii) and (iii) imply (i), so it suffices to show that (i) implies (ii) and that (ii) implies (iii). \\

We first prove that (i) implies (ii). Suppose that \( \Delta \) is coupling-resistant, and fix countable vector spaces \( V_1, \dots, V_r \) over \( \bK_1, \dots, \bK_r \), respectively. Let \( \theta \) be a \( \Delta \)-invariant and \( \Delta \)-ergodic probability measure on \( \widehat{V}_1 \times \cdots \times \widehat{V}_r \). \\

Following the proof of Lemma \ref{Lemma_MeasureRigidEquivalence}, consider the pairings
\[
\widehat{V}_i \times V_i \to \widehat{\bK}_i, \quad (\xi,v) \mapsto \xi_v, \quad \text{where } \xi_v(b) = \xi(bv) \text{ for } b \in \bK_i, \enskip i = 1,\dots,r.
\]
For a fixed \( \underline{v} = (v_1, \dots, v_r) \in V_1 \times \cdots \times V_r \), let \( \theta_{\underline{v}} \) be the pushforward of \( \theta \) under the map
\[
\widehat{V}_1 \times \cdots \times \widehat{V}_r \to \widehat{\bK}_1 \times \cdots \times \widehat{\bK}_r, 
\quad (\xi_1, \dots, \xi_r) \mapsto ((\xi_1)_{v_1}, \dots, (\xi_r)_{v_r}).
\]
Then, for all \( (b_1, \dots, b_r) \in \bK_1 \times \cdots \times \bK_r \),
\[
\widehat{\theta}(b_1v_1, \dots, b_r v_r)
= \int_{\widehat{\bK}_1} \cdots \int_{\widehat{\bK}_r} 
\eta_1(b_1) \cdots \eta_r(b_r) \, d\theta_{\underline{v}}(\eta_1, \dots, \eta_r).
\]
Since the map above is \( \bK_1^* \times \cdots \times \bK_r^* \)-equivariant, it follows that \( \theta_{\underline{v}} \) is \( \Delta \)-invariant and \( \Delta \)-ergodic. Consequently, the pushforward \( (\theta_{\underline{v}})_i \) of \( \theta_{\underline{v}} \) to each coordinate \( i \) must be either \( \delta_1 \) or \( m_{\widehat{\bK}_i} \). Define
\[
I_\theta(\underline{v}) := \{ i \in \{1,\dots,r\} \mid (\theta_{\underline{v}})_i = \delta_1 \}.
\]
If \( I_\theta(\underline{v}) \) is empty, then \( \theta_{\underline{v}} \) is a \( \Delta \)-invariant coupling of \( (\widehat{\bK}_1, m_{\widehat{\bK}_1}), \dots, (\widehat{\bK}_r, m_{\widehat{\bK}_r}) \), which, by the coupling-resistance of \( \Delta \), must be the product measure \( m_{\widehat{\bK}_1} \otimes \cdots \otimes m_{\widehat{\bK}_r} \). If \( I_\theta(\underline{v}) \) is nonempty, then \( \theta_{\underline{v}} \) decomposes as a product of delta measures at the trivial characters for \( i \in I_\theta(\underline{v}) \) and a \( \Delta \)-invariant coupling of \( \smash{(\widehat{\bK}_i, m_{\widehat{\bK}_i})} \) over \( i \notin I_\theta(\underline{v}) \), where \( \Delta \) acts via the projection onto the subspace indexed by \( i \notin I_\theta(\underline{v}) \). Again, by the coupling-resistance of \( \Delta \), this coupling must be the product of Haar measures. We conclude that
\[
\widehat{\theta}(b_1 v_1, \dots, b_r v_r)
= 
\begin{cases}
1, & \text{if } b_i \in \bK_i \text{ for } i \in I_\theta(\underline{v}) \text{ and } b_i = 0 \text{ for } i \notin I_\theta(\underline{v}), \\[0.2cm]
0, & \text{otherwise}.
\end{cases}
\]
Thus, \( \widehat{\theta} \) is the indicator function of the set
\[
U := \{ \underline{v} \in V_1 \times \cdots \times V_r \mid \widehat{\theta}(\underline{v}) = 1 \}.
\]
Since \( \widehat{\theta} \) is invariant under multiplication by \( (b_1, \dots, b_r) \in \bK_1 \times \cdots \times \bK_r \), and using the same argument as in Lemma \ref{Lemma_MeasureRigidEquivalence}, we conclude that \( U \) is also closed under addition. Consequently, \( U \) splits as a product of subspaces \( W_1, \dots, W_r \), where \( W_i < V_i \), so that
\[
\theta = m_{W_1^{\perp}} \otimes \cdots \otimes m_{W_r^\perp}.
\]
Conversely, every measure of this form is \( \Delta \)-invariant and \( \Delta \)-ergodic. This proves that (i) implies (ii). \\

Next, we prove that (ii) implies (iii). Suppose each vector space \( V_i \) is finite-dimensional over the field \( \bK_i \) for \( i=1,\dots,r \). Since there are only countably many linear subspaces of \( V_1, \dots, V_r \), the unique ergodic decomposition of any \( \Delta \)-invariant probability measure on \( \smash{\widehat{V}_1 \times \cdots \times \widehat{V}_r} \) is a sum of measures as in (ii), completing the proof.
\end{proof}

\subsection{Absence of higher-order decorrelation implies linear relations}

In what follows, let \(\Delta < \bK_1^* \times \cdots \times \bK_r^*\) be a coupling-resistant subgroup such that each projection \(\Delta_i = \pr_i(\Delta) < \bK_i^*\) is measure-rigid for \(i = 1, \dots, r\). We write
\[
A = \widehat{\bK}_1 \times \cdots \times \widehat{\bK}_r \qand m_A = m_{\widehat{\bK}_1} \otimes \cdots \otimes m_{\widehat{\bK}_r}.
\]
Let \(N \geq 2\) be an integer, and pick a $\Delta$-invariant $N$-coupling $\theta$ of
$(A,m_A)$. \\

For \(\alpha \geq 0\) and an \(N\)-tuple \(\underline{f} = (f_1, \dots, f_N) \in L^\infty(A)^N\), define the (possibly empty) set
\[
S_\alpha(\theta, \underline{f}) = \Big\{ (\delta_1, \dots, \delta_N) \in \Delta^N \, : \, 
\big|\theta(f_1 \circ \delta_1 \otimes \cdots \otimes f_N \circ \delta_N) \big| \geq \alpha \Big\}.
\]
If \( m_{A}(f_i) = 0 \) for all \( i = 1, \dots, N \), then the size of the set \( S_\alpha(\theta, \underline{f}) \) quantifies the failure of higher-order decorrelation of \( \theta \) under the action of $\Delta$ on $(A,m_A)$.\\

Our next theorem establishes that for \( \alpha > 0 \), the elements of \( S_\alpha(\theta, \underline{f}) \) must satisfy a finite set of linear relations, implying that this set is not excessively curved when viewed within the  \emph{additive} group \( \bK^N_1 \times \cdots \times \bK^N_r \cong \widehat{A}^N \). \\

To state the theorem, we introduce some notation. Given an $N$-tuple 
\[
(u_1,\ldots,u_N) \in (\bK_1 \times \cdots \bK_r)^N,
\]
where $u_i = (u_{i,1},\ldots,u_{i,r})$, with $u_{i,j} \in \bK_j$ for every $i=1,\ldots,N$
and  $j=1,\ldots,r$, we define its \textbf{transpose} $(v_1,\ldots,v_r) \in \bK_1^N \times \cdots \bK_r^N$ by
\[
(v_1,\ldots,v_r) = (u_1,\ldots,u_N)^{\mathrm{tr}}, \quad v_j = (u_{1,j},\ldots,u_{N,j}) \in \bK_j^N.
\]
One readily checks that the transpose operation is bijective and respects both addition and coordinate-wise multiplication. We now state our theorem.

\begin{theorem}[Forcing linear relations]
\label{Thm_F}
Let \( \Delta \) and \( \theta \) be as above. Then, for every \( \alpha > 0 \) and every \( N \)-tuple \( \underline{f} = (f_1, \dots, f_N) \) satisfying \( m_A(f_i) = 0 \) for all \( i \), there exists a finite set  
\[
\cF = \cF(\theta, \underline{f}, \alpha) \subset \smash{\mathscr{L}(\bK_1^N) \times \cdots \times \mathscr{L}(\bK_r^N)}
\]  
such that for every \( (W_1,\ldots,W_r) \in \cF \), at least one index \( i \in \{1,\ldots,r\} \) satisfies \( W_i \neq \bK_i^N \), and  
\[
\forall \, (\delta_1, \dots, \delta_N) \in S_\alpha(\theta,\underline{f}), \quad \exists \, (W_1, \dots, W_r) \in \cF \quad \text{such that} \quad \lambda_j \in W_j \text{ for all } j = 1, \dots, r,
\]  
where \( (\lambda_1, \dots, \lambda_r) = (\delta_1, \dots, \delta_N)^{\mathrm{tr}} \). 
\end{theorem}

\subsection{Proof of Theorem \ref{Thm_F}}

Let $\Delta < A$ and $\theta$ be as above. 

\begin{definition}
The \textbf{transpose measure} of $\theta$ is the probability measure 
\[
\tau_\theta \in \Prob(\widehat{\bK}_1^N \times \cdots \times \widehat{\bK}_r^N)
\] 
defined by
\[
\widehat{\tau}_\theta(v_1,\dots,v_r) = \widehat{\theta}(u_1,\dots,u_N),
\]
where $(v_1,\dots,v_r) = (u_1,\dots,u_N)^{\mathrm{tr}} \in \bK_1^N \times \cdots \times \bK_r^N$.
\end{definition}

\begin{remark}
Since the transpose operation is bijective and preserves addition, $\widehat{\tau}_\theta$ is a well-defined positive definite function on $\bK_1^N \times \cdots \times \bK_r^N$. By Bochner's theorem, this function uniquely determines a probability measure on $\smash{\widehat{\bK}_1^N \times \cdots \times \widehat{\bK}_r^N}$.
\end{remark}

Since $\theta$ is $\Delta$-invariant, we have
\[
\widehat{\theta}(\delta u_1,\dots, \delta u_N) = \widehat{\theta}(u_1,\dots,u_N)
\]
for all $(u_1,\dots,u_N) \in A^N$ and $\delta = (\delta^{(1)},\dots,\delta^{(r)}) \in \Delta$, where
\[
\delta u_j = (\delta^{(1)} u_{j,1},\dots,\delta^{(r)} u_{j,r}), \quad j = 1,\dots,N.
\]
Since $\widehat{\tau}_\theta(\delta^{(1)} v_1,\dots,\delta^{(r)} v_r) = \widehat{\theta}(\delta u_1,\dots, \delta u_N)$, it follows that $\tau_\theta$ is also $\Delta$-invariant. \\

For $i = 1,\dots,r$ and $j = 1,\dots,N$, define the projection map
\[
\pi_{i,j} : \widehat{\bK}_1^N \times \cdots \times \widehat{\bK}_r^N \to \widehat{\bK}_i, \quad (\xi_1,\dots,\xi_r) \mapsto \xi_{i,j},
\]
where each $\xi_i \in \widehat{\bK}_i^N$ satisfies
\[
\xi_{i}(v_i) = \xi_{i,1}(v_{i,1}) \cdots \xi_{i,N}(v_{i,N}), \quad v_i = (v_{i,1},\dots,v_{i,N}) \in \bK_i^N.
\]
Since $\theta$ is an $N$-coupling of $(A, m_A)$, it follows that
\[
(\pi_{i,j})_* \tau_\theta = m_{\bK_i}, \quad \forall i,j.
\]
By Lemma \ref{Lemma_CouplingResistantEquivalence}, applied to the vector spaces $V_i = \bK_i^N$ for $i=1,\ldots,r$, there exists a unique probability measure $p_\theta$ on $\mathscr{L}(V_1) \times \cdots \times \mathscr{L}(V_r)$ such that
\[
\tau_\theta = \sum_{U_1 \leq \bK_1^N} \cdots \sum_{U_r \leq \bK_r^N} p_\theta(U_1,\dots,U_r) m_{U_1^{\perp}} \otimes \cdots \otimes m_{U_r^{\perp}}.
\]

\begin{lemma}
\label{Lemma_c}
For every $r$-tuple $(U_1,\ldots,U_r) \in \supp(p_\theta)$ and $i=1\,\ldots,r$, there 
exists an $N$-tuple $c_i = (c_{i,1},\ldots,c_{i,N}) \in \smash{(\bK_i^*)^N}$, depending on $U_i$, such that
\[
\sum_{j=1}^N c_{i,j} v_{i,j} = 0, \quad \textrm{for all $v_i = (v_{i,1},\ldots,v_{i,N}) \in U_i$}.
\]
\end{lemma}

\begin{proof}
We show that, after possibly permuting coordinates, we cannot write
\[
U_i = \bK_i \times U_i' \quad \text{for some linear subspace } U_i' < \bK_i^{N-1}.
\]  
Otherwise, 
\[
U_i^{\perp} = \{1\} \times (U_i')^{\perp}.
\]  
As a result, \(\pi_{i,1}(U_i^{\perp}) = \{1\}\), contradicting our earlier observation that  
\[
(\pi_{i,j})_* \tau_\theta = m_{\widehat{\bK}_i} \quad \text{for all } i \text{ and } j.
\]  
\end{proof}

Let \((\delta_1, \dots, \delta_N) \in \Delta^N\), where each \(\delta_j = (\delta_{j,1}, \dots, \delta_{j,r})\). Then, for all 
\[
(u_1, \dots, u_N) \in (\bK_1 \times \cdots \times \bK_r)^N,
\]
we have  
\[
(\delta_1 u_1, \dots, \delta_N u_N)^{\mathrm{tr}}
= (\lambda_1 v_1, \dots, \lambda_r v_r),
\]
where  
\[
(v_1, \dots, v_r) = (u_1, \dots, u_N)^{\mathrm{tr}},  
\quad \text{and} \quad  
(\lambda_1, \dots, \lambda_r) = (\delta_1, \dots, \delta_N)^{\mathrm{tr}}.
\]
Moreover, the components satisfy  
\[
(\lambda_i v_i)_j = \delta_{j,i} u_{j,i},  
\quad \text{for } i = 1, \dots, r, \quad j = 1, \dots, N.
\]
In particular, since \(\widehat{m}_{U_i^{\perp}}(\lambda_i v_i) = \chi_{U_i}(\lambda_i v_i)\) for all \(i\), we obtain  
\begin{align*}
\widehat{\theta}(\delta_1 u_1, \dots, \delta_N u_N)
&= \widehat{\tau}_\theta(\lambda_1 v_1, \dots, \lambda_r v_r) \\[0.2cm]
&= \sum_{U_1 \leq \bK_1^N} \cdots \sum_{U_r \leq \bK_r^N}
p(U_1, \dots, U_r) \chi_{U_1}(\lambda_1 v_1) \cdots \chi_{U_r}(\lambda_r v_r).
\end{align*}

\begin{lemma}
\label{Lemma_W}
For every \(\cU = (U_1,\dots,U_r) \in \supp(p_\theta)\) and  
\[
u = (u_1,\dots,u_N) \in (\bK_1 \times \cdots \times \bK_r)^N, \quad u_j \neq 0 \text{ for all } j,
\]
there exists a tuple \((W_1,\dots,W_r) \in \mathscr{L}(\bK_1^N) \times \cdots \times \mathscr{L}(\bK_r^N)\) such that at least one \( W_i \neq \bK_i^N \) for some \( i \in \{1,\dots,r\} \), and for every \((\delta_1,\dots,\delta_N) \in \Delta^N\),  
\[
(\lambda_1 v_1, \dots, \lambda_r v_r) \in U_1 \times \cdots \times U_r
\implies (\lambda_1,\dots,\lambda_r) \in W_1 \times \cdots \times W_r,
\]
where  
\[
(v_1,\dots,v_r) = (u_1,\dots,u_N)^{\mathrm{tr}} \quad \text{and} \quad (\lambda_1,\dots,\lambda_r) = (\delta_1,\dots,\delta_N)^{\mathrm{tr}}.
\]
\end{lemma}

\begin{remark}
The \( r \)-tuple \( (W_1,\dots,W_r) \) depends on both \( \cU \) and \( u \).
\end{remark}

\begin{proof}
Fix \( (U_1,\dots,U_r) \in \supp(p_\theta) \) and let \( u_1,\dots,u_N \) be nonzero elements in \( \bK_1 \times \cdots \times \bK_r \). Since the transpose map  
\[
(u_1,\dots,u_N) \mapsto (v_1,\dots,v_r)
\]  
is bijective, the set  
\[
I = \{ i \in \{1,\dots,r\} \mid v_i \neq 0 \}
\]  
is nonempty.  For each \( i \in I \), let \( c_i = (c_{i,1},\dots,c_{i,N}) \in (\bK_i^*)^N \) be as given by Lemma~\ref{Lemma_c}, so that  
\[
\sum_{j=1}^N c_{i,j} w_{i,j} = 0, \quad \text{for all } w_i = (w_{i,1},\dots,w_{i,N}) \in U_i.
\]  
Now, for any \( (\delta_1,\dots,\delta_N) \in \Delta^N \), define  
\[
(\lambda_1,\dots,\lambda_r) = (\delta_1,\dots,\delta_N)^{\mathrm{tr}}.
\]  
If  
\[
(\lambda_1 v_1,\dots,\lambda_r v_r) \in U_1 \times \cdots \times U_r,
\]  
then for all \( i \in I \), we have  
\[
\sum_{j=1}^N c_{i,j} \lambda_{i,j} v_{i,j} = 0.
\]  
Since \( v_i \neq 0 \) for each \( i \in I \) and \( c_i \) consists of nonzero elements, it follows that  
\[
d_i = (c_{i,1} v_{i,1},\dots,c_{i,N} v_{i,N}) \neq 0.
\]  
Notably, \( d_i \) is independent of \( (\delta_1,\dots,\delta_N) \). Define the proper linear subspace  
\[
W_i = \Big\{ w_i = (w_{i,1},\dots,w_{i,N}) \in \bK_i^N \mid  
\sum_{j=1}^N d_{i,j} w_{i,j} = 0 \Big\} < \bK_i^N, \quad i \in I.
\]  
For \( i \notin I \), set \( W_i = \bK_i^N \). Thus, we conclude that  
\[
(\lambda_1,\dots,\lambda_r) \in W_1 \times \cdots \times W_r,
\]  
which completes the proof.
\end{proof}

\textbf{Proof of Theorem \ref{Thm_F}:} Fix \(0 < \varepsilon < \alpha\), and let
\[
\underline{f} = (f_1,\ldots,f_N) \in L^\infty(A)^N, \qquad m_A(f_j)=0 \ \text{for all } j=1,\ldots,N.
\]
Choose a parameter \(\varepsilon_0 > 0\), to be specified later, and pick functions 
\(g_1,\ldots,g_N \in L^\infty(A)\) such that
\[
m_A(g_j)=0, \qquad \|f_j - g_j\|_{L^1(A)} < \varepsilon_0, \qquad 
\|g_j\|_\infty \le \|f_j\|_\infty, \qquad 
G_j := \supp(\widehat{g}_j) \ \text{is finite},
\]
for each \(j=1,\ldots,N\). Such functions always exist. Indeed, take an approximate identity 
\((\rho_n) \) in  \(L^1(A)\) whose Fourier transforms \(\widehat{\rho}_n\) have finite 
support in \(\smash{\widehat{A}}\) for every \(n\) (see \cite[Theorem 33.12]{HeRo}).  
Setting \(g_j = \rho_n * f_j\), $j=1,\ldots,r$, we obtain functions satisfying all the conditions 
above for sufficiently large \(n\). \\

For every $(\delta_1,\ldots,\delta_N) \in \Delta^N$, we expand  
\begin{align*}
\theta(f_1 \circ \delta_1 \otimes \cdots \otimes f_N \circ \delta_N)
&= \theta(g_1 \circ \delta_1 \otimes \cdots \otimes g_N \circ \delta_N) \\
&\quad + \sum_{k=1}^{N} \theta\Big( \Big(\bigotimes_{l=1}^{k-1} g_l \circ \delta_l\Big) \otimes (f_k - g_k) \circ \delta_k \otimes \Big(\bigotimes_{l=k+1}^{N} f_l \circ \delta_l\Big) \Big),
\end{align*}  
where tensor products within the sum are ignored if the index ranges are empty. Since $\theta$ is an $N$-coupling of $(A,m_A)$, we deduce  
\[
\Big| \theta(f_1 \circ \delta_1 \otimes \cdots \otimes f_N \circ \delta_N)
- \theta(g_1 \circ \delta_1 \otimes \cdots \otimes g_N \circ \delta_N) \Big|
\leq N \eps_o \cdot \prod_{k=1}^N \|f_k\|_\infty.
\]  
Since each $g_j$ has finitely many nonzero Fourier coefficients and $m_A(g_j) = 0$, we can express  
\begin{align*}
\theta(g_1 \circ \delta_1 \otimes \cdots \otimes g_N \circ \delta_N)
&= \sum_{u_1 \neq 0} \cdots \sum_{u_N \neq 0} 
\widehat{g}_1(u_1) \cdots \widehat{g}_N(u_N) \widehat{\theta}(\delta_1u_1,\ldots,\delta_N u_N) \\
&= \sum_{u_1 \neq 0} \cdots \sum_{u_N \neq 0} 
\widehat{g}_1(u_1) \cdots \widehat{g}_N(u_N) \widehat{\tau}_\theta(\lambda_1 v_1,\ldots,\lambda_r v_r),
\end{align*}  
where the summation is finite. Using previous discussions, we rewrite this as  
\[
\sum_{u_1 \neq 0} \cdots \sum_{u_N \neq 0} \sum_{(U_1,\ldots,U_r)}
p_\theta(U_1,\ldots,U_r) \widehat{g}_1(u_1) \cdots \widehat{g}_N(u_N)
\chi_{U_1}(\lambda_1 v_1) \cdots \chi_{U_r}(\lambda_r v_r).
\]  
Since $p_\theta$ is a probability measure on a countable set, we can choose a finite subset $\cH \subset \supp(p_\theta)$ such that  
\[
\sum_{(U_1,\ldots,U_r) \notin \cH} p_\theta(U_1,\ldots,U_r) < 
\frac{\eps_o}{|G_1| \cdots |G_N|}.
\]  
and thus $\theta(g_1 \circ \delta_1 \otimes \cdots \otimes g_N \circ \delta_N)$ differs from the sum
\[
\sum_{u_1 \neq 0} \cdots \sum_{u_N \neq 0} \sum_{(U_1,\ldots,U_r) \in \cH}
p_\theta(U_1,\ldots,U_r) \widehat{g}_1(u_1) \cdots \widehat{g}_N(u_N)
\chi_{U_1}(\lambda_1 v_1) \cdots \chi_{U_r}(\lambda_r v_r).
\]
by at most 
\[
\eps_o \cdot \|f_1\|_\infty \cdots \|f_N\|_\infty,
\]
in absolute value, uniformly in $(\delta_1,\ldots,\delta_N)$. We conclude that if 
we choose 
\[
\eps_o < \frac{\eps}{(N+1)\|f_1\|_\infty \cdots \|f_N\|_\infty},
\]
then 
\[
\theta(f_1 \circ \delta_1 \otimes \cdots \otimes f_N \circ \delta_N)
\]
differs from
\[
\sum_{u_1 \neq 0} \cdots \sum_{u_N \neq 0} \sum_{(U_1,\ldots,U_r) \in \cH}
p_\theta(U_1,\ldots,U_r) \widehat{g}_1(u_1) \cdots \widehat{g}_N(u_N)
\chi_{U_1}(\lambda_1 v_1) \cdots \chi_{U_r}(\lambda_r v_r) 
\]
by at most $\eps$ in absolute value, uniformly over $(\delta_1,\ldots,\delta_N)$. \\

Now suppose that $(\delta_1,\ldots,\delta_N) \in S_\alpha(\theta,\underline{f})$, so
that
\[
\Big| \theta(f_1 \circ \delta_1 \otimes \cdots \otimes f_N \circ \delta_N) \big| \geq \alpha.
\]
Since $0 < \eps < \alpha$, we conclude that
\[
\sum_{u_1 \neq 0} \cdots \sum_{u_N \neq 0} \sum_{(U_1,\ldots,U_r) \in \cH}
p_\theta(U_1,\ldots,U_r) \widehat{g}_1(u_1) \cdots \widehat{g}_N(u_N)
\chi_{U_1}(\lambda_1 v_1) \cdots \chi_{U_r}(\lambda_r v_r) \neq 0,
\]
and thus there exist
\[
u = (u_1,\ldots,u_N) \in G_1 \times \cdots \times G_N \qand \cU = (U_1,\ldots,U_r) \in \cH
\] 
such that $u_j \neq 0$ for all $j$, and $(\lambda_1 v_1,\ldots,\lambda_r v_r) \in U_1 \times \cdots U_r$. By Lemma
\ref{Lemma_W}, there is an $r$-tuple 
\[
\cW(u,\cU) = (W_1,\ldots,W_r) \in \mathscr{L}(\bK_1^N \times \cdots \times \bK_r^N)
\]
such that $W_i \neq \bK_i^N$ for some $i$, and $(\lambda_1,\ldots,\lambda_r) \in W_1 \times \cdots \times W_r$. Define the finite set
\[
\cF(\theta,\underline{f},\alpha)
= \Big\{ \cW(u,\cU) \, : \, u \in G_1 \times \cdots G_N, \enskip (U_1,\ldots,U_r) \in \cH \Big\}.
\]
We have then shown that for all $(\delta_1,\ldots,\delta_N) \in S_\alpha(\theta,\underline{f})$, there exists $(W_1,\ldots,W_r) \in \cF(\theta,\underline{f},\alpha)$
such that $(\delta_1,\ldots,\delta_N)^{\mathrm{tr}} \in \cF(\theta,\underline{f},\alpha)$, finishing the proof of Theorem \ref{Thm_F}.

\subsection{An inductive criterion for coupling-resistance}

Throughout this sub-section, let $\bK_1,\ldots,\bK_{r+1}$ be countable fields and
let $\Gamma < \bK_1^* \times \cdots \times \bK_{r+1}^*$ be a subgroup. We denote by
$\Delta$ the projection of $\Gamma$ to $\bK_1^* \times \cdots \times \bK^*_{r}$, and
for each $i = 1,\ldots,r+1$, we write $\Gamma_i$ for the projection of $\Gamma$ to $\bK_i$. \\

Let $N \geq 2$. For $\underline{c} = (c_1,\ldots,c_N) \in (\bK_{r+1}^*)^N$, we define the complicated-looking set
\[
\Lambda_N(\underline{c}) 
= \left\{ 
(\lambda_1,\ldots,\lambda_r) \in \bK_1^N \times \cdots \times \bK_r^N \,\enskip  \middle| \, \enskip
\begin{aligned}
    &\exists \, (\delta_1,\ldots,\delta_N) \in \Delta^N, \enskip (a_1,\ldots,a_N) \in \bK^*_{r+1} \\[0.4cm]
    &\text{ such that } (\delta_j,a_j) \in \Gamma \quad \text{for all $j=1,\ldots,N$} , \\
    & (\lambda_1,\ldots,\lambda_r) = (\delta_1,\ldots,\delta_N)^{\textrm{tr}}, \qand 
    \sum_{j=1}^N c_j a_j = 0
\end{aligned}
\right\}.
\]

\begin{definition}  
A subgroup \(\Gamma < \bK_1^* \times \cdots \times \bK_{r+1}^{*}\) is said to be \textbf{left-scrambled} if there exist an integer \(N \geq 2\) and a tuple \(\underline{c} = (c_1,\ldots,c_N) \in (\bK_{r+1}^*)^N\) such that no finite collection  
\[
\cF \subset \mathscr{L}(\bK_1^N) \times \cdots \times \mathscr{L}(\bK_r^N)
\]  
satisfies the following property: for every \((W_1,\ldots,W_r) \in \cF\), there exists at least one index \(i\) for which \(W_i \neq \bK_i^N\), and  
\[
\Lambda_N(\underline{c}) \subseteq \bigcup_{(W_1,\ldots,W_r) \in \cF} W_1 \times \cdots \times W_r.
\] 
We refer to $(N,\underline{c})$ as the \textbf{parameters} of $\Gamma$. 
\end{definition}

\begin{theorem}
\label{Thm_CriterionCR}
Let \(\Gamma < \bK_1^* \times \cdots \times \bK_{r+1}^*\) be a subgroup satisfying the following conditions:
\begin{itemize}
    \item[$(i)$] Each projection \(\Gamma_i < \bK_i^*\) is measure-rigid for all \(i = 1,\ldots,r+1\). 
    \item[$(ii)$] The projected subgroup \(\Delta < \bK_1^* \times \cdots \times \bK^*_r\) is coupling-resistant. 
    \item[$(iii)$] \(\Gamma\) is left-scrambled.
\end{itemize}
Then \(\Gamma\) is coupling-resistant.
\end{theorem}

\subsection{Proof of Theorem \ref{Thm_CriterionCR}}

Let $\Gamma < \bK_1^* \times \cdots \bK_{r+1}^*$ be a subgroup satisfying the following conditions:
\begin{itemize}
\item[$(i)$] Each $\Gamma_i$ is measure-rigid for all $i=1,\ldots,r+1$. \vspace{0.1cm}
\item[$(ii)$] $\Delta < \bK_1^* \times \cdots \times \bK_r^*$ is coupling-resistant.\vspace{0.1cm}
\item[$(iii)$] $\Gamma$ is \emph{not} coupling-resistant.
\end{itemize} 
Our aim is to show that $\Gamma$ cannot be left-scrambed, meaning that for any choice of parameters $(N,\underline{c})$, there exists a finite collection
\[
\cF \subset \mathscr{L}(\bK_1^N) \times \cdots \times \mathscr{L}(\bK_r^N)
\]  
such that for every $(W_1,\ldots,W_r) \in \cF$, at least one index $i$
satisfies $W_i \neq \bK_i^N$, and
\[
\Lambda_N(\underline{c}) \subseteq \bigcup_{(W_1,\ldots,W_r) \in \cF} W_1 \times \cdots \times W_r.
\] 
Towards the proof of this, let us fix $N \geq 2$ and an $N$-tuple $\underline{c} = (c_1,\ldots,c_N) \in (\bK_{r+1}^*)^N$. \\

As above, we write $A = \widehat{\bK}_1 \times \cdots \times \widehat{\bK}_r$ and denote by $m_A$ the Haar probability measure on the compact group $A$. Since $\Gamma$ is assumed to not be coupling-resistant, we can find a $\Gamma$-invariant coupling $\eta$ of 
\[
(\widehat{\bK}_1,m_{\widehat{\bK}_1}),\ldots,(\widehat{\bK}_{r+1},m_{\widehat{\bK}_{r+1}})
\]
different than the product measure $\smash{m_A \otimes m_{\widehat{\bK}_{r+1}}}$. However, since $\Delta$ is assumed to be coupling-resistant, the projection of $\eta$ to $A$
must equal $m_A$, and thus $\eta$ is a $\Gamma$-invariant coupling between 
$\smash{(A,m_A)}$ and $\smash{(\widehat{\bK}_{r+1},m_{\widehat{\bK}_{r+1}})}$ which
is not a product. By Lemma \ref{Lemma_CouplingMarkov}, we can thus find 
a non-constant Markov intertwiner $S : L^2(A) \ra \smash{L^2(\widehat{\bK}_{r+1})}$ such that
\begin{equation}
\label{Seq}
S(f \circ \delta) = Sf \circ a, \quad \textrm{for all $\gamma = (\delta,a) \in \Gamma$}
\end{equation}
and 
\[
\eta(f \otimes \varphi) = \int_{\widehat{\bK}_{r+1}} Sf \cdot \varphi \, dm_{\widehat{\bK}_{r+1}}, \quad \textrm{for all $f \in L^\infty(A)$ and $\varphi \in L^\infty(\widehat{\bK}_{r+1})$}.
\]
Define $\theta_{\underline{c}} \in \Prob(A^N)$ by
\[
\theta_{\underline{c}}(f_1 \otimes \cdots \otimes f_N) 
= \int_{\widehat{\bK}_{r+1}} Sf_1(c_1^*\xi_{r+1}) \cdots Sf_N(c_N^*\xi_{r+1}) \, dm_{\widehat{\bK}_{r+1}}(\xi_{r+1}),
\]
for all $f_1,\ldots,f_N \in L^\infty(A)$. The equivariance property of $S$ above readily shows that $\theta_{\underline{c}}$ is a $\Delta$-invariant $N$-coupling 
of $(A,m_A)$. \\

The following key lemma will be proved in the next sub-section. 
\begin{lemma}
\label{Lemma_alphafo}
For any $(N,\underline{c})$, there exist $\alpha > 0$ and $f_o \in L^\infty(A)$ such that $m_A(f_o) = 0$
and
\[
\big| \theta_{\underline{c}}(f_o \circ \delta_1 \otimes \cdots \otimes f_o \circ \delta_N)\big| \geq \alpha, \quad \textrm{for all $(\delta_1,\ldots,\delta_N)^{\mathrm{tr}} \in \Lambda_N(\underline{c})$}.
\]
\end{lemma}
We now prove Theorem \ref{Thm_CriterionCR} using Lemma \ref{Lemma_alphafo}. This lemma can be concisely restated as the inclusion
\[
\Lambda_N(\underline{c}) \subseteq S_\alpha(\theta_{\underline{c}},(f_o,\ldots,f_o))^{\mathrm{tr}}. 
\]
Theorem \ref{Thm_F} tells us that there is a finite collection
\[
\cF = \cF(\theta,(f_o,\ldots,f_o),\alpha) \subset \mathscr{L}(\bK_1^N) \times \cdots \times \mathscr{L}(\bK_r^N)
\]
such that for every $(W_1,\ldots,W_r) \in \cF$, there is at least one index $i$ such 
that $W_i \neq \bK_i^N$, and 
\[
\Lambda_N(\underline{c}) \subseteq S_\alpha(\theta,(f_o,\ldots,f_o))^{\mathrm{tr}} \subseteq  \bigcup_{(W_1,\ldots,W_r)} W_1 \times \cdots \times W_r.
\]
Since $(N,\underline{c})$ is arbitrary, we conclude that the subgroup $\Gamma$ is not left-scrambed, finishing the proof of Theorem \ref{Thm_CriterionCR}. 

\subsection{Proof of Lemma \ref{Lemma_alphafo}}

We recall some notation and conventions from Subsection \ref{subsec:FT}. If $\bK$ is a field and $b \in \bK$, we define the character $\smash{b^{\vee}}$ on $\smash{\widehat{\bK}}$ by 
\[
b^{\vee}(\xi) = \xi(b), \quad \xi \in \widehat{\bK}.
\]
In particular, if $\tau$ is a probability measure on $\widehat{\bK}$ and $b_1,b_2 \in \bK$, then
\[
\widehat{\tau}(b_1,b_2) = \tau(b_1^{\vee} \otimes b_2^{\vee}),
\]
with our convention for the Fourier transform of a measure (see Subsection \ref{subsec:FT}). \\

We now present two general Fourier-analytic lemmas for a countably infinite field $\bK$, from which we then derive Lemma \ref{Lemma_alphafo}. 

\begin{lemma}
\label{Lemma_PosCorr}
Let \( \varphi_1, \dots, \varphi_N \in L^\infty(\widehat{\bK}) \) be functions whose Fourier transforms \( \widehat{\varphi}_1, \dots, \widehat{\varphi}_N \) are non-negative on \( \bK \). Then, for all \( (s_1, \dots, s_N) \in (\bK^*)^N \) satisfying \( s_1 + \dots + s_N = 0 \), we have  
\[
\int_{\widehat{\bK}} \varphi_1(s_1^* \xi) \cdots \varphi_N(s_N^* \xi) \, dm_{\widehat{\bK}}(\xi)
\geq \sum_{b \in \bK} \widehat{\varphi}_1(b) \cdots \widehat{\varphi}_N(b).
\]
\end{lemma}

\begin{proof}
We first assume that $\widehat{\varphi}_1,\ldots,\widehat{\varphi}_N \in \ell^1(\bK)^{+}$.  
Then
\[
\varphi_j(\xi) = \sum_{b_j \in \bK} \widehat{\varphi}_j(b_j)\, (-b_j)^{\vee}(\xi),
\qquad j=1,\ldots,N.
\]
For any $s_1,\ldots,s_N \in \bK^*$ with $s_1 + \cdots + s_N = 0$, we have
\begin{align*}
\int_{\widehat{\bK}} 
\varphi_1(s_1^*\xi)\cdots \varphi_N(s_N^*\xi)\, dm_{\widehat{\bK}}(\xi)
&=
\sum_{\mathclap{b_1,\ldots,b_N \in \bK}} \, 
\widehat{\varphi}_1(b_1)\cdots \widehat{\varphi}_N(b_N)
\int_{\widehat{\bK}}
\xi\!\left(-(b_1 s_1 + \cdots + b_N s_N)\right) dm_{\widehat{\bK}}(\xi)
\\[0.2cm]
&=
\sum_{b_1 s_1 + \cdots + b_N s_N = 0}
\widehat{\varphi}_1(b_1)\cdots \widehat{\varphi}_N(b_N) \\[0.2cm]
&\ge 
\sum_{b \in \bK} \widehat{\varphi}_1(b)\cdots \widehat{\varphi}_N(b),
\end{align*}
where the inequality follows from the non-negativity of 
$\widehat{\varphi}_1,\ldots,\widehat{\varphi}_N$ by restricting to the diagonal
$b_1=\cdots=b_N=b$.
The general case (where the Fourier transforms are still non-negative but not necessarily in 
$\ell^1(\bK)$) follows by a standard approximation argument.
\end{proof}

\begin{lemma}
\label{Lemma_FourierPosMarkov}
Let $G$ be a countable discrete group and let $\Gamma < G \times \bK^*$ be a subgroup
whose projection to $\bK^*$ is measure-rigid. Let $(X,\mu)$ be a Borel $G$-space and 
suppose that there is a non-product $\Gamma$-invariant coupling $\eta$ between $(X,\mu)$ and 
$(\smash{\widehat{\bK}},m_{\widehat{\bK}})$. Let 
\[
S : L^\infty(X,\mu) \ra L^\infty(\widehat{\bK}) \qand S^* : L^\infty(\widehat{\bK})
\ra L^\infty(X,\mu)
\]
denote the corresponding Markov intertwiners such that
\[
\eta(f_1 \otimes f_2) = \int_{\widehat{\bK}} Sf_1 f_2 \, dm_{\widehat{\bK}} = \int_X f_1 S^* f_2 \, d\mu,
\]
for all $f_1 \in L^\infty(X,\mu)$ and $f_2 \in L^\infty(\widehat{\bK})$. Then there exist $b_o, b_1 \in \bK \setminus \{0\}$ with the property that the function $f_o = S^*(b_o^{\vee}) \in L^\infty(X,\mu)$ satisfies
\[
m_{\widehat{\bK}}(Sf_o) = 0 \qand \widehat{Sf_o}(b_1) > 0 \qand \widehat{Sf_o}(b) \geq 0, \quad \textrm{for all $b \in \bK$}.
\]
\end{lemma}

\begin{proof}
Let $\Delta$ denote the projection of $\Gamma$ to $\bK^*$. Since $\eta$ is $\Gamma$-invariant, we
have
\begin{equation}
\label{Sstareq}
S^*(\varphi \circ \delta) = S^*(f) \circ g, \quad \textrm{for all $\gamma = (g,\delta) \in \Gamma$ and $\varphi \in L^\infty(\widehat{\bK})$}.
\end{equation}
Hence, if we define $\tau_S \in \Prob(\widehat{\bK}^2)$ by
\[
\tau_S(\varphi_1 \otimes \varphi_2) = \int_{X} S^*\varphi_1 \cdot S^*\varphi_2 \, d\mu, \quad \varphi_1, \varphi_2 \in L^\infty(\widehat{\bK}),
\]
then $\tau_S$ is a $\Delta$-invariant $2$-coupling of $(\widehat{\bK},m_{\widehat{\bK}})$. Indeed, 
since $S^*(1) = 1$ and $\mu(S^*\varphi) = m_{\widehat{\bK}}(\varphi)$ for all $\varphi \in \smash{L^\infty(\widehat{\bK})}$, we see that $\tau_S$ is a $2$-coupling. Moreover, fix $\delta \in \Delta$,
and pick $g \in G$ such that $\gamma = (g,\delta) \in \Gamma$. Then, by \eqref{Sstareq},
\begin{align*}
\tau_S(\varphi_1 \circ \delta \otimes \varphi_2 \circ \delta) 
&= \int_X S^*(\varphi_1)(g.x) 
S^*(\varphi_2)(g.x) \, d\mu(x) = \int_X S^*(\varphi_1)(x) 
S^*(\varphi_2)(x) \, d\mu(x) \\
&= \tau_S(\varphi_1  \otimes \varphi_2), 
\end{align*}
for all $\varphi_1,\varphi_2 \in L^\infty(\widehat{\bK})$. Since $\delta \in \Delta$ is arbitrary, $\tau_S$ is $\Delta$-invariant. \\

We claim that $\tau_S$ is not the product joining. Indeed, assume otherwise. Then,
\[
\int_X S^*(\varphi_1) S^*(\varphi_2) \, d\mu = \int_X S^*(\varphi_1) \, d\mu \int_X S^*(\varphi_2) \, d\mu, \quad \textrm{for all $\varphi_1,\varphi_2 \in L^\infty(\widehat{\bK})$}.
\]
In particular, if we take $\varphi_1 = \varphi_2$ and assume $\varphi$ is real, this forces
\[
\int_X S^*(\varphi_1)^2  \, d\mu = \left( \int_X S^*(\varphi_1) \, d\mu \right)^2,
\]
and thus $S^*(\varphi_1)$ is $\mu$-almost surely constant, and equal to $m_{\widehat{\bK}}(\varphi_1)$. Since $\varphi_1$ is arbitrary, this implies that
\[
\eta(f \otimes \varphi_1) = \int_X f S^*(\varphi) \, d\mu =  \mu(f) m_{\widehat{\bK}}(\varphi) = \mu \otimes m_{\widehat{\bK}}(f \otimes \varphi),
\]
contradicting our assumption that $\eta$ is non-product. \\

Since $\Delta$ is a measure-rigid subgroup of $\bK^*$, Lemma \ref{Lemma_MeasureRigidEquivalence} now guarantees that there is a unique probability measure $p$ on $\mathscr{L}(\bK^2)$ such that $p(\{0\}) = 0$ and 
\[
\tau_S = \sum_{W \leq \bK^2} p(W) m_{W^\perp}. 
\]
In particular, for all $b_o, b \in \bK$,
\[
\widehat{\tau}_S(b_o,b) = \sum_{W \leq \bK^2} p(W) \chi_W(b_o,b) \geq 0.
\]
Since $p(\{0\}) = 0$ and $p$ gives zero measure to the linear subspaces $\bK \times \{0\}$ and $\{0\} \times \bK$ (since $\tau_S$ is a $2$-coupling of the Haar measure), there exist $b_o, b_1 \in \bK \setminus \{0\}$ such that
\[
\widehat{\tau}(b_o,b_1) = \sum_{W \leq \bK^2} p(W) \chi_W(b_o,b_1) > 0.
\]
Fix such $b_o$ and $b_1$ and define $f_o = S^*(b_o^{\vee}) \in L^\infty(X,\mu)$. Then, for all $b \in \bK$, 
\begin{align*}
\widehat{S(f_o)}(b) &= \int_{\widehat{\bK}} S(f_o)(\xi) \xi(b) \, dm_{\widehat{\bK}}(\xi) = \eta(f_o \otimes b^{\vee}) = \int_X S^*(b_o^{\vee}) S^*(b^{\vee}) \, d\mu \\[0.2cm]
&= \tau_S(b_o^{\vee} \otimes b^{\vee}) = \sum_{W \leq \bK^2} p(W) \chi_W(b_o,b) \geq 0,
\end{align*}
and
\[
\widehat{S(f_o)}(b_1) =  \sum_{W \leq \bK^2} p(W) \chi_W(b_o,b_1) > 0.
\]
Finally, $m_{\widehat{\bK}}(Sf_o) = \mu(f_o) = m_{\widehat{\bK}}(b_o^{\vee}) = 0$, and the proof is done.
\end{proof}

\begin{proof}[Proof of Lemma \ref{Lemma_alphafo} assuming Lemma \ref{Lemma_PosCorr} and Lemma \ref{Lemma_FourierPosMarkov}]
We retain the notation from the previous subsections. In particular, let \( S \) and \( \theta_{\underline{c}} \) be as defined in Lemma \ref{Lemma_alphafo}. Applying Lemma \ref{Lemma_FourierPosMarkov} to \((X,\mu) = (A,m_A)\) and \(\bK = \bK_{r+1}\), we obtain the existence of a function \( f_o \in L^\infty(A) \) and an element \( b_1 \in \bK \setminus \{0\} \) such that  
\[
m_A(f_o) = 0, \quad \widehat{Sf_o}(b_1) > 0, \quad \text{and} \quad \widehat{Sf_o} \geq 0.
\]  

Now, fix \((a_1,\ldots,a_N) \in \Lambda_N(\underline{c})\). By construction, there exist \(\delta_1,\ldots,\delta_N \in \Delta\) such that  
\[
c_1 a_1 + \dots + c_N a_N = 0, \quad \text{and} \quad (\delta_j,a_j) \in \Gamma \text{ for all } j = 1,\ldots,N.
\]  

Using \eqref{Seq}, we obtain  
\[
\theta_{\underline{c}}(f_o \circ \delta_1 \otimes \cdots \otimes f_o \circ \delta_N)
= \int_{\bK_{r+1}} Sf_o((c_1a_1)^*.\xi_{r+1}) \cdots Sf_o((c_Na_N)^*.\xi_{r+1}) \, dm_{\widehat{\bK}_{r+1}}(\xi_{r+1}).  
\]  

Applying Lemma \ref{Lemma_PosCorr} with \( s_j = c_j a_j \) for \( j = 1,\ldots,N \), we deduce that  
\[
\theta_{\underline{c}}(f_o \circ \delta_1 \otimes \cdots \otimes f_o \circ \delta_N)
\geq \sum_{b \neq 0} \widehat{Sf_o}(b) > 0.
\]  

The term \( b = 0 \) is excluded from the sum since \( \widehat{Sf_o}(0) = m_A(f_o) = 0 \). Since the $N$-tuple \((a_1,\ldots,a_N) \in \Lambda_N(\underline{c})\) was arbitrary, we conclude that  
\[
\Lambda_{N}(\underline{c}) \subset S_\alpha(\theta,(f_o,\ldots,f_o))^{\mathrm{tr}},  
\]  
where \(\alpha = \sum_{b \neq 0} \widehat{Sf_o}(b) > 0\).  
\end{proof}

\section{Coupling resistance between non-isomorphic fields: Proof of Theorem \ref{Thm_main1}}
\label{Sec:ThmMain1}

Throughout this section, let \(\bK\) and \(\bL\) be countably infinite fields, equipped with their discrete topologies, and assume that \(\rho : \bK^* \to \bL^*\) is an \emph{injective} homomorphism. Our goal is to establish the following theorem, which, in combination with the results from Appendix \ref{App:FieldHomo}, leads to Theorem \ref{Thm_main1}. After stating the theorem, we provide a brief derivation of this implication.

\begin{theorem}
\label{Thm_FieldHomo}
Suppose that \(\Gamma_{\bL} = \rho(\bK^*) < \bL^*\) is a measure-rigid subgroup and that the two ergodic \(\bK^*\)-actions  
\[
\bK^* \curvearrowright (\widehat{\bK},m_{\widehat{\bK}}) \quad \text{and} \quad \bK^* \curvearrowright_{\rho} (\widehat{\bL},m_{\widehat{\bL}})
\]  
are not disjoint. Then, there exist finite-valued functions \( u,v : \bK^* \to \bL^* \) such that  
\[
\rho(1+x) = u(x) + v(x) \rho(x), \quad \text{for all } x \in \bK^*,
\]  
where we adopt the convention that $\rho(0) = 0$.
\end{theorem}

\subsection{Proof of Theorem \ref{Thm_main1} assuming Theorem \ref{Thm_FieldHomo}}

Let $\bK$ and $\bL$ be as above, and suppose $\rho : \bK^* \to \bL^*$ is an injective homomorphism such that the ergodic $\bK^*$-actions
\[
\bK^* \curvearrowright (\widehat{\bK},m_{\widehat{\bK}}) \quad \text{and} \quad \bK^* \curvearrowright_{\rho} (\widehat{\bL},m_{\widehat{\bL}})
\]
are not disjoint. Then, by Theorem \ref{Thm_FieldHomo} and Theorem \ref{Thm_Fields}, there exists a finite-valued homomorphism $w : \bK^* \to \bL^*$ such that the map $\kappa : \bK \to \bL$ defined by
\[
\kappa(0) = 0, \quad \kappa(x) = w(x) \rho(x), \quad \text{for all } x \in \bK^*
\]
is an injective field homomorphism. Moreover, if $\rho$ is surjective, then $\kappa$ is a field isomorphism. \\

Assume now that $\rho$ is a bijection, so that $\kappa$ is a field isomorphism, and let 
$\kappa^* : \widehat{\bL} \to \widehat{\bK}$ denote the induced continuous isomorphism between the Pontryagin duals (of the additive groups). 
Our goal is to show that every coupling on $\widehat{\bK} \times \widehat{\bL}$ which is invariant and ergodic under the $\bK^*$-action
\begin{equation}
\label{rho-action}
a.(\xi,\eta) = (a^*\xi,\rho(a)^*\eta), \qquad (\xi,\eta) \in \widehat{\bK} \times \widehat{\bL},
\end{equation}
is of the form $(e,b^*)_*\mu$ for some $b \in \bL^*$, where
\[
\mu(f)
= \frac{1}{|\im(w)|}
\sum_{c \in \im(w)}
\int_{\widehat{\bL}} f(\kappa^*(\eta),c^*\eta)\, dm_{\widehat{\bL}}(\eta),
\qquad f \in C(\widehat{\bK}\times\widehat{\bL}).
\]
It is straightforward to verify that $\mu$ is $\bK^*$-invariant and $\bK^*$-ergodic under this action; we leave this to the reader.

\medskip

Let $\theta$ be an invariant and ergodic coupling on $\widehat{\bK}\times\widehat{\bL}$.  
Define $\Delta=\ker(w) < \bK^*$. Since $w$ has finite image, $\Delta$ has finite index in $\bK^*$, and therefore $\Delta$ is measure-rigid by Lemma~\ref{Lemma_FiniteIndexMeasureRigid}.  
Moreover, for every $\delta\in\Delta$,
\[
\kappa(\delta)=w(\delta)\rho(\delta)=\rho(\delta),
\]
and hence for all $b\in\bK$ and $\eta\in\widehat{\bL}$,
\begin{equation}
\label{delta-equi}
(\delta^*\cdot \kappa^*(\eta))(b)
= \eta(\kappa(\delta b))
= \eta(\rho(\delta)\kappa(b))
= (\kappa^*(\rho(\delta)^*\eta))(b).
\end{equation}
Thus $\kappa^*$ is equivariant with respect to the $\Delta$-actions on $\widehat{\bL}$ and $\widehat{\bK}$. Consequently,
\[
\widetilde{\theta} := (\mathrm{id}\times \kappa^*)_* \theta
\]
is a $\Delta$-invariant $2$-coupling of $(\widehat{\bK},m_{\widehat{\bK}})$.  
Since $\Delta$ is measure-rigid, Lemma~\ref{Lemma_MeasureRigidEquivalence} ensures the existence of a unique probability measure  
$p_\theta$ on $\mathscr{L}(\bK^2)$ such that
\[
\widetilde{\theta}
= \sum_{W\le\bK^2} p_\theta(W)\, m_{W^\perp}.
\]
Because $\widetilde{\theta}$ is a $2$-coupling, the support of $p_\theta$ consists only of the trivial subspace $\{(0,0)\}$ and the lines
\[
W_a=\{(b,-ab)\mid b\in\bK\}, \qquad a\in\bK^*.
\]
Thus,
\[
\theta(f)
= p_\theta(\{(0,0)\})\cdot (m_{\widehat{\bK}}\otimes m_{\widehat{\bL}})(f)
+ \sum_{a\in\bK^*} p_\theta(W_a)
\int_{\widehat{\bL}} f(a^*\kappa^*(\eta),\eta)\, dm_{\widehat{\bL}}(\eta),
\qquad f\in C(\widehat{\bK}\times\widehat{\bL}).
\]

Since $\theta$ is assumed to be a non-trivial $\bK^*$-invariant and $\bK^*$-ergodic coupling, and since 
the product measure $m_{\widehat{\bK}}\otimes m_{\widehat{\bL}}$ is itself $\bK^*$-invariant and $\bK^*$-ergodic, we must have  
$p_\theta(\{(0,0)\})=0$.  
Using the isomorphism $\kappa$, we may rewrite the sum over $\bK^*$ as a sum over $\bL^*$ and obtain
\[
\theta = \sum_{b\in\bL^*} p(b)\, (e,b^*)_*\tau,
\]
where $p(b)=p_\theta(W_{\kappa^{-1}(b)^{-1}})$ is a probability measure on $\bL^*$, and
\[
\tau(f)=\int_{\widehat{\bL}} f(\kappa^*(\eta),\eta)\, dm_{\widehat{\bL}}(\eta).
\]
Indeed, using $(\kappa^{-1}(b))^*\kappa^*(\eta)=\kappa^*(b^*\eta)$ for all $b\in\bL^*$ and $\eta\in\widehat{\bL}$,
\begin{align*}
\sum_{a\in\bK^*} p_\theta(W_a)
\int_{\widehat{\bL}} f(a^* \kappa^*(\eta),\eta)\,dm_{\widehat{\bL}}(\eta)
&=
\sum_{b\in\bL^*} p_\theta(W_{\kappa^{-1}(b)})
\int_{\widehat{\bL}}
f((\kappa^{-1}(b))^*\kappa^*(\eta), \eta)\, dm_{\widehat{\bL}}(\eta)
\\[0.15cm]
&=
\sum_{b\in\bL^*} p_\theta(W_{\kappa^{-1}(b)})
\int_{\widehat{\bL}}
f(\kappa^*(b^*\eta), \eta)\, dm_{\widehat{\bL}}(\eta)
\\[0.15cm]
&=
\sum_{b\in\bL^*} p_\theta(W_{\kappa^{-1}(b^{-1})})
\int_{\widehat{\bL}}
f(\kappa^*(\eta), b^*\eta)\, dm_{\widehat{\bL}}(\eta),
\end{align*}
where we used the $\bL^*$-invariance of $m_{\widehat{\bL}}$ in the last step.

\medskip

A direct computation shows that
\[
(a^*,\rho(a)^*)\tau(f)
= \int_{\widehat{\bL}} f(\kappa^*(\eta),(w(a)^{-1})^*\eta)\, dm_{\widehat{\bL}}(\eta),
\]
equivalently,
\[
(a^*,\rho(a)^*)\tau = (e,(w(a)^{-1})^*)_*\tau,
\qquad a\in\bK^*.
\]
Since $\theta$ is invariant under all $(a^*,\rho(a)^*)$, the probability measure $p$ on $\bL^*$ must be invariant under the finite subgroup $\im(w)$. Consequently, $p$ descends to a probability measure on the quotient group $\bL^*/\im(w)$. Thus,
\[
\theta
= \sum_{[b]\in \bL^*/\im(w)} p(b)\, |\im(w)| \, (e,b^*)_*\mu.
\]
Finally, as observed above, $\mu$ and hence each $(e,b^*)_*\mu$ is invariant and ergodic for the $\bK^*$-action \eqref{rho-action}.  
Since $\theta$ is itself invariant and ergodic, it must coincide with one such extremal measure.  
Therefore,
\[
\theta = (e,b^*)_*\mu
\]
for some $b\in \bL^*/\im(w)$, completing the proof.

\subsection{Proof of Theorem \ref{Thm_FieldHomo}}

Let \(\bK_1 = \bL\) and \(\bK_2 = \bK\), and define the subgroup  
\[
\Gamma_\rho = \{ (\rho(a),a) \mid a \in \bK^* \}.
\]
Observe that \((\Gamma_\rho)_1 = \Gamma_{\bL}\) and \(\Gamma_2 = \bK^*\) are both measure-rigid. Moreover, \(\Delta = (\Gamma_\rho)_1\) is trivially coupling-resistant. By assumption, the \(\bK^*\)-actions  
\[
\bK^* \curvearrowright (\widehat{\bK},m_{\widehat{\bK}}) \quad \text{and} \quad \bK^* \curvearrowright_{\rho} (\widehat{\bL},m_{\widehat{\bL}})
\]  
are not disjoint, implying that \(\Gamma_\rho\) is not coupling-resistant. Therefore, by Theorem \ref{Thm_CriterionCR}, \(\Gamma_\rho\) cannot be left-scrambled. \\

In particular, setting \(N = 3\) and \(\underline{c} = (1,1,1)\), we conclude that there exists a finite collection \(\mathcal{F} = \{ W^{(1)},\ldots,W^{(M)} \}\) of \emph{proper} subspaces of \(\bK^3\) such that  
\[
\Lambda_3(1,1,1) = \{ (\rho(a_1),\rho(a_2),\rho(a_3)) \in (\bL^*)^3 \mid a_1 + a_2 + a_3 = 0 \} 
\subseteq \bigcup_{k=1}^M W^{(k)}.
\]
Since each \(W^{(k)}\) is a proper linear subspace, there exists a nonzero vector  
\[
(\beta_1^{(k)},\beta_2^{(k)},\beta_3^{(k)}) \in \bL^3
\]
such that  
\[
\beta_1^{(k)} w_1 + \beta_2^{(k)} w_2 + \beta_3^{(k)} w_3 = 0, 
\quad \text{for all } (w_1,w_2,w_3) \in W^{(k)}.
\]
For every \(x \in \bK^* \setminus \{-1\}\), choose \(k_x \in \{1,\ldots,M\}\) such that  
\[
(1,\rho(x),-\rho(1+x)) \in W^{(k_x)},
\]
which implies  
\[
\beta_1^{(k_x)} + \beta_2^{(k_x)} \rho(x) = \beta_3^{(k_x)} \rho(1+x).
\]
Note that at least one of \(\beta_1^{(k_x)}, \beta_2^{(k_x)}, \beta_3^{(k_x)}\) must be nonzero, and since \(\rho\) never takes the value zero, at most one of them can be zero. \\

Define the set
\[
P = \{ x \in \bK^* \setminus \{-1\} \mid \beta_1^{(k_x)}, \beta_2^{(k_x)}, \beta_3^{(k_x)} \neq 0 \},
\]
and for \(x \in P\), 
\[
u(x) = \frac{\beta_1^{(k_x)}}{\beta_3^{(k_x)}} \quad \text{and} \quad v(x) = \frac{\beta_2^{(k_x)}}{\beta_3^{(k_x)}},
\]
so that \(u(x), v(x) \neq 0\) and $\rho(1+x) = u(x) + v(x) \rho(x)$. Since \(\mathcal{F}\) consists of only \(M\) subspaces, \(u\) and \(v\) each attain at most \(M^2\) different values on \(P\). \\

For \(j=1,2,3\), define  
\[
N_j = \{ x \in \bK^* \setminus \{-1\} \mid \beta_j^{(k_x)} = 0, \quad \beta_i^{(k_x)} \neq 0 \text{ for } i \neq j \}.
\]
If \(x \in N_1\), then  
\[
\rho\Big( \frac{x}{1+x} \Big) = \frac{\beta_3^{(k_x)}}{\beta_2^{(k_x)}}.
\]
Since the right-hand side takes only finitely many values, \(N_1\) must be finite. A similar argument applies to \(N_2\) and \(N_3\), so setting \(N = N_1 \cup N_2 \cup N_3\), we conclude that \(N\) is finite and  
\[
\bK^* = P \cup N \cup \{-1\}.
\]
For \(x \in N \), define  
\[
u(x) = \rho(1+x) - \rho(x), \quad v(x) = 1,
\]
so that \(\rho(1+x) = u(x) + v(x) \rho(x)\). Since \(\rho\) is injective, \(u(x) \neq 0\), and since \(N\) is finite, \(u\) and \(v\) attain only finitely many values on \(N\). Finally, for \(x = -1\), set \(u(-1) = v(-1) = 1\), ensuring that  
\[
\rho(0) = 0 = u(-1) + v(-1) \rho(-1).
\]
Thus, the functions \(u,v\) take only finitely many values and satisfy  
\[
\rho(1+x) = u(x) + v(x) \rho(x), \quad \text{for all } x \in \bK^*,
\]
with the convention that \(\rho(0) = 0\).

\section{Coupling resistance of moment curves: Proof of Theorem \ref{Thm_DiagonalMain}}
\label{Sec:ThmDiagonal}

In view of Lemma \ref{Lemma_CouplingResistantEquivalence}, Theorem \ref{Thm_DiagonalMain} immediately follows from the following result.

\begin{theorem}
\label{Thm_RelPrime}
Let $n_1,\ldots,n_r$ be non-zero integers, not divisible by the characteristics of $\bK$. Then the subgroup 
\[
\Gamma = \{ (a^{n_1},\ldots,a^{n_r}) \, \mid \, a \in \bK^* \} < (\bK^*)^r
\]
is coupling-resistant. 
\end{theorem}

\subsection{Proof of Theorem \ref{Thm_RelPrime}}

The following lemma is the key to the proof of Theorem \ref{Thm_RelPrime}. 

\begin{lemma}
\label{Lemma_p1pN}
Let \( n \) be a positive integer, not divisible by the characteristic of $\bK$. There exists a positive integer \( N \), depending only on \( n \) and the characteristic of \( \bK \), along with polynomials \( p_1, \dots, p_N \) over \( \bK \) and \( c_1, \dots, c_N \in \bK^* \), such that  
\[
c_1 p_1(t)^n + \dots + c_N p_N(t)^n = 0, \quad \text{for all } t \in \bK.
\]
Moreover, for every integer \( m \), either negative or greater than \( n \) and not divisible by the characteristic of $\bK$, the rational functions  
\[
p_1^m, \dots, p_N^m : \bK \setminus Z \to \bK,
\]
where \( Z \) is the union of the zero sets of \( p_1, \dots, p_N \), are linearly independent as \(\bK\)-valued functions on \( \bK \setminus Z \).
\end{lemma}

\begin{proof}[Proof of Theorem \ref{Thm_RelPrime} assuming Lemma \ref{Lemma_p1pN}]
We proceed by induction. The base case $r = 1$ is trivial. Suppose the theorem holds for some $r \geq 1$ and all distinct non-zero integers $n_1, \dots, n_r$, not divisible by the characteristic of $\bK$. Our goal is to establish the result for any $(r+1)$-tuple $n_1, \dots, n_{r+1}$ of distinct nonzero integers that are not divisible by the characteristic of $\bK$. \\

Define the subgroup  
\[
\Gamma = \{ (a^{n_1},\dots,a^{n_r},a^{n_{r+1}}) \mid a \in \bK^* \} \subset (\bK^*)^{r+1}.
\]
We aim to show that $\Gamma$ is coupling-resistant. To do this, we introduce the index sets  
\[
I_{+} = \{ k \in \{1,\dots,r+1\} \mid n_k > 0 \}, \quad I_{-} = \{ k \in \{1,\dots,r+1\} \mid n_k < 0 \}.
\]
Since $\Gamma$ remains unchanged if we replace $(n_1, \dots, n_{r+1})$ with $(-n_1, \dots, -n_{r+1})$, we may assume without loss of generality that $I_+$ is nonempty. Moreover, by permuting indices if necessary, we may further assume that $n_{r+1}$ is the smallest element in $I_+$. \\

Define  
\[
\Delta = \{ (a^{n_1},\dots,a^{n_r}) \mid a \in \bK^* \} \subset (\bK^*)^r
\]
as the projection of $\Gamma$ onto the first $r$ coordinates, and let $\Gamma_i$ denote the projection of $\Gamma$ onto the $i$th coordinate for each $i = 1, \dots, r+1$. We note:
\begin{itemize}
    \item By Lemma \ref{Lemma_KnIsMeasureRigid}, each $\Gamma_i$ is measure-rigid for $i=1, \dots, r+1$.
    \item By assumption, $\Delta$ is coupling-resistant.
\end{itemize}

By Theorem \ref{Thm_CriterionCR}, to prove that $\Gamma$ is coupling-resistant, it suffices to show that $\Gamma$ is left-scrambled. \\
 
Applying Lemma \ref{Lemma_p1pN} with $n = n_{r+1}$, we obtain a positive integer $N$, polynomials $p_1, \dots, p_N$, and coefficients $\underline{c} = (c_1, \dots, c_N) \in (\bK^*)^N$ satisfying  
\[
c_1 p_1(t)^{n_{r+1}} + \dots + c_N p_N(t)^{n_{r+1}} = 0, \quad \forall t \in \bK.
\]
Furthermore, for each $i = 1, \dots, r$, the rational functions  
\[
\{ p_1(t)^{n_i}, \dots, p_N(t)^{n_i} \}, \quad t \in \bK \setminus Z,
\]
where $Z$ is the union of the zero sets of $p_1, \dots, p_N$, are linearly independent over $\bK$. We claim that $(N, \underline{c})$ serve as parameters for $\Gamma$. \\

To prove this, we introduce some notation. Define the rational functions  
\[
p_{i,j}(t) = p_j(t)^{n_i}, \quad i=1,\dots,r, \quad j=1,\dots,N.
\]
Set  
\[
\delta_j(t) = (p_{1,j}(t), \dots, p_{r,j}(t)), \quad a_j(t) = p_{r+1,j}(t),
\]
for $j=1,\dots,N$. Then, for all $t \in \bK \setminus Z$,  
\[
\delta_j(t) \in \Delta, \quad \gamma_j(t) := (\delta_j(t), a_j(t)) \in \Gamma.
\]
Moreover, define  
\[
(\lambda_1(t), \dots, \lambda_r(t)) = (\delta_1(t), \dots, \delta_N(t))^{\text{tr}},
\]
where  
\[
\lambda_i(t) = (p_{i,1}(t), \dots, p_{i,N}(t)) \in (\bK^*)^N, \quad \forall t \in \bK \setminus Z.
\]
Since  
\[
c_1 a_1(t) + \dots + c_N a_N(t) = 0, \quad \forall t \in \bK \setminus Z,
\]
it follows that  
\[
(\lambda_1(t), \dots, \lambda_r(t)) \in \Lambda_N(\underline{c}), \quad \forall t \in \bK \setminus Z.
\]
Assume, for contradiction, that $\Gamma$ is not left-scrambled with parameters $(N, \underline{c})$. Then there exists a finite collection  
\[
\cF \subset \mathscr{L}(\bK^N) \times \dots \times \mathscr{L}(\bK^N)
\]
such that, for every $(W_1, \dots, W_r) \in \cF$, there exists at least one index $i$ with $W_i \neq \bK^N$ and  
\[
\Lambda_N(\underline{c}) \subseteq \bigcup_{(W_1, \dots, W_r)} W_1 \times \dots \times W_r.
\]
Since $(\lambda_1(t), \dots, \lambda_r(t)) \in \Lambda_N(\underline{c})$ for all $t \in \bK \setminus Z$, there must exist some index $i$ and a \emph{proper} linear subspace $W \subset \bK^N$ such that $\lambda_i(t) \in W$ for infinitely many $t$. However, as each component of $\lambda_i$ is a rational function, this implies  
\[
\lambda_i(t) \in W, \quad \forall t \in \bK \setminus Z.
\]
Since $W$ is a proper subspace, there exist nonzero coefficients $b_1, \dots, b_N$ such that  
\[
\sum_{j=1}^N b_j p_{i,j}(t) = \sum_{j=1}^N b_j p_j(t)^{n_i} = 0, \quad \forall t \in \bK \setminus Z.
\]
This contradicts the assumption that $\{ p_1(t)^{n_i}, \dots, p_N(t)^{n_i} \}$ are linearly independent. Thus, $\Gamma$ is left-scrambled, completing the proof. 
\end{proof}

\subsection{Proof of Lemma \ref{Lemma_p1pN}}

Note that
\begin{equation}
\label{tnexp}
(t^n + 1)^n = \sum_{k=0}^n \binom{n}{k} t^{kn}, \quad \textrm{for all $t \in \bK$},
\end{equation}
and let $P = \{ k = 0,\ldots,n \, \mid \, \binom{n}{k} \neq 0\}$ and $N = |P|+1$, where the binomial coefficients are viewed as elements in $\bK$. If $n = 1$, $P = \{0,1\}$ and $N = 3$, and if $n \geq 2$, then $P \supseteq \{0,1,n-1,n\}$ and thus $N \geq 5$, since $n$ is not divisible by the characteristic of $\bK$. If the characteristic is zero, then $|P| = n+1$ and $N = n+2$. Let 
\[
0 = k_1 < k_2 < \ldots < k_{N-1} = n 
\]
denote the elements in $P$, and define the polynomials 
\[
p_1(t) = 1, \enskip p_j(t) = t^{k_j}, \quad j=2,\ldots,N-1, \quad p_N(t) = t^n + 1,
\]
and
\[
c_j = \binom{n}{k_j}, \quad j=1,\ldots,N-1, \quad c_N = -1.
\]
Clearly, $c_1,\ldots,c_N \in \bK^*$ and it follows from \eqref{tnexp} that
\[
c_1 p_1(t)^n + \ldots + c_N p_N(t)^n = 0, \quad \textrm{for all $t \in \bK$}.
\]
Let $Z = \{0\} \cup \{ t \in \bK^* \, : \, t^n = -1\}$, and note that the polynomials 
$p_1,\ldots,p_N$ are all non-zero off the set $Z$. Let \( m \) be a non-zero integer, either negative or greater than \( n \) and not divisible by the characteristic of $\bK$. Our aim is to show that the rational functions  
\[
p_1^m, \dots, p_N^m : \bK \setminus Z \to \bK,
\]
are linearly independent as \(\bK\)-valued functions on \( \bK \setminus Z \). To do this, let us fix $b_1,\ldots,b_N$ such that
\begin{equation}
\label{bjs}
\sum_{j=1}^N b_j p_j(t)^m = b_N (t^n + 1)^m + b_{N-1} t^{nm} + \sum_{j=1}^{N-2} b_j t^{k_j m} = 0, \quad \textrm{for all $t \in \bK \setminus Z$}.
\end{equation}
We want to show that $b_1 = \ldots = b_N = 0$. We first note that if $b_N = 0$, then the remaining $b_1,\ldots,b_{N-1}$ must be zero too, since the rational functions 
$\{t^i\}_{i \in \bZ}$ on $\bK \setminus \{0\}$ are clearly linearly independent, since 
$\bK$ is infinite. So from now on we may assume that $b_{N} \neq 0$, and derive a contradiction. Our analysis splits into two cases. \\

\noindent \textbf{Case 1: $m < 0$.} \\

Multiply both sides of \eqref{bjs} by \( t^{n|m|}(t^n+1)^{|m|} \) to obtain the polynomial identity  
\[
b_N t^{n|m|} + b_{N-1}(t^n + 1)^{|m|} + \sum_{j=1}^{N-2} b_j t^{|m|(n-k_j)} (t^n+1)^{|m|} = 0, \quad \text{for all } t \in \bK \setminus Z.
\]
Since \( \bK \) is infinite and the left-hand side is a polynomial, this identity must hold for all \( t \in \bK \). Substituting \( t = 0 \) gives \( b_{N-1} = 0 \), since \( k_j < n \) for all \( j = 1, \dots, N-2 \), reducing the equation to  
\[
b_N t^{n|m|} + \sum_{j=1}^{N-2} b_j t^{|m|(n-k_j)} (t^n+1)^{|m|} = 0, \quad \text{for all } t \in \bK.
\]
Dividing both sides by \( t^{|m|(n-k_{N-2})} \) yields  
\[
b_N t^{k_{N-2}|m|} + \sum_{j=1}^{N-2} b_j t^{|m|(k_{N-2}-k_j)} (t^n+1)^{|m|} = 0, \quad \text{for all } t \in \bK.
\]
Setting \( t = 0 \) now forces \( b_{N-2} = 0 \). If \( N = 3 \), this establishes \( b_1 = b_2 = 0 \), implying \( b_3 = 0 \). If \( N > 3 \), we proceed iteratively: at each step, we divide by the monomial of the smallest degree and substitute \( t = 0 \), eventually concluding that \( b_{N-3}, \dots, b_1 = 0 \), which ultimately forces \( b_N = 0 \). \\

\noindent \textbf{Case 2: $m > n$.} \\

Since \( m \) is not divisible by the characteristic of \( \bK \), the coefficient of \( t^{n(m-1)} \) in the binomial expansion of \( (t^n + 1)^m \) is nonzero. For \eqref{bjs} to hold, this \( t^{n(m-1)} \)-term must be canceled by another term of the same degree in the sum over \( j = 1, \dots, N-2 \). That is, the equation  
\[
k_j m = n(m-1)
\]
must have a solution for some \( j \in \{1, \dots, N-2\} \). Rewriting this equation as  
\[
n = m(n - k_j),
\]
we see that no such solution exists, since \( m > n \) and \( n - k_j > 1 \) for all \( j = 1, \dots, N-2 \).

\section{Multiplicative asymmetry of F\o lner sequences: Proof of Theorem \ref{Thm_MultAssym_Main}}
\label{Sec:MultAssym}

Let $\bK$ be a countable discrete field. Our proof of Theorem \ref{Thm_MultAssym_Main} relies on a classical result stating that the discrete group $\PGL_2(\bK)$ is amenable if and only if $\bK$ is a union of finite subfields. This result is well-documented, with proofs available in \cite[Proposition 9]{HR} and \cite[Proposition 11]{Bekka}.
Given this result, it is evident that the following theorem directly implies Theorem \ref{Thm_MultAssym_Main}.

\begin{theorem}[Multiplicative asymmetry]
\label{Thm_MultASymm2}
Let \(\bK\) be a countable field and suppose that there exists a F\o lner sequence \((F_n)\) in \((\bK,+)\) such that \(((F_n \setminus \{0\})^{-1})\) is also a F\o lner sequence. Then the group $\PGL_2(\bK)$ is amenable. 
\end{theorem}

\subsection{Inducing F\o lner sequences in $\PGL_2(\bK)$: Proof of Theorem \ref{Thm_MultASymm2}}

Let us introduce some notation. If $X$ is a set, let $X^{(3)}$ denote the set of 
all $3$-tuples of distinct points in $X$. If $X$ is finite, then
\[
|X^{(3)}| = |X| \cdot (|X|-1) \cdot (|X|-2).
\]
Let $\bP(\bK^2)$ denote the projective space of $\bK^2$, and recall that $\PGL_2(\bK)$
is the image of $\GL_2(\bK)$ under the map
\[
g \mapsto [g], \quad 
\textrm{where $[g].[v] = [gv]$ for $g \in \GL_2(\bK)$ and $[v] \in \bP(\bK^2)$},
\]
with values inside the set of bijections on $\bP(\bK^2)$. 
It is well-known that $\PGL_2(\bK)$ acts freely and transitively on the set
$\bP(\bK^2)^{(3)}$, and that the group $\PGL_2(\bK)$ is generated by elements of the form
\[
u_{+}(b) = 
\left[ 
\left(
\begin{matrix}
1 & b \\
0 & 1
\end{matrix}
\right)
\right]
\qand
u_{-}(b) = 
\left[ 
\left(
\begin{matrix}
1 & 0 \\
b & 1
\end{matrix}
\right)
\right], \quad b \in \bK.
\]
Define the \emph{injective} map $\iota : \bK^{(3)} \ra \bP(\bK^2)^{(3)}$ by
\[
\iota(x_1,x_2,x_3)
= \left( 
\left[ 
\left(
\begin{matrix}
x_1  \\
1 
\end{matrix}
\right)
\right], 
\left[ 
\left(
\begin{matrix}
x_2  \\
1 
\end{matrix}
\right)
\right], 
\left[ 
\left(
\begin{matrix}
x_3  \\
1 
\end{matrix}
\right)
\right]
\right),
\]
and for $b \neq 0$, we define the map $\varphi_b : \bK \setminus \{-\frac{1}{b}\} \ra \bK$ by
\[
\varphi_b(x) = \frac{x}{1+bx}, \quad x \in \bK \setminus \left\{-\frac{1}{b}\right\}.
\]
We see that $\iota$ satisfies the equivariance properties:
\begin{align*}
u_{+}(b).\iota(x_1,x_2,x_3) &= \iota(x_1 + b,x_2+b,x_3+b) \\[0.2cm]
u_{-}(b).\iota(y_1,y_2,y_3) &= \iota(\varphi_b(y_1),\varphi_b(y_2),\varphi_b(y_3)),
\end{align*}
for all $b \in \bK$ and $(x_1,x_2,x_3), (y_1,y_2,y_3) \in \bK^{(3)}$ with $
y_j \neq -\frac{1}{b}$ for all $j=1,2,3$. Finally, let
\[
\ell_1 = \left[ 
\left(
\begin{matrix}
1  \\
0 
\end{matrix}
\right)
\right], \enskip
\ell_2 =  
\left[ 
\left(
\begin{matrix}
0  \\
1 
\end{matrix}
\right)
\right],
\enskip
\ell_3 =  
\left[ 
\left(
\begin{matrix}
1  \\
1 
\end{matrix}
\right)
\right].
\]
Given a finite set $A \subset \bK$, define
\[
Q(A) = \{ [g] \in \PGL_2(\bK) \, : \, ([g].\ell_1,[g].\ell_2,[g].\ell_3) \in \iota(A^{(3)}) \}.
\]
Since $\iota$ is injective and $\PGL_2(\bK)$ acts transitively and freely on $\bP(\bK^{(3)})$, we see that
\[
|Q(A)| = |A| \cdot (|A| - 1) \cdot (|A| - 2), \quad \textrm{for every finite set $A \subset \bK$},
\]
and for all finite sets $A, B \subset \bK$, we have
\[
Q(A \cap B) = Q(A) \cap Q(B). 
\]
Moreover, it follows from the equivariance properties for $\iota$ above that
\[
u_{+}(b)Q(A) = Q(A+b) \qand
u_{-}(b)(A) \supseteq Q\left(\varphi_b\left(A \setminus \left\{ -\frac{1}{b} \right\}\right)\right),
\quad \textrm{for all $b \in \bK \setminus \{0\}$}.
\]
We need the following simple lemma.

\begin{lemma}
Let \( A \subset \bK \setminus \{0\} \) be a finite set satisfying \( A^{-1} = A \). Then, for all \( b \in \bK \setminus \{0\} \), we have  
\[
\big| A \cap \varphi_b\big(A \setminus \{-\frac{1}{b}\}\big) \big|  
= \big| (A \setminus \{0\}) \cap (A + b) \big|.
\]
\end{lemma}

\begin{proof}
Fix \( b \in \bK \setminus \{0\} \). Since \( A^{-1} = A \), we have  
\[
A_b := A \cap \varphi_b\left(A \setminus \left\{ -\frac{1}{b} \right\}\right)
= 
A^{-1} \cap \varphi_b\left(A^{-1} \setminus \left\{ -\frac{1}{b} \right\}\right).
\]  
Choose \( x' \in A_b \). Then there exist \( x \in A \setminus \{0\} \) and \( y \in A \setminus \{-b\} \) such that  
\[
x' = \frac{1}{x} = \varphi_b\left(\frac{1}{y}\right) = \frac{1}{y+b}.
\]  
In particular, this implies that \( x \in A^b := (A \setminus \{0\}) \cap (A+b) \), and therefore \( |A_b| \leq |A^b| \).  Conversely, suppose \( x \in A^b \). Then there exists \( y \in A \) such that \( x = y+b \). Since \( x \neq 0 \), we must have \( y \neq -b \), and thus  
\[
\frac{1}{x} = \frac{1}{y+b} = \varphi_{b}\Big( \frac{1}{y} \Big).
\]  
Since \( A = A^{-1} \), it follows that  
\[
x' = \frac{1}{x} \in A, \quad y' = \frac{1}{y} \in A \setminus \left\{ -\frac{1}{b} \right\}.
\]  
Hence, \( x' \in A_b \), which gives \( |A^b| \leq |A_b| \).  
\end{proof}

We now proceed with the proof of Theorem \ref{Thm_MultASymm2}. Our goal is to show that if \( (F_n) \) is a F\o lner sequence in \( (\bK,+) \) such that \( ((F_n \setminus \{0\}^{-1}) \) is also a F\o lner sequence, then the sequence \( (\widetilde{F}_n) \), defined by  
\[
\widetilde{F}_n = Q((F_n \setminus \{0\}) \cup (F_n \setminus \{0\})^{-1}),
\]  
forms a F\o lner sequence in \( \PGL_2(\bK) \), thereby establishing the amenability of \( \PGL_2(\bK) \). To achieve this, it suffices to verify that  
\begin{equation} \label{limits}
\lim_{n \to \infty} \frac{|\widetilde{F}_n \cap u_{+}(b)\widetilde{F}_n|}{|\widetilde{F}_n|}
= 
\lim_{n \to \infty} \frac{|\widetilde{F}_n \cap u_{-}(b)\widetilde{F}_n|}{|\widetilde{F}_n|} = 1, \quad \text{for all } b \in \bK \setminus \{0\}.
\end{equation}
To simplify notation, let 
\[
A_n = (F_n \setminus \{0\}) \cup (F_n \setminus \{0\})^{-1}, \quad n \geq 1,
\]
and note that since both $(F_n)$ and $((F_n \setminus \{0\})^{-1})$ are 
F\o lner sequences in $(\bK,+)$, the sequence $(A_n)$ is F\o lner too, and satisfies 
$A_n = A_n^{-1}$ for all $n$. \\

To prove the first limit in \eqref{limits}, we observe that  
\[
\widetilde{F}_n \cap u_{+}(b)\widetilde{F}_n
= Q(A_n \cap (A_n + b)), 
\]
and thus
\[
\frac{|\widetilde{F}_n \cap u_{+}(b)\widetilde{F}_n|}{|\widetilde{F}_n|}
= \frac{|A_n \cap (A_n + b)| \cdot (|A_n \cap (A_n + b)| - 1) \cdot (|A_n \cap (A_n + b)| - 2)|}{|A_n| \cdot (|A_n| - 1) \cdot (|A_n|-2)}, 
\]
which clearly tends to $1$ as $n \ra \infty$ for all $b \in \bK$ since $(A_n)$ is a 
F\o lner sequence.  \\

For the second limit in \eqref{limits}, we first note that  
\[
|\widetilde{F}_n \cap u_{-}(b)\widetilde{F}_n| \geq 
\left|Q\left(A_n \cap \varphi_b\left(A_n \setminus \left\{ -\frac{1}{b} \right\}\right)\right)\right|.
\]
Since \( A_n^{-1} = A_n \), we obtain by the previous lemma: 
\[
\left|A_n \cap \varphi_b\left(A_n \setminus \big\{ -\frac{1}{b} \big\}\right)\right| = |A_n \cap (A_n + b)|,
\]
and thus
\[
1 \geq \frac{|\widetilde{F}_n \cap u_{-}(b)\widetilde{F}_n|}{|\widetilde{F}_n|}
\geq \frac{|Q(A_n \cap (A_n + b))|}{|A_n| \cdot (|A_n|-1) \cdot (|A_n| - 2)}, \quad \textrm{for all $n$}.
\]
Since \( (A_n) \) is a F\o lner sequence in \( (\bK,+) \), the right-hand side tends to $1$ as $n \ra \infty$, completing the proof.

\section{General F\o lner-Kloosterman sums: Proof of Theorem \ref{Thm_Equin1nrMain}}

\subsection{Untwisted case}

We begin the proof of Theorem \ref{Thm_Equin1nrMain} by first establishing the untwisted case. In the following subsection, we will then reduce the general case to this special case.

\begin{theorem}
\label{Thm_Equin1nr}
Let $n_1,\ldots,n_r$ be non-zero integers, not divisible by the characteristics of $\bK$. Suppose that the set $I = \{ i=1,\ldots,r \, \mid \, n_i > 0\}$ is non-empty. Then, for every double F\o lner sequence $(F_k)$ in $\bK$ and 
$\xi_1,\ldots,\xi_r \in \smash{\widehat{\bK}}$ satisfying $\xi_i \neq 1$ for all $i \in I$, we have
\[
\lim_{k \ra \infty} \frac{1}{|F_k|} \sum_{a \in F_k \setminus \{0\}} \xi_1(a^{n_1}) \cdots \xi_r(a^{n_r}) = 0.
\]
\end{theorem}

The proof is very similar to the proof of Theorem \ref{Thm_main2}, but we present it here for completeness.

\begin{proof}
Fix \( \xi_1, \dots, \xi_r \in \widehat{\mathbb{K}} \setminus \{1\} \), and let \( (F_k) \) be a double F\o lner sequence in \( (\mathbb{K}, +) \). Consider the sequence \( (\theta^{(k)}) \) of sub-probability measures on \( \smash{\widehat{\mathbb{K}}}^r  \), defined by  
\[
\theta^{(k)} = \frac{1}{|F_k|} \sum_{a \in F_k \setminus \{0\}} \delta_{(a^{n_1})^* \xi_1} \otimes \cdots \otimes \delta_{(a^{n_r})^* \xi_r}, \quad k \geq 1.
\]
Let \( \smash{(\theta^{(k_j)})} \) be a weak*-convergent subsequence with limit \( \theta \). Since \( (F_k) \) is a double F\o lner sequence, it is asymptotically invariant under \( \mathbb{K}^* \), implying that \( \theta \) is a probability measure on \( \smash{\widehat{\mathbb{K}}^r} \) that is invariant under the subgroup  
\[
\Gamma = \{ (a^{n_1},\dots,a^{n_r}) \mid a \in \mathbb{K}^* \}.
\]
By Theorem \ref{Thm_RelPrime}, this subgroup is coupling-resistant. Moreover, by Lemma \ref{Lemma_KnIsMeasureRigid}, each coordinate projection \(\Gamma_i < \mathbb{K}^*\) is measure-rigid. Applying Lemma \ref{Lemma_CouplingResistantEquivalence} (iii), it follows that there exists a unique probability measure \( p_\theta \) on \( \mathscr{L}(\mathbb{K}) \times \cdots \times \mathscr{L}(\mathbb{K}) \) such that  
\[
\theta = \sum_{(W_1,\ldots,W_r)} p_\theta(W_1,\ldots,W_r) m_{W_1^{\perp}} \otimes \cdots \otimes m_{W_r^{\perp}}.
\]
Note that each $W_i$ is either $\{0\}$ or $\bK$. Consequently,  
\begin{align*}
\lim_{j \to \infty} \frac{1}{|F_{k_j}|} \sum_{a \in F_{k_j} \setminus \{0\}} \xi_1(a^{n_1}) \cdots \xi_r(a^{n_r}) 
&= \widehat{\theta}(1,\dots,1) \\[0.2cm]
&= \sum_{(W_1,\ldots,W_r)} p_\theta(W_1,\ldots,W_r) \chi_{W_1}(1) \cdots \chi_{W_r}(1).
\end{align*}
Now, since the index set \( I \) is assumed to be nonempty, Lemma \ref{Lemma_KnIsMeasureRigid} implies that, for each \( i \in I \), the projection \( \theta_i \) onto the \( i \)th coordinate satisfies  
\[
\widehat{\theta}_i(b) = \lim_{k \to \infty} \frac{1}{|F_k|} \sum_{a \in F_k \setminus \{0\}} \xi_i(a^{n_i} b) = 0, \quad \text{for all } b \in \mathbb{K} \setminus \{0\}.
\]
Thus, \( \theta_i = m_{\widehat{\mathbb{K}}} \) for all \( i \in I \), ensuring that any \( r \)-tuple \( (W_1, \dots, W_r) \) in the support of \( p_\theta \) satisfies \( W_i = \{0\} \) for all \( i \in I \). Consequently,  
\[
\sum_{(W_1,\ldots,W_r)} p_\theta(W_1,\ldots,W_r) \chi_{W_1}(1) \cdots \chi_{W_r}(1) = 0,
\]
which establishes  
\[
\lim_{j \to \infty} \frac{1}{|F_{k_j}|} \sum_{a \in F_{k_j} \setminus \{0\}} \xi_1(a^{n_1}) \cdots \xi_r(a^{n_r}) = 0.
\]
Since the weak*-convergent subsequence \( (\theta^{(k_j)}) \) was chosen arbitrarily, the result follows.  
\end{proof}

\subsection{Proof of Theorem \ref{Thm_Equin1nrMain}}
\label{subsec:Equi1nrMain}
Let \(\bK\) be a countably infinite discrete field, and let \((F_k)\) be a double F\o lner sequence in \(\bK\). Consider nonzero integers \(n_1, \ldots, n_r\) such that none of them are divisible by the characteristic of \(\bK\). Assume that the set
\[
I = \{ i \in \{1, \ldots, r\} \mid n_i > 0 \}
\]
is nonempty. \\

Fix additive characters \(\xi_1, \ldots, \xi_r \in \widehat{\bK}\) with \(\xi_i \neq 1\) for all \(i \in I\), and fix a multiplicative character \(\eta \in \widehat{\bK^*}\). Define the function \(u : \bK^* \to \mathbb{C}\) by
\[
u(a) = \eta(a) \, \xi_1(a^{n_1}) \cdots \xi_r(a^{n_r}), \quad a \in \bK^*.
\]
Our goal is to prove that
\[
\lim_{k \to \infty} \frac{1}{|F_k|} \sum_{a \in F_k \setminus \{0\}} u(a) = 0.
\]
To this end, we will use the following special case of the classical van der Corput lemma (see \cite[Theorem 2.15]{BM}):

\begin{lemma} 
If \(u : \bK^* \to \mathbb{C}\) satisfies
\[
\lim_{k \to \infty} \frac{1}{|F_k|} \sum_{a \in F_k \setminus \{0\}} u(ah) \overline{u(a)} = 0
\]
for all but finitely many \(h \in \bK^*\), then
\[
\lim_{k \to \infty} \frac{1}{|F_k|} \sum_{a \in F_k \setminus \{0\}} u(a) = 0.
\]
\end{lemma}

Now, define the set
\[
D := \{ h \in \bK^* \mid h^{n_j} = 1 \text{ for some } j = 1, \ldots, r \}.
\]
Since each equation \(h^{n_j} = 1\) describes the roots of a polynomial of degree \(n_j\), the set \(D\) is a finite union of finite sets, hence finite. \\

By applying the lemma and using the multiplicativity of \(\eta\), it clearly suffices to show that for all \(h \notin D\),
\[
\lim_{k \to \infty} \frac{1}{|F_k|} \sum_{a \in F_k \setminus \{0\}} u(ah) \overline{u(a)} = 0.
\]
Note that
\[
u(ah) \overline{u(a)} = \eta(h) \cdot \prod_{j=1}^r \xi_{j,h}(a^{n_j}),
\]
where
\[
\xi_{j,h}(a) := \xi_j\big((h^{n_j} - 1) a \big), \quad j=1,\ldots,r, \quad a \in \bK^*.
\]
Since \(h \notin D\), none of the values \(h^{n_j}\) equals 1, so \(\xi_{i,h} \neq 1\) for every \(i \in I\). Therefore, the limit
\[
\lim_{k \to \infty} \frac{1}{|F_k|} \sum_{a \in F_k \setminus \{0\}} \xi_{1,h}(a^{n_1}) \cdots \xi_{r,h}(a^{n_r}) = 0
\]
follows directly from Theorem \ref{Thm_Equin1nr} in the untwisted case.

\section{Failure of equidistribution: Proof of Theorem \ref{Thm_WeirdInverseMain}}
\label{Sec:Failure}

For the proof, we will need the following simple lemma. We emphasize that $\widehat{\bR}$ denote the set of \emph{continuous} characters on $\bR$.

\begin{lemma}
\label{Lemma_RealSubfield}
For any pair \((\xi_1, \xi_2) \in \widehat{\bR}^2\),
\[
Z(\xi_1,\xi_2) := \overline{\{(\xi_1(t), \xi_2(1/t)) \, \mid \, t \in \bR^*\}} \neq \bT^2.
\]  
\end{lemma}

\begin{proof}
The lemma is trivial if either \(\xi_1\) or \(\xi_2\) is trivial. We therefore assume that both \(\xi_1\) and \(\xi_2\) are non-trivial, meaning they can be expressed as  
\[
\xi_1(t) = e^{2\pi i x_1 t}, \quad \xi_2(t) = e^{2\pi i x_2 t}, \quad t \in \mathbb{R},  
\]  
for unique \(x_1, x_2 \in \mathbb{R} \setminus \{0\}\). Our goal is to show that for each \(\lambda_1 = e^{2\pi i \theta} \in \mathbb{T} \setminus \{1\}\), 
where we can assume that $\theta \in (0,1)$, the closed set  
\[
Z(\xi_1,\xi_2)_{\lambda_1} = \{ \lambda_2 \in \mathbb{T} \mid (\lambda_1,\lambda_2) \in Z(\xi_1,\xi_2) \}
\]  
is countable and accumulates only at \(1\). This immediately implies that \(Z(\xi_1,\xi_2) \neq \mathbb{T}^2\). \\

To establish this, pick $\lambda_2 \in Z(\xi_1,\xi_2)_{\lambda_1}$ and choose a sequence \((t_n)\) in \(\mathbb{R}^*\) satisfying  
\[
\xi_1(t_n) \ra \lambda_1 \qand \xi_2(1/t_n) \ra \lambda_2,
\]
or equivalently, 
\[
e^{2\pi i x_1 t_n} \to e^{2\pi i \theta} \quad \text{and} \quad e^{2\pi i x_2/t_n} \to \lambda_2. 
\]  
The first limit implies that we can write  
\[
x_1 t_n = \theta_n + k_n,  
\]  
where \(\theta_n \to \theta\) for some sequence \((\theta_n)\) in $[0,1)$ and \((k_n)\) is a sequence of integers. Consequently,  
\[
e^{2\pi i x_2/t_n} = e^{2\pi i x_1 x_2 / (\theta_n + k_n)} \ra \lambda_2, \quad \textrm{as $n \ra \infty$}.
\]  
Since \(\theta \neq 0\), we conclude that either \(\lambda_2 = 1\) (if \(|k_n| \to \infty\)) or \(\lambda_2 = e^{2\pi i x_1 x_2 / (\theta + k)}\) for some \(k \in \mathbb{Z}\). In other words,
\[
Z(\xi_1,\xi_2)_{\lambda_1}
= \{1\} \cup \{ e^{2\pi i x_1 x_2 / (\theta + k)} \, \mid \, k \in \bZ \},
\]
completing the proof.  
\end{proof}

\subsection{Proof of Theorem \ref{Thm_WeirdInverseMain}}

In what follows, let $\bK \subset \bR$ be a countable sub-field. Clearly, $\bK$ 
is a dense subset of $\bR$. \\

Fix $\xi \in \widehat{\mathbb{R}} \setminus \{1\}$.  
By applying Lemma \ref{Lemma_RealSubfield} to $\xi_1 = \xi_2 = \xi$, we obtain open sets $U_1, U_2 \subset \mathbb{T}$ such that  
\begin{equation}
\label{xaainv}
\{(\xi(a), \xi(a^{-1})) \mid a \in \mathbb{K}^*\} \cap (U_1 \times U_2) = \emptyset.
\end{equation}  
Choose open sets $V_1, V_2, W_1, W_2 \subset \mathbb{T}$ satisfying  
\[
V_1 W_1^{-1} \subset U_1 \quad \text{and} \quad V_2 W_2^{-1} \subset U_2.
\]  
Then,  
\[
\int_{\mathbb{T}^2} \chi_{V_1}(\xi(a) \lambda_1) \chi_{V_2}(\xi(a^{-1}) \lambda_2) \chi_{W_1}(\lambda_1) \chi_{W_2}(\lambda_2) \, dm_{\mathbb{T}^2}(\lambda_1,\lambda_2) = 0, \quad \text{for all } a \in \mathbb{K}^*.
\]  
To see why, assume for contradiction that there exists $a \in \mathbb{K}^*$ along with $\lambda_1 \in W_1$ and $\lambda_2 \in W_2$ such that $\xi(a) \lambda_1 \in V_1$ and $\xi(a^{-1}) \lambda_2 \in V_2$. This would imply  
\[
(\xi(a), \xi(a^{-1})) \in V_1 W_1^{-1} \times V_2 W_2^{-1},
\]  
contradicting \eqref{xaainv}. \\

On the other hand, using the $L^2$-Fourier expansions  
\[
\chi_{V_i}(\lambda_i) = \sum_{m_i \in \mathbb{Z}} \widehat{\chi}_{V_i}(m_i) \lambda_i^{m_i}  
\quad \text{and} \quad  
\chi_{W_i}(\lambda_i) = \sum_{n_i \in \mathbb{Z}} \widehat{\chi}_{W_i}(n_i) \lambda_i^{n_i},
\]  
for $i=1,2$, we obtain  
\begin{align*}
&\int_{\mathbb{T}^2} \chi_{V_1}(\xi(a) \lambda_1) \chi_{V_2}(\xi(a^{-1}) \lambda_2) \chi_{W_1}(\lambda_1) \chi_{W_2}(\lambda_2) \, dm_{\mathbb{T}^2}(\lambda_1,\lambda_2) \\[0.2cm]
&= 
\int_{\mathbb{T}^2} \chi_{V_1}(\xi(a) \lambda_1) \chi_{V_2}(\xi(a^{-1}) \lambda_2) \overline{\chi_{W_1}(\lambda_1)} \overline{\chi_{W_2}(\lambda_2)} \, dm_{\mathbb{T}^2}(\lambda_1,\lambda_2) \\[0.2cm] 
&= \sum_{(m_1, m_2) \in \mathbb{Z}^2} \widehat{\chi}_{V_1}(m_1) \overline{\widehat{\chi}_{W_1}(m_1)} \widehat{\chi}_{V_2}(m_2) \overline{\widehat{\chi}_{W_2}(m_2)} \xi(a)^{m_1} \xi(a^{-1})^{m_2} \\[0.2cm]
&=\sum_{(m_1, m_2) \in \mathbb{Z}^2} c(m_1, m_2) \xi(m_1 a) \xi(m_2 a^{-1}) = 0,
\end{align*}  
for all $a \in \mathbb{K}^*$, where  
\[
c(m_1, m_2) = \widehat{\chi}_{V_1}(m_1) \overline{\widehat{\chi}_{W_1}(m_1)} \widehat{\chi}_{V_2}(m_2) \overline{\widehat{\chi}_{W_2}(m_2)}, \quad (m_1, m_2) \in \mathbb{Z}^2,
\]  
and where we in the second to last identity have used that $\xi$ is an additive character on $\bR$, and
thus $\xi(b)^m = \xi(mb)$ for all $m \in \mathbb{Z}$ and $b \in \mathbb{K}$. Moreover, Parseval's formula ensures that the map $(m_1,m_2) \mapsto c(m_1, m_2)$ is absolutely summable on $\bZ^2$. \\ 

Now, let $(F_n)$ be a double F\o lner sequence for $(\mathbb{K},+)$. The above identity implies that  
\[
\sum_{(m_1, m_2) \in \mathbb{Z}^2} c(m_1, m_2) \left( \frac{1}{|F_n|} \sum_{a \in F_n} \xi(a)^{m_1} \xi(a^{-1})^{m_2} \right) = 0,  
\]  
for all $n$. We decompose this sum into four terms:  
\begin{align*}
A_n &= m_{\mathbb{T}}(V_1) m_{\mathbb{T}}(V_2) m_{\mathbb{T}}(W_1) m_{\mathbb{T}}(W_2), \\[0.2cm]
B_n &= \sum_{m_1, m_2 \neq 0} c(m_1, m_2) \left( \frac{1}{|F_n|} \sum_{a \in F_n} \xi(m_1 a) \xi(m_2 a^{-1}) \right), \\[0.2cm]
C_n &= \sum_{m_1 \neq 0} c(m_1, 0) \left( \frac{1}{|F_n|} \sum_{a \in F_n} \xi(m_1 a) \xi(a^{-1})^{m_2} \right), \\[0.2cm]
D_n &= \sum_{m_2 \neq 0} c(0, m_2) \left( \frac{1}{|F_n|} \sum_{a \in F_n} \xi(m_2 a^{-1}) \right).
\end{align*}  
By Theorem \ref{Thm_Equin1nr}, we have $B_n, C_n \to 0$ as $n \to \infty$, while $A_n$ is independent of $n$ and strictly positive. Since $A_n + B_n + C_n + D_n = 0$, we conclude that $(D_n)$ cannot converge to zero. Since $m_2 \mapsto c(0, m_2)$ is absolutely summable, this implies that at least one nonzero $m_2$ must satisfy  
\[
\frac{1}{|F_n|} \sum_{a \in F_n} \xi(m_2 a^{-1}) \not\to 0.
\]  
Thus, using the asymptotic invariance of $(F_n)$ under multiplication by $1/m_2$, we conclude that  
\[
\frac{1}{|F_n|} \sum_{a \in F_n} \xi(a^{-1}) \not\to 0,
\]  
completing the proof.  

\section{Orbital Poincar\'e pairs: Proofs of Theorem \ref{Thm_CombMain} and Theorem \ref{Thm_FSthm}}

Let $A$ be countable discrete abelian group. 

\begin{definition}
A subset \( \mathscr{S} \subset A \setminus \{0\}  \) is called a \textbf{Poincar\'e set} if, for every probability measure-preserving action \( A \curvearrowright (X, \mu) \) and every measurable set \( B \subset X \) with \( \mu(B) > 0 \), there exists \( s \in \mathscr{S} \) such that $\mu(B \cap sB) > 0$.
\end{definition}

\begin{remark}
Poincar\'{e}'s Recurrence Theorem (or more precisely, its standard proof) implies
that if $S$ is an \emph{infinite} subset of $A$, then $\mathscr{S} = (S-S)\setminus \{0\}$ is a Poincar\'{e} set.
\end{remark}

The following lemma is easy.

\begin{lemma}
\label{Lemma_s1sk}
Let $k \geq 1$ and let $\mathscr{S}_1,\ldots,\mathscr{S}_k \subset A$ be Poincar\'e sets. Then, for every probablity measure-preserving action $A \acts (X,\mu)$ and every measurable set $B \subset X$ with $\mu(B) > 0$, there exist elements $s_1 \in \mathscr{S}_1,\ldots,s_k \in \mathscr{S}_k$ such that
\[
\mu(B \cap s_1B \cap \cdots \cap s_k B) > 0.
\]
\end{lemma}

\begin{proof}
The case $k=1$ holds by definition, so let us assume $k \geq 2$. Since $\mathscr{S}_1$ is a Poincar\'{e} set, we can find $s_1 \in \mathscr{S}$ such that $B_1 := B \cap s_1 B$
satisfies $\mu(B_1) > 0$. We can now use that $\mathscr{S}_2$ is a Poincar\'{e} set
to find $s_2 \in \mathscr{S}_2$ such that $B_2 := B_1 \cap s_2 B_2 \subset B \cap s_1 B \cap s_2 B$ has positive measure. Continuing this process iteratively, we construct a sequence of elements \( s_1 \in \mathscr{S}_1, \dots, s_k \in \mathscr{S}_k \) such that  
\[
B_k := B_{k-1} \cap s_k B_{k-1} \subset B \cap s_1 B \cap \dots \cap s_k B
\]
has positive measure. This completes the proof.
\end{proof}

\begin{definition}
Let $\Gamma$ be a subgroup of the automorphism group $\Aut(A)$. We say that $(A,\Gamma)$ is an \textbf{orbital Poincar\'{e} pair} if for every $a \in A \setminus \{0\}$, the orbit
\[
\cO_a = \{ \gamma(a) \, \mid \, \gamma \in \Gamma \} \subseteq A \setminus \{0\}
\]
is a Poincar\'{e} set.
\end{definition}

The following lemma provides a criterion for when $(A,\Gamma)$ is an orbital Poincar\'{e} pair. Recall that if $\gamma \in \Aut(A)$, then $\gamma^* \in \Aut(\widehat{A})$ is defined by
\[
(\gamma^*. \,\xi)(a) = \xi(\gamma(a)), \quad \xi \in \widehat{A}, \enskip a \in A.
\]
We emphasize that the map  
\[
\Gamma \times \widehat{A} \to \widehat{A}, \quad (\gamma, \xi) \mapsto \gamma^* \cdot \xi
\]
defines a \emph{right} action, and only a left-action if $\Gamma$ is abelian. In most parts of this paper, we assume $\Gamma$ is abelian, so this distinction is not significant. However, the concept of a probability measure invariant under a right-action remains well-defined. We denote by $\Prob_\Gamma(\widehat{A})$ the space of probability measures on $\widehat{A}$ that are invariant under the right-action of $\Gamma$. Since $A$ is countable and discrete, $\widehat{A}$ is compact and metrizable, and the Haar probability measure $m_A$ is always $\Gamma$-invariant.

\begin{lemma}
\label{Lemma_PosFourier}
Suppose $\Gamma$ is amenable. If every \( \mu \in \Prob_\Gamma(\widehat{A}) \) has a non-negative Fourier transform, then \( (A, \Gamma) \) forms an orbital Poincar\'{e} pair.
\end{lemma}

\begin{remark}
The lemma extends to a general \(\Gamma\) if \(\Gamma\)-invariance is replaced by stationarity with respect to a fixed random walk on \(\Gamma\). For similar results, we refer the reader to \cite{BB,BF1}, which also address cases where the Fourier transform is not necessarily non-negative. Instead, a (potentially negative) lower bound is imposed on the Fourier transform of a stationary probability measure on \(\smash{\widehat{A}}\), which can be made arbitrarily close to zero by restricting to a sufficiently small finite-index subgroup of \(A\).
\end{remark}

\begin{proof}
Let \( A \acts (X,\mu) \) be a probability measure-preserving action, and let \( B \subset X \) be a measurable subset with positive \(\mu\)-measure. By Bochner's theorem, there exists a unique non-negative Borel measure \( \sigma_B \) on \( \widehat{A} \) such that  
\[
\mu(B \cap aB) = \int_{\widehat{A}} \xi(a) \, d\sigma_B(\xi), \quad \text{for all } a \in A,
\]
with \( \sigma_B(\{1\}) > 0 \). This implies that \( \sigma_B \) can be decomposed as  
\[
\sigma_B = \sigma_B(\{1\}) \delta_1 + \sigma^o_B,
\]
where \( \sigma^o_B \) is a finite non-negative Borel measure satisfying \( \sigma^o_B(\{1\}) = 0 \). Now, fix \( a \in A \setminus \{1\} \) and a F\o lner sequence \( (F_n) \) in \( \Gamma \). Our goal is to show that  
\[
\limsup_{n \to \infty} \frac{1}{|F_n|} \sum_{\gamma \in F_n} \mu(B \cap \gamma(a)B) > 0.
\]
In particular, this implies that there exists at least one $\gamma \in \Gamma$ such that $\mu(B \cap \gamma(a)B) > 0$. \\

To this end, we observe that  
\begin{align*}
\mu(B \cap \gamma(a)B) 
&= \int_{\widehat{A}} \xi(\gamma(a)) \, d\sigma_B(\xi)
=  \int_{\widehat{A}} (\gamma^*. \, \xi)(a) \, d\sigma_B(\xi) = \widehat{(\gamma^*)_*\sigma}_B(a) \\[0.2cm]
&= \sigma_B(\{1\}) + \widehat{(\gamma^*)_*\sigma^o_B}(a),
\end{align*}
and thus
\[
\frac{1}{|F_n|} \sum_{\gamma \in F_n} \mu(B \cap \gamma(a)B) = 
\sigma_B(\{1\}) + \frac{1}{|F_n|} \sum_{\gamma \in F_n} \widehat{(\gamma^*)_*\sigma^o_B}(a).
\]
Consider the sequence $(\tau_n)$ of bounded non-negative Borel measures on $\widehat{A}$ defined by
\[
\tau_n = \frac{1}{|F_n|} \sum_{\gamma \in F_n} (\gamma^*)_*\sigma^o_B.
\]
We select a weak*-convergent subsequence with weak*-limit \( \tau_\infty \), which then must be $\Gamma$-invariant, and note that
\[
\limsup_{n \to \infty} \frac{1}{|F_n|} \sum_{\gamma \in F_n} \mu(B \cap \gamma(a)B)
\geq \sigma_B(\{1\}) + \widehat{\tau}_\infty(a).
\]
Since \( \tau_\infty \) is a positive scaling of a \(\Gamma\)-invariant probability measure on \( \widehat{A} \), our assumption imply that $\widehat{\tau}_\infty(a) \geq 0$ for all $a \in A$, and the result follows.  
\end{proof}

Using Furstenberg's Correspondence Principle, we can now deduce the following combinatorial consequence in the presence of orbital Poincar\'e pairs. We recall
that the \textbf{upper Banach density} $d^*(E)$ of a set $E \subset A$ is defined by 
\[
d^*(E) = \sup_{(F_n)} \left( \varlimsup_{n \ra \infty} \frac{|E \cap F_n|}{|F_n|} \right),
\]
where the supremum is taken over all F\o lner sequences $(F_n)$ in $A$.

\begin{theorem}
\label{Thm_Apsi}
Let $(A,\Gamma)$ be an orbital Poincar\'{e} pair and let $Y$ be a set. Suppose 
$\psi : A \ra Y$ is a $\Gamma$-invariant function, and let $E \subset A$ be a subset
with positive upper Banach density. Then, for every finite set $F \subset A$, 
there exists $a_o = a_o(F) \in E$ such that 
\[
\psi(F) \subset \psi(E-a_o).
\]
In particular, $\psi(A) = \psi(E-E)$.
\end{theorem}

\begin{proof}
Let \( E \subset A \) be a set with positive upper Banach density. By Furstenberg's Correspondence Principle, there exists a probability measure-preserving action \( A \acts (X, \mu) \) and a measurable subset \( B \subset X \) with positive measure such that  
\[
d^*(E \cap (E - a_1) \cap \cdots \cap (E - a_k)) \geq \mu(B \cap a_1B \cap \cdots \cap a_k B),
\]
for all \( a_1, \ldots, a_k \in A \). \\

Let us now fix a finite set \( F = \{a_1, \ldots, a_k\} \) in \( A \setminus \{0\} \). Since \( (A, \Gamma) \) is an orbital Poincar\'{e} pair, the orbits \( \mathcal{O}_{a_1}, \ldots, \mathcal{O}_{a_k} \) are  Poincar\'{e} sets in \( A \). Hence, by Lemma \ref{Lemma_s1sk}, there exist \( \gamma_1, \ldots, \gamma_k \in \Gamma \) such that  
\[
\mu(B \cap \gamma_1(a_1)B \cap \ldots \cap \gamma_k(a_k)B) > 0,
\]
and consequently,  
\[
E \cap (E - \gamma_1(a_1)) \cap \ldots \cap (E - \gamma_k(a_k)) \neq \emptyset.
\]
Choose any element \( a_0 \) from this intersection. Then, we have  
\[
a_0 + \gamma_j(a_j) \in E, \quad \text{for all } j = 1, \ldots, k.
\]
Thus, for any \( \Gamma \)-invariant function \( \psi : A \to Y \),  
\[
\psi(a_j) = \psi(\gamma_j(a_j)) \in \psi(E - a_0), \quad \text{for all } j = 1, \ldots, k,
\]
so $\psi(F) \subset \psi(E-a_0)$. Finally, since $a_0 \in E$, we have $0 \in E-a_0$, so this inclusion also works if we assume $0 \in F$. 
This completes the proof.
\end{proof}

\subsection{Orbital Poincar\'e sets from diagonalizable representations}

Let $\bK$ be a countably infinite field and let $U$ be a countable vector space over $\bK$. Let $\rho_U : \bK^* \ra \GL(U)$ be a representation and suppose $U$ admits a direct sum decomposition of the form
\[
U = \bigoplus_{m \in I} U_m,
\]
for some (possibly infinite) set $I \subset \bZ \setminus \{0\}$ such that every $m \in I$, 
\[
\rho(a)u = a^m u, \quad \textrm{for all $u \in U_m$}.
\]
\begin{lemma}
\label{Lemma_Urho}
$(U,\rho_U)$ is an orbital Poincar\'e pair.
\end{lemma}

\begin{proof}
Our goal is to show that for every nonzero vector \( u \in U \setminus \{0\} \), the orbit
\[
\mathcal{O}_u = \{ \rho(a)u : a \in \bK^* \} \subset U
\]
is a Poincar\'e set. Given any fixed \( u \in U \), there exists a \emph{finite} subset \( J \subset I \) such that
\[
u = \sum_{m \in J} u_m, \quad \text{with } u_m \in U_m \text{ for each } m \in J.
\]
Therefore, without loss of generality, we may restrict the representation \(\rho\) to the finite-dimensional subspace of \( U \) spanned by the weight spaces \(\{ U_m : m \in J \}\). After this restriction, we may assume that the pair \((U, \rho)\) is diagonalizable. \\

To prove that the orbit \(\mathcal{O}_u\) is indeed a Poincar\'e set, Lemma \ref{Lemma_PosFourier} reduces the problem to showing that every \(\Gamma\)-invariant probability measure on the dual group \(\smash{\widehat{U}}\) has a non-negative Fourier transform. Since the ergodic decomposition of \(\Gamma\)-invariant probability measures commutes with the Fourier transform, it suffices to verify this property for \(\Gamma\)-invariant and \(\Gamma\)-ergodic probability measures on \(\smash{\widehat{U}}\). \\

By Theorems \ref{Thm_ReducedModels} and \ref{Thm_DiagonalMain}, each such measure \(\theta\) can be expressed as
\[
\theta = \Psi^*\big(m_{W_1^\perp} \otimes \cdots \otimes m_{W_s^\perp}\big),
\]
for some linear subspaces \( W_1 \leq V_1, \ldots, W_s \leq V_s \), where \(\Psi = \Phi^{-1}\) corresponds to a model \((V, \rho_V, \Phi)\) of \((U, \rho_U)\) with weight spaces \(V_1, \ldots, V_s\). \\

Since the measure \( m_{W_1^\perp} \otimes \cdots \otimes m_{W_s^\perp} \) has a Fourier transform taking values only in \(\{0,1\}\), and the pullback \(\smash{\Psi^*}\) preserves the range of Fourier transforms, the claim follows.
\end{proof}

\subsection{Proofs of Theorem \ref{Thm_CombMain} and Theorem \ref{Thm_FSthm}}
\label{subsec:FSthm}

Theorem \ref{Thm_CombMain} follows directly from Theorem \ref{Thm_Apsi} and Lemma \ref{Lemma_Urho}. \\

To prove Theorem \ref{Thm_FSthm}, let $\bK$ be a countably infinite field, and consider the countable vector space 
\[
U \subset \bK[X, X^{-1}]
\]
consisting of all Laurent polynomials of the form
\[
p(X) = \sum_{k=-\infty}^\infty p_k X^k,
\]
where the coefficients $(p_k)$ have finite support and satisfy $p_0 = 0$. Define a representation
\[
\rho_U : \bK^* \to \mathrm{GL}(U)
\]
by
\[
(\rho_U(a) p)(X) = p(aX) = \sum_{k=-\infty}^\infty a^k p_k X^k, \quad a \in \bK^*, \quad p \in U.
\]
For each $k \neq 0$, the one-dimensional subspace 
$U_k = \mathrm{span}_\bK \{ X^k \} \subset U$
is a weight space for $\rho_U$ with weight $k$. Consequently,
\[
U = \bigoplus_{k \neq 0} U_k,
\]
so the pair $(U, \rho_U)$ satisfies the assumptions of Lemma \ref{Lemma_Urho}. It follows that for every $p \in U$, the orbit
\[
\mathcal{O}_p = \{ \rho_U(a) p : a \in \bK^* \} \subset U
\]
forms a Poincar\'e set with respect to every probability measure-preserving action of the additive group $(U, +)$. \\

Now, let $E \subset \bK$ be a set with positive upper Banach density. By Furstenberg's Correspondence Principle, there exists an ergodic $\bK$-space $(X, \mu)$ and a measurable set $B \subset X$ with $\mu(B) > 0$ such that
\[
d^*(E \cap (E - b)) \geq \mu(B \cap b.B) \quad \text{for all } b \in \bK.
\]
In particular, for every Laurent polynomial $p \in U$,
\[
d^*(E \cap (E - p(a))) \geq \mu(B \cap p(a).B) \quad \text{for all } a \in \bK^*.
\]
The crucial observation is that the evaluation map
\[
\varphi : U \to \bK, \quad p \mapsto p(1)
\]
is a group homomorphism with respect to addition. Using $\varphi$, we can view $(X, \mu)$ as an ergodic $U$-space via $\varphi$.
Thus,
\[
\mu(B \cap p(a).B) = \mu(B \cap \varphi(\rho_U(a) p).B) \quad \text{for all } a \in \bK^*, p \in U.
\]
Since $\mathcal{O}_p$ is a Poincar\'e set for this $U$-action, there exists some $a \in \bK^*$ such that
\[
0 < \mu(B \cap \varphi(\rho_U(a) p).B) = \mu(B \cap p(a).B) \leq d^*(E \cap (E - p(a))).
\]
Hence, $E \cap (E - p(a))$ is nonempty, implying that $p(a) \in E - E$, 
as claimed.

\subsection{Hyperbolas do not contain infinite difference sets}
\label{subsec:hyperbola}

Let \(\bK\) be a countably infinite field. By Lemma \ref{Lemma_Urho}, for every \(t \neq 0\), the \textbf{hyperbola}
\[
\mathcal{H}_t = \left\{ \left(a, \frac{t}{a}\right) : a \in \bK^* \right\} \subset \bK^2
\]
is a Poincar\'e set. \\

If the set \(\cH_t\) were to contain the difference set of an infinite subset of \(\bK^2\), the result would follow immediately from Poincar\'e's classical recurrence theorem. However, as we now show, this is very far from the case: hyperbolas cannot contain the difference set of any subset with more than $3$ elements.

\begin{lemma}
\label{Lemma_Htnotcontain}
Suppose $\Char(\bK) \neq 2$ and let $F \subset \bK^2$ be a set with at least $4$ elements. Then,
\[
(F - F) \setminus \{0\} \nsubseteq \cH_t, \quad \textrm{for every $t \in \bK^*$}.
\]
\end{lemma}

\begin{remark}
The assumption $\operatorname{char}(\bK)\neq 2$ is essential for the lemma below; however, we have not investigated whether it is also necessary for Lemma~\ref{Lemma_Htnotcontain}.
\end{remark}

In the proof we will use the following simple result.

\begin{lemma}
\label{Lemma_stupid}
Suppose $\Char(\bK) \neq 2$. Then there does not exist $x,y,z \in \bK^*$ such that 
\[
x + y \neq 0, \quad y + z \neq 0, \quad x + y + z \neq 0,
\]
satisfying
\[
\left\{
\begin{array}{ll}
\frac{1}{x+y} &= \frac{1}{x} + \frac{1}{y} \\[0.2cm]
\frac{1}{y+z} &= \frac{1}{y} + \frac{1}{z} 
\end{array}
\right.
\qand
\frac{1}{x+y+z} = \frac{1}{x} + \frac{1}{y} + \frac{1}{z}. \\
\]
\end{lemma}

\begin{remark}
Observe that if $\Char(\bK) = 2$ and $\bK$ contains a primitive third root of unity $\alpha_0$, then for every $s \in \bK^*$, the triple
\[
(x, y, z) = (s\alpha_0,\, s,\, s\alpha_0)
\]
is a solution to the system of equations in the lemma. In particular, the system admits infinitely many distinct solutions.
\end{remark}

\begin{proof}[Proof of Lemma \ref{Lemma_stupid}]
Suppose \(x, y, z \in \bK^*\) are elements satisfying the conditions of the lemma. The two given equations on the left imply that
\[
(x + y)^2 = xy, \quad (y + z)^2 = yz.
\]
Since \(y\) and \(z\) are nonzero, define
\[
u := \frac{x}{y}, \quad v := \frac{y}{z}.
\]
Each of these satisfies the quadratic equation
\[
\alpha^2 + \alpha + 1 = 0.
\]
If this equation has no solution in \(\bK\), then the claim follows immediately. Otherwise, let \(\alpha_o \in \bK\) be a solution in $\bK$; note that \(1/\alpha_o\) is also a solution (in particular, $\alpha_o \neq 0$). We must then consider the four possible cases:
\begin{equation}
\label{uvalpha}
u = \alpha_0^{\pm 1}, \quad v = \alpha_0^{\pm 1}.
\end{equation}
Rewriting the equation
\[
\frac{1}{x + y + z} = \frac{1}{x} + \frac{1}{y} + \frac{1}{z}
\]
in terms of \(u\) and \(v\), we get
\[
\frac{1}{1 + v + uv} = 1 + \frac{1}{v} + \frac{1}{uv}.
\]
The exponents for \(u\) and \(v\) in \eqref{uvalpha} must agree; otherwise, we would have 
\[
\frac{x}{z} = \frac{x}{y} \cdot \frac{y}{z} = uv = 1,
\]
and thus
\[
\frac{1}{2+v} = 2 + \frac{1}{v},
\]
leading to the equation
\[
v = (2+v)(2v+1) = 4v + 2 + 2v^2 + v,
\]
and thus
\[
2\underbrace{(v^2 + v + 1)}_{=0} + 2v = 0,
\]
which is impossible, since $\Char(\bK) \neq 2$ and $v \neq 0$. Hence, only two cases remain:
\[
u = v = \alpha_0 \quad \text{or} \quad u = v = \frac{1}{\alpha_o}.
\]
Substituting the above cases into
\[
\frac{1}{1 + v + uv} = 1 + \frac{1}{v} + \frac{1}{uv}
\]
and simplifying, they both reduce to the same condition:
\[
0 = (1 + \alpha_o + \alpha_o^2)^2 = \alpha_0^2,
\]
forcing $\alpha_o = 0$, which is clearly a contradiction, and the proof is complete.
\end{proof}

\begin{proof}[Proof of Lemma \ref{Lemma_Htnotcontain}]
Let \( F = \{v_1, v_2, v_3, v_4\} \) be a set of four elements. Suppose, for the sake of contradiction, that for some \( t \in \bK^* \), we have the inclusion
\[
F - F \subset \cH_t.
\]
Then there exist \( a_1, a_2, a_3 \in \bK^* \) such that
\[
v_2 - v_1 = \left(a_1, \frac{t}{a_1}\right), \qquad
v_3 - v_2 = \left(a_2, \frac{t}{a_2}\right), \qquad
v_4 - v_3 = \left(a_3, \frac{t}{a_3}\right),
\]
and \( b_1, b_2, b_3 \in \bK^* \) such that
\[
v_3 - v_1 = \left(b_1, \frac{t}{b_1}\right), \qquad
v_4 - v_2 = \left(b_2, \frac{t}{b_2}\right), \qquad
v_4 - v_1 = \left(b_3, \frac{t}{b_3}\right).
\]
We have:
\begin{align*}
v_3 - v_1 &= (v_3 - v_2) + (v_2 - v_1), \\[0.2cm]
v_4 - v_2 &= (v_4 - v_3) + (v_3 - v_2), \\[0.2cm]
v_4 - v_1 &= (v_4 - v_3) + (v_3 - v_2) + (v_2 - v_1).
\end{align*}
Substituting the expressions for these differences yields:
\[
b_1 = a_1 + a_2, \qquad
b_2 = a_2 + a_3, \qquad
b_3 = a_1 + a_2 + a_3.
\]
Since \( b_1, b_2, b_3 \in \bK^* \), we must have:
\begin{equation}
\label{aisnotnull}
a_1 + a_2 \neq 0, \qquad a_2 + a_3 \neq 0, \qquad a_1 + a_2 + a_3 \neq 0.
\end{equation}
On the other hand, from the second coordinates, we obtain:
\[
\frac{1}{a_1 + a_2} = \frac{1}{a_1} + \frac{1}{a_2}, \qquad
\frac{1}{a_2 + a_3} = \frac{1}{a_2} + \frac{1}{a_3}, \qquad
\frac{1}{a_1 + a_2 + a_3} = \frac{1}{a_1} + \frac{1}{a_2} + \frac{1}{a_3}.
\]
However, Lemma~\ref{Lemma_stupid} asserts that no triple \( (a_1, a_2, a_3) \in (\bK^*)^3 \) satisfies \eqref{aisnotnull} and all three of these identities simultaneously. This contradiction completes the proof.
\end{proof}

\newpage

\appendix

\section{Field embedddings from embeddings of multiplicative groups}
\label{App:FieldHomo}

Our aim in this appendix is to prove the following theorem.

\begin{theorem}
\label{Thm_Fields}
Let $\bK$ and $\bL$ be countably infinite fields and suppose there exists an injective 
homomorphism $\rho : \bK^* \ra \bL^*$ along with finite-valued maps $u,v : \bK^* \ra \bL^*$ satisfying
\[
\rho(1 + x) = u(x) + v(x) \rho(x), \quad \textrm{for all $x \in \bK^*$},
\]
where we adopt the convention that $\rho(0) = 0$. Then,
\begin{itemize}
\item[$(i)$] The map $w = \frac{v}{u} : \bK^* \ra \bL^*$ is a homomorphism with finite image. \vspace{0.2cm}
\item[$(ii)$] The map $\kappa : \bK \ra \bL$ defined by $\kappa(0) = 0$ and $\kappa(x) = w(x) \rho(x)$ for $x \in \bK^*$ is an injective field homomorphism.
\end{itemize}
Furthermore, if $\rho$ is a bijection, then $\kappa$ is a field isomorphism between $\bK$ and $\bL$.
\end{theorem}

The proof is lengthy and technical, so we will break it down into several steps.

\subsection{Patchwise homomorphisms}

The core of the proof is to show that the maps \( w \) and \( \kappa \) from Theorem \ref{Thm_Fields} act as homomorphisms on large "patches" of the groups \((\bK^*,\cdot)\) and \((\bK,+)\), respectively. The following auxiliary definition and Lemma \ref{Lemma_Patchwise} capture the essential properties needed.

\begin{definition}
\label{Def_Patchwise}
Let \( G \) and \( H \) be groups. A map \( \eta: G \to H \) is called a \emph{patchwise homomorphism} if the following conditions hold:  
\vspace{0.2cm}
\begin{enumerate}
    \item There exists a finite set \( A \subset G \) such that for all \( x \notin A \),  
    \vspace{0.2cm}
    \begin{itemize}
        \item \( \eta(x^{-1}) = \eta(x)^{-1} \), and  \vspace{0.2cm}
        \item there exists a finite set \( B_x \subset G \) such that  
        \[
        \eta(xy) = \eta(x) \eta(y), \quad \text{for all } y \notin B_x.
        \]  
    \end{itemize}

    \item For every \( (x,y) \in G \times G \), there exists a finite set \( C_{x,y} \subset G \) such that  
    \[
    zy \notin B_{xz^{-1}}, \quad \text{for all } z \notin C_{x,y}.
    \]
\end{enumerate}
\end{definition}

\begin{lemma}
\label{Lemma_Patchwise}
Let $G$ and $H$ be groups. If $G$ is infinite, then every patchwise homomorphism from $G$ to $H$ is a homomorphism.
\end{lemma}

\begin{proof}
Let \( \eta: G \to H \) be a patchwise homomorphism, and suppose \( G \) is infinite. Pick \( x, y \in G \). We aim to show  
\begin{equation}
\label{xyy-1}
\eta(x) = \eta(xy) \eta(y^{-1}).
\end{equation}
Because \(x\) and \(y\) are arbitrary, substituting \(x=st\) and \(y=t^{-1}\) yields
\[
\eta(st)=\eta(s)\eta(t), \quad \textrm{for all $s,t\in G$},
\]
and hence \(\eta\) is a homomorphism. To prove \eqref{xyy-1}, note that by Condition (2), there exists a finite set \( C_{x,y} \subset G \) such that \( zy \notin B_{xz^{-1}} \) for all \( z \notin C_{x,y} \), and thus, by Condition (1),  
\[
\eta(xy) = \eta((xz^{-1})(zy)) = \eta(xz^{-1}) \eta(zy), \quad \text{for all } z \notin C_{x,y} \cup xA^{-1}.
\]  
If further \( z^{-1} \notin B_x \), then \( \eta(xz^{-1}) = \eta(x) \eta(z^{-1}) \). Similarly, if \( zy \notin A \) and \( z^{-1} \notin B_{y^{-1}}\), then  
\[
\eta(zy) = \eta(y^{-1}z^{-1})^{-1} = (\eta(y^{-1})\eta(z^{-1}))^{-1} = \eta(z^{-1})^{-1} \eta(y^{-1})^{-1}.
\]  
Defining the finite set  
\[
D_{x,y} := C_{x,y} \cup xA^{-1} \cup B_x^{-1} \cup Ay^{-1} \cup B_{y^{-1}}^{-1},
\]  
we conclude that for all \( z \notin D_{x,y} \),  
\[
\eta(xy) = \eta(x) \eta(z^{-1}) \eta(z^{-1})^{-1} \eta(y^{-1})^{-1} = \eta(x) \eta(y^{-1})^{-1}.
\]  
Since \( G \) is infinite, \( D_{x,y} \neq G \), and thus the identity \eqref{xyy-1} holds. This completes the proof.
\end{proof}

\subsection{Relations for $u$ and $v$}

We collect here some relations for the functions $u$ and $v$ which will be used in the proof of Theorem \ref{Thm_Fields}. Recall that the map \( w : \mathbb{K}^* \to \mathbb{L}^* \) is defined by
\[
w(x) = \frac{v(x)}{u(x)}, \quad x \in \mathbb{K}^*).
\]

\begin{lemma}
\label{Lemma_Inversew}
There is a finite set $A_o \subset \bK^*$ such that $w\big(\frac{1}{x}\big) = \frac{1}{w(x)}$ for all $x \notin A_o$.
\end{lemma}

\begin{proof}
Fix \( x \in \bK^* \) and observe that  
\[
\rho(1+x) = \rho(x) \rho\Big(1 + \frac{1}{x}\Big) = \rho(x) \, u\Big(\frac{1}{x}\Big) + v\Big(\frac{1}{x}\Big) = u(x) + v(x) \rho(x).
\]
Thus,  
\[
a(x) + b(x) \rho(x) = 0, \quad \text{for all } x \in \bK^*,
\]
where  
\[
a(x) = v\Big(\frac{1}{x}\Big) - u(x), \quad b(x) = u\Big(\frac{1}{x}\Big) - v(x).
\]
Since \( u \) and \( v \) take finitely many values, so do \( a \) and \( b \). Let \( Q \subset (\bL^*)^2 \setminus \{(0,0)\} \) be the (possibly empty) finite set of values for \( (a,b) \), and enumerate its elements as  
\[
Q = \{ (\alpha_1,\beta_1),\ldots,(\alpha_n,\beta_n) \}.
\]
Define  
\[
E_j = \{ x \in \bK^* \mid a(x) = \alpha_j, \, b(x) = \beta_j \}, \quad j=1,\dots,n,
\]
and  
\[
E_0 = \{ x \in \bK^* \mid a(x) = b(x) = 0 \}.
\]
Then \( \bK^* = \bigcup_{j=0}^n E_j \), with \( \bK^* = E_0 \) if \( Q \) is empty. \\

Assuming \( Q \neq \emptyset \), if any \( E_j \) contains two distinct elements \( x_1, x_2 \), then  
\[
\alpha_j + \beta_j \rho(x_1) = \alpha_j + \beta_j \rho(x_2) = 0,
\]
implying \( \alpha_j = \beta_j = 0 \), a contradiction. Hence, each \( E_j \) contains at most one element, so the union \( A_o = E_1 \cup \dots \cup E_n \) has at most \( n \) points. For \( x \notin A_o \), we have \( a(x) = b(x) = 0 \), yielding  
\[
v\Big(\frac{1}{x}\Big) = u(x), \quad u\Big(\frac{1}{x}\Big) = v(x).
\]
Thus, \( w\big(\frac{1}{x}\big) = \frac{1}{w(x)} \) for all \( x \notin A_o \).
\end{proof}

\begin{lemma}
\label{Lemma_wminus}
There is a finite set $A_1 \subset \bK^*$ such that $w(-x) = w(x)$ for all 
$x \notin A_1$.
\end{lemma}

\begin{proof}
For \( x \in \bK^* \), we can express \( \rho(1 - x^2) \) in two ways:  
\[
\rho(1 - x^2) = u(-x^2) - v(-x^2) \rho(x^2),
\]
or, alternatively, using that \( \rho(1 - x^2) = \rho(1 + x) \rho(1 - x) \), 
\[
\rho(1 - x^2) = u(x) u(-x) - (u(x) v(-x) - u(-x) v(x)) \rho(x) - v(x) v(-x) \rho(x)^2.
\]
Equating these expressions yields  
\[
a(x) + b(x) \rho(x) + c(x) \rho(x)^2 = 0, \quad \forall x \in \bK^*,
\]
where  
\[
a(x) = u(x) u(-x) - u(-x^2), \quad b(x) = u(-x) v(x) - u(x) v(-x), \quad c(x) = v(-x^2) - v(x) v(-x).
\]
Since \( u \) and \( v \) take finitely many values, so do \( a, b, c \). Fixing any  \((\alpha, \beta, \gamma) \in \bL^* \times \bL^* \times \bL^*\), the set  
\[
E = \{ x \in \bK^* \mid a(x) = \alpha, \, b(x) = \beta, \, c(x) = \gamma \}
\]
contains at most two elements. Indeed, assume for the sake of contradiction that $E$
contains three distinct points $x_1,x_2,x_3$. Then,  
\[
\left(
\begin{matrix}
1 & \rho(x_1) & \rho(x_1)^2 \\
1 & \rho(x_2) & \rho(x_2)^2 \\
1 & \rho(x_3) & \rho(x_3)^2
\end{matrix}
\right)
\left(
\begin{matrix}
\alpha \\ \beta \\ \gamma
\end{matrix}
\right)
=
\left(
\begin{matrix}
0 \\ 0 \\ 0
\end{matrix}
\right)
\]
Since $\rho$ is injective, $\rho(x_1), \rho(x_2)$ and $\rho(x_3)$ are distinct, so the standard fact that Vandermonde determinants are non-zero ensures \((\alpha, \beta, \gamma) = (0,0,0)\), a contradiction. Defining \( A_1 \) as the finite union of all nonempty such \( E \), we obtain 
\[
a(x) = b(x) = c(x) = 0 \quad \textrm{for \( x \notin A_1 \)}, 
\]
which implies \( u(x) v(-x) = u(-x) v(x) \), yielding \( w(x) = w(-x) \) for all \( x \notin A_1 \).  
\end{proof}

\begin{lemma}
\label{Lemma_uv_are_compatible}
There is a finite set $D \subset \bK^*$ such that 
\[
w(1+x) u(x) = 1 \qand v(1+x) v(x) = 1, \quad \textrm{for all $x \notin D$}.
\]
\end{lemma}

\begin{proof}
Fix \( x \in \bK^* \) and compute:  
\begin{align*}
\rho(x) 
&= \rho(x+1-1) = -\rho(1-(1+x)) = -\big(u(1+x) - v(1+x)\rho(1+x) \big) \\[0.2cm]
&= -\big(u(1+x) - v(1+x)(u(x) + v(x) \rho(x)) \big) \\[0.2cm]
&= v(1+x)u(x) - u(1+x) + v(1+x)v(x) \rho(x).
\end{align*}
Thus,  
\[
a(x) + b(x) \rho(x) = 0, \quad \text{for all } x \in \bK^*,
\]
where  
\[
a(x) = v(1+x)u(x) - u(1+x) \quad \textrm{and} \quad b(x) = v(1+x)v(x) - 1.
\]
Following the same argument as in Lemma \ref{Lemma_Inversew}, we conclude that there exists a finite set \( D \subset \bK^* \) such that \( a(x) = b(x) = 0 \) for all \( x \notin D \). Consequently, we obtain  
\[
v(1+x)u(x) = u(1+x) \qand v(1+x)v(x) = 1, \quad \text{for all } x \notin D.
\]
\end{proof}

\subsection{Proof of Theorem \ref{Thm_Fields}}

Let \(\bK, \bL, \rho, u, v, w, \kappa\) be as defined in Theorem \ref{Thm_Fields} above. In this subsection, we derive Theorem \ref{Thm_Fields} from the following lemma, whose proof will be the focus of the remainder of the section.

\begin{lemma}
\label{Lemma_w_is_PWhomo}
The map $w : \bK^* \ra \bL^*$ is a patchwise homomorphism.
\end{lemma}

\begin{proof}[Proof of Theorem \ref{Thm_Fields} assuming 
Lemma \ref{Lemma_w_is_PWhomo}]
(i) Since $\bK^*$ is an infinite group, Lemmas \ref{Lemma_w_is_PWhomo} and \ref{Lemma_Patchwise} imply that $w$ is a homomorphism from $\bK^*$ to $\bL^*$. Moreover, since both $u$ and $v$ are finite-valued, so is $w$. \\

(ii) Since $\rho$ and $w$ are both homomorphisms from $\bK^*$ to $\bL^*$, so is $\kappa$. To show that $\kappa$ is also additive, recall from Lemma \ref{Lemma_uv_are_compatible} that there exists a finite set $D \subset \bK^*$ such that
\[
w(1+x) u(x) = 1  \quad \text{for all } x \notin D.
\]
Thus, for all $x \notin D$,
\begin{align*}
\kappa(1 + x) &= w(1+x) \rho(1+x) = w(1+x)u(x) + w(1+x)v(x) \rho(x) \\[0.2cm]
&= 1 + w(1+x)u(x) w(x) \rho(x) = 1 + \kappa(x).
\end{align*}
For all $x \in \bK^*$ and $y \notin Dx \cup \{0\} \cup \{-x\}$, we obtain
\begin{align*}
\kappa(x+y) &= \kappa(x) \kappa\Big(1 + \frac{y}{x}\Big) = \kappa(x) \Big(1 + \kappa\Big(\frac{y}{x}\Big) \Big) \\[0.2cm]
&= \kappa(x) + \kappa(x) \kappa\Big(\frac{y}{x}\Big) = \kappa(x) + \kappa(y).
\end{align*}
If either $x$ or $y$ is zero, then $\kappa(x+y) = \kappa(x) + \kappa(y)$ trivially. By Lemma \ref{Lemma_wminus}, there exists a finite set $A_1 \subset \bK^*$ such that $w(-x) = w(x)$ for all $x \notin A_1$. Hence, if $x + y = 0$ with $x \notin A_1$, then $\kappa(x+y) = 0$, and  
\[
\kappa(x) + \kappa(-x) = (w(x) - w(-x))\rho(x) = 0.
\]
We conclude that we may assume $-1 \notin D$, and
\[
\kappa(-x) = -\kappa(x) \quad \text{and} \quad \kappa(x+y) = \kappa(x) + \kappa(y) \quad \text{for all } x \notin A_1, \ y \notin Dx.
\]
We claim that $\kappa$ is a patchwise homomorphism between the additive groups $(\bK,+)$ to $(\bL,+)$. Define  
\[
A = A_1, \quad B_x = Dx, \quad C_{x,y} = \Big\{ \frac{dx-y}{1 + d} \ : \ d \in D \Big\}.
\]
Since $D$ is finite and does not contain $-1$, both $B_x$ and $C_{x,y}$ are finite, and $C_{x,y}$ is well-defined for all $x,y \in \bK$. We have already proved that
\[
\kappa(-x) = -\kappa(x) \quad \text{and} \quad \kappa(x+y) = \kappa(x) + \kappa(y) \quad \text{for all } x \notin A, \ y \notin B_x,
\]
so it remains to show that $z + y \notin B_{x-z}$ for all $z \notin C_{x,y}$. Suppose otherwise, so that there exists an element $d \in D$ such that $z + y = d(x-z)$. Then,  
\[
z = \frac{dx-y}{1 + d} \in C_{x,y},
\]
contradicting our assumption. Thus, $\kappa$ is a patchwise homomorphism. Since $\bK$ is infinite, Lemma \ref{Lemma_Patchwise} implies that $\kappa$ is additive.

To prove injectivity of $\kappa$, note that since neither $w$ nor $\rho$ attains $0$, we have $\kappa(x) = 0$ if and only if $x = 0$. Since $\kappa$ is additive, this establishes injectivity, making $\kappa$ an injective field homomorphism.

Finally, to show that $\kappa$ is a field isomorphism when $\rho$ is surjective, pick $t \in \bL$. We seek $x \in \bK$ such that $\kappa(x) = t$. Let $\Gamma = \ker(w)$. Since $w$ is finite-valued, $\Gamma$ is a finite-index subgroup of the multiplicative group $\bK^*$, so $\rho(\Gamma)$ is a finite-index subgroup of the multiplicative group $\bL^*$, as $\rho$ is surjective.  
By \cite[Theorem 1.7]{BS}, we have $\rho(\Gamma) - \rho(\Gamma) = \bL$, so there exist $\gamma_1, \gamma_2 \in \Gamma$ such that $t = \rho(\gamma_1) - \rho(\gamma_2)$. Setting $x = \gamma_1 - \gamma_2$, we obtain  
\begin{align*}
\kappa(x) &= \kappa(\gamma_1 - \gamma_2) = \kappa(\gamma_1) - \kappa(\gamma_2) \\
&= w(\gamma_1) \rho(\gamma_1) - w(\gamma_2) \rho(\gamma_2) \\
&= \rho(\gamma_1) - \rho(\gamma_2) = t.
\end{align*}
Thus, $\kappa$ is a field isomorphism.
\end{proof}

\subsection{Additive relations for $\rho$}

We now proceed with the proof of Lemma \ref{Lemma_w_is_PWhomo}. The central component of the argument is the following technical lemma, which roughly asserts that there are very few additive relations among $\rho(x)$, $\rho(y)$, and $\rho(xy)$ for $x,y \in \bK^*$.

\begin{lemma}
\label{Lemma_alpha}
For \( (\alpha,\beta,\gamma,\delta) \in (\bL^*)^4 \setminus \{(0,0,0,0)\} \), define  
\[
E = \{ (x,y) \in \bK^* \times \bK^* \mid \alpha + \beta \rho(x) + \gamma \rho(y) + \delta \rho(x) \rho(y) = 0 \}.
\]  
Then either $E$ is empty, or \( E \) takes one of two forms:  
\vspace{0.2cm}
\begin{enumerate}
\item[$(i)$] \textsc{Graph of an Injective Function:}  
There exists a nonempty set \( S \subset \bK^* \) and an injective map 
\( S \to \bL^*, \, x \mapsto y_x \) such that  
\[
E = \{(x, y_x) \mid x \in S\}.
\] 
Furthermore, we have $\gamma + \delta \rho(x) \neq 0$ and $\alpha + \beta \rho(x) \neq 0$ for all $x \in S$, and $x \mapsto y_x$ is uniquely determined by the formula:
\[
\rho(y_x) = -\frac{\alpha + \beta \rho(x)}{\gamma + \delta \rho(x)}.
\]

\item[$(ii)$] \textsc{Union of Vertical and/or Horizontal Lines:}  
There exist \( x_o, y_o \in \bK^* \) such that  
\[
E = \{x_o\} \times \bK^* \quad \text{or} \quad E = \bK^* \times \{y_o\}
\quad \text{or} \quad E = (\{x_o\} \times \bK^*) \cup (\bK^* \times \{y_o\}).
\]
\end{enumerate}
\end{lemma}

\begin{proof}
Assume that $E$ is non-empty, and pick $x \in \bK^*$ such that  
\[
E_x = \{ y \in \bK^* \mid (x,y) \in E \}
\]
is non-empty. We now distinguish between two cases.

\subsection*{Case 1: $|E_x| \geq 2$}
Suppose there exist distinct elements $y_1, y_2 \in \bK^*$ such that $(x, y_1), (x, y_2) \in E$. Then, by the defining equation of $E$,
\begin{align*}
(\alpha + \beta \rho(x)) + (\gamma + \delta \rho(x)) \rho(y_1) &= 0, \\
(\alpha + \beta \rho(x)) + (\gamma + \delta \rho(x)) \rho(y_2) &= 0.
\end{align*}
Subtracting these two equations yields
\[
(\gamma + \delta \rho(x)) (\rho(y_1) - \rho(y_2)) = 0.
\]
Since $ y_1 \neq y_2 $ and $\rho$ is injective, it follows that $\rho(y_1) \neq \rho(y_2)$, which forces  
\[
\tag{$\alpha\beta\gamma\delta$}
\gamma + \delta \rho(x) = 0, \quad \textrm{and thus}, \quad \alpha + \beta \rho(x) = 0.
\]
Hence, for all $ y \in \bK^* $,  
\[
(\alpha + \beta \rho(x)) + (\gamma + \delta \rho(x)) \rho(y) = 0,
\]
which implies $ E_x = \bK^* $. \\

Since we assume that $(\alpha, \beta, \gamma, \delta) \neq (0,0,0,0)$ and $\rho$ takes values in $\bL^*$, at least one of $\beta$ or $\delta$ must be nonzero; otherwise the equations ($\alpha\beta\gamma\delta)$ would force $\alpha = \beta = \gamma = \delta = 0$. This ensures that either $\rho(x) = -\alpha/\beta$ or $\rho(x) = -\gamma/\delta$, uniquely determining $ x $ in terms of $\alpha, \beta, \gamma, \delta$. Consequently, there is at most one $ x \in \bK^* $ such that $ |E_x| \geq 2 $, and for this $ x $, we have $ E_x = \bK^* $.

\subsection*{Case 2: $|E_x| = 1$}
If $ E_x $ contains only one element, we denote it by $ y_x $, which is uniquely determined by the equation
\[
(\alpha + \beta \rho(x)) + (\gamma + \delta \rho(x)) \rho(y_x) = 0.
\]
The uniqueness of $ y_x $ requires that $\gamma + \delta \rho(x) \neq 0$, while the existence of $ y_x $ ensures $\alpha + \beta \rho(x) \neq 0$, leading to the explicit formula:
\[
\rho(y_x) = -\frac{\alpha + \beta \rho(x)}{\gamma + \delta \rho(x)}.
\]

\subsection*{Remark on Case 2: Structure of $E$}
Define  
\[
S = \{ x \in \bK^* \mid |E_x| = 1 \}.
\]
We claim that the function  
\[
S \to \bK^*, \quad x \mapsto y_x
\]
is either \textbf{injective} or \textbf{constant}. To see this, suppose it is not injective, so there exist distinct $ x_1, x_2 \in S $ such that $ y_{x_1} = y_{x_2} $. Then,
\begin{align*}
(\alpha + \gamma \rho(y_{x_1})) + (\beta + \delta \rho(y_{x_1})) \rho(x_1) &= 0, \\
(\alpha + \gamma \rho(y_{x_1})) + (\beta + \delta \rho(y_{x_1})) \rho(x_2) &= 0.
\end{align*}
Hence, since $\rho(x_1) \neq \rho(x_2)$, 
\[
\alpha + \gamma \rho(y_{x_1})  = \beta + \delta \rho(y_{x_1}) = 0,
\]
and thus, for all $ x \in \bK^* $,
\[
(\alpha + \gamma \rho(y_{x_1})) + (\beta + \delta \rho(y_{x_1})) \rho(x) = 0,
\]
which implies $ (x, y_{x_1}) \in E $ for all $ x \in \bK^* $, so the map $x \mapsto y_x $ is constant (equal to $y_{x_1}$) on $ S $.

\subsection*{Summary: Possible Forms of $E$}

We conclude that the set \( E \) must take one of the following three forms:

\begin{enumerate}
    \item \( E = \{x_o\} \times \bK^* \) for some \( x_o \in \bK^* \). \vspace{0.2cm}
    
    \item There exists a non-empty subset \( S \subset \bK^* \) such that  
    \[
    E = \{ (x,y_x) \mid x \in S\},
    \]
    where the mapping \( x \mapsto y_x \) is either injective or constant, and is uniquely determined by 
    \[
\rho(y_x) = -\frac{\alpha + \beta \rho(x)}{\gamma + \delta \rho(x)}.
\]
    
    \item There exist \( x_o \in \bK^* \) and a non-empty subset \( S \subset \bK^* \setminus \{x_o\} \) such that  
    \[
    E = (\{x_o\} \times \bK^*) \cup \{ (x,y_x) \mid x \in S\}.
    \]
\end{enumerate}

To further analyze these cases, we examine the structure of the \( y \)-sections of \( E \), given by
\[
E^y = \{ x \in \bK^* \mid (x,y) \in E \}.
\]
For every \( y \in \bK^* \), arguing exactly as above, we see that the set \( E^y \) is either empty, contains exactly one element, or equals \( \bK^* \). This classification leads to the following refinements: 

\begin{itemize}
    \item If the map \( x \mapsto y_x \) is injective, then Case 3 is impossible. Indeed, if we select any \( x \in S \) and set \( y = y_x \), we obtain \( E^y = \{x_o, x\} \). However, since \( x_o \notin S \) and \( \bK^* \) is infinite, this contradicts the structure of \( y \)-sections described above. Thus, Case 3 cannot occur under injectivity. \vspace{0.2cm}
    
    \item If \( x \mapsto y_x \) is constant on \( S \), say \( y_x = y_o \) for all \( x \in S \), then Cases 2 and 3 can be rewritten as
    \[
    E = S \times \{y_o\} \quad \text{or} \quad E = (\{x_o\} \times \bK^*) \cup (S \times \{y_o\}).
    \] 
    Considering the \( y \)-section at \( y = y_o \), we see that in the first case, \( S \) must either be a singleton, say \( S = \{x_1\} \) for some \( x_1 \in \bK^* \), or the entire set \( \bK^* \). In the second case, since \( E^{y_o} = \{x_o\} \sqcup S \) and \( |S| \geq 2 \), we must have \( S = \bK^* \setminus \{x_o\} \). Thus, if \( x \mapsto y_x \) is constant on \( S \), Cases 2 and 3 reduce to three possibilities:
    \[
    E = \{x_1\} \times \{y_o\}, \quad E = \bK^* \times \{y_o\}, \quad E = (\{x_o\} \times \bK^*) \cup (\bK^* \times \{y_o\}).
    \]
    However, the first case is already contained within Case 2 by taking \( S = \{x_1\} \) and setting \( y_{x_1} = y_o \), which trivially satisfies injectivity.
\end{itemize}
\end{proof}

\subsection{Proof of Lemma \ref{Lemma_w_is_PWhomo}}

We will need the following lemma. 

\begin{lemma}
\label{Lemma_abcd}
For all \((x,y) \in \bK^* \times \bK^*\), the following identity holds:
\[
a(x,y) + b(x,y) \rho(x) + c(x,y) \rho(y) + d(x,y) \rho(x) \rho(y) = 0,
\]
where the coefficients are given by
\begin{align*}
a(x,y) &= u(x) u(y) - u(xy + x + y), \\[0.2cm]
d(x,y) &= v(x)v(y) - v(xy + x + y) u\Big(\frac{x+y}{xy}\Big),
\end{align*}
and
\begin{align*}
b(x,y) &= v(x)u(y) - v(xy + x + y) v\Big(\frac{x+y}{xy}\Big) u\Big(\frac{y}{x}\Big), \\[0.2cm]
c(x,y) &= u(x)v(y) - v(xy + x + y) v\Big(\frac{x+y}{xy}\Big) v\Big(\frac{y}{x}\Big).
\end{align*}
\end{lemma}

\begin{proof}
Fix \( x, y \in \bK^* \) and express \( \rho(1+x)\rho(1+y) \) in two different ways. Expanding directly using the definition of \( \rho \), we obtain
\begin{align*}
\rho(1+x) \rho(1+y) &= (u(x) + v(x)\rho(x))(u(y) + v(y)\rho(y)) \\[0.2cm]
&= u(x)u(y) + v(x)u(y) \rho(x) + u(x)v(y)\rho(y) + v(x)v(y) \rho(x) \rho(y).
\end{align*}

Alternatively, using the identity \( (1+x)(1+y) = 1 + xy + x + y \), we write:
\begin{align*}
\rho(1 + xy + x + y) &= u(xy + x + y) + v(xy + x + y) \rho(xy + x + y) \\[0.2cm]
&= u(xy + x + y) + v(xy + x + y) \rho(xy) \rho\Big(1 + \frac{x+y}{xy}\Big) \\[0.2cm]
&= u(xy + x + y) + v(xy + x + y) u\Big(\frac{x+y}{xy}\Big) \rho(xy) \\[0.2cm]
&\quad + v(xy + x + y) v\Big(\frac{x+y}{xy}\Big) \rho(x+y) \\[0.2cm]
&= u(xy + x + y) + v(xy + x + y) u\Big(\frac{x+y}{xy}\Big) \rho(xy) \\[0.2cm]
&\quad + v(xy + x + y) v\Big(\frac{x+y}{xy}\Big) u\Big(\frac{y}{x}\Big) \rho(x) \\[0.2cm]
&\quad + v(xy + x + y) v\Big(\frac{x+y}{xy}\Big) v\Big(\frac{y}{x}\Big) \rho(y).
\end{align*}

Equating these expressions yields
\[
a(x,y) + b(x,y) \rho(x) + c(x,y) \rho(y) + d(x,y) \rho(x) \rho(y) = 0,
\]
where the coefficients \( a(x,y), b(x,y), c(x,y), d(x,y) \) are as stated in the lemma.
\end{proof}

\begin{proof}[Proof of Lemma \ref{Lemma_w_is_PWhomo}]
Since $u$ and $v$ take only finitely many values, the coefficient functions $a, b, c, d$ in Lemma~\ref{Lemma_abcd} also take on finitely many values. Define  
\[
Q \subset (\bL^*)^4 \setminus \{(0,0,0,0)\}
\]
as the (possibly empty) finite set of values of $(a,b,c,d)$, excluding $(0,0,0,0)$. If $Q \neq \emptyset$, enumerate its elements as  
\[
Q = \{ (\alpha_1,\beta_1,\gamma_1,\delta_1),\ldots,(\alpha_n,\beta_n,\gamma_n,\delta_n) \}.
\]
Now, define the sets  
\[
E_j = \{ (x,y) \in \bK^* \times \bK^* \mid \alpha_j + \beta_j \rho(x) + \gamma_j \rho(y) + \delta_j \rho(x) \rho(y) = 0 \}, \quad j = 1,\dots,n,
\]
and  
\[
E_0 = \{ (x,y) \in \bK^* \times \bK^* \mid a(x,y) = b(x,y) = c(x,y) = d(x,y) = 0 \}.
\]
These sets form a partition:
\[
\bK^* \times \bK^* = \bigcup_{j=0}^n E_j,
\]
where $\bK^* \times \bK^* = E_0$ if $Q = \emptyset$. By Lemma~\ref{Lemma_alpha}, each set $E_j$ for $j=1,\dots,n$ takes one of the following forms:
\begin{enumerate}
    \item[$(i)$] There exists a non-empty set $S_j \subset \bK^*$ and an injective map $x \mapsto y_x^{(j)}$ such that  
    \[
    E_j = \{ (x, y_x^{(j)}) \mid x \in S_j \}.
    \]
    Moreover, $x \mapsto y_x^{(j)}$ is uniquely determined by 
    \[
    \rho(y^{(j)}_x) = -\frac{\alpha_j + \beta_j \rho(x)}{\gamma_j + \delta_j \rho(x)}, \quad x \in S_j.
    \]
    \item[$(ii)$] There exist $x_j, y_j \in \bK^*$ such that  
    \[
    E_j = \{x_j\} \times \bK^*, \quad E_j = \bK^* \times \{y_j\}, \quad \text{or} \quad E_j = (\{x_j\} \times \bK^*) \cup (\bK^* \times \{y_j\}).
    \]
\end{enumerate}
Let $A'$ and $B_o$ be the sets of all such $x_j$ and $y_j$ from (ii), respectively, each containing at most $n$ elements. For $x \notin A'$, define  
\[
J_x = \{ j=1,\ldots,n \mid x \in S_j, \text{ with } E_j \text{ of type (i)} \},
\]
and set  
\[
B_o(x) = B_o \cup \{ y_x^{(j)} \mid j \in J_x \}.
\]  
Clearly, $|B_o(x)| \leq 2n$. If $x \notin A'$ and $y \notin B_o(x)$, then $(x,y) \in E_0$, implying  
\[
a(x,y) = b(x,y) = c(x,y) = d(x,y) = 0.
\]
In particular, from Lemma \ref{Lemma_abcd} we see that the equations $b(x,y) = c(x,y) = 0$ yield  
\begin{align*}
v(x)u(y) &= v(xy + x + y) v\Big(\frac{x+y}{xy}\Big) u\Big(\frac{y}{x}\Big), \\
u(x)v(y) &= v(xy + x + y) v\Big(\frac{x+y}{xy}\Big) v\Big(\frac{y}{x}\Big).
\end{align*}
Dividing these identities gives  
\begin{equation} \label{wdiv}
\frac{w(y)}{w(x)} = w\Big(\frac{y}{x}\Big), \quad \text{for all } x \notin A', \quad y \notin B_o(x).
\end{equation}
By Lemma~\ref{Lemma_Inversew}, there exists a finite (symmetric) set $A_o \subset \bK^*$ such that  
\begin{equation} \label{wC}
w\Big(\frac{1}{x}\Big) = \frac{1}{w(x)}, \quad \text{for all } x \notin A_o.
\end{equation}
Hence, if we let 
\[
A = \left\{ \frac{1}{x} \, \mid \, x \notin A_o \cup A' \right\} \qand B_x = B_o(1/x), \quad \textrm{for $x \notin A$}
\]
then $A$ is a finite set, and $B_x$ is a finite set for all $x \notin A$, and 
\[
w\left(\frac{1}{x}\right) = \frac{1}{w(x)} \qand w(xy) = w(x) w(y), \quad \textrm{for all $x \notin A$ and $y \notin B_x$}.
\]
To prove that $w : \bK^* \ra \bL^*$ is a patchwise homomorphism it thus remains to show that for every pair \((x,y) \in \bK^* \times \bK^*\), there exists a finite set \(C_{x,y} \subset \bK^*\) such that  
\[
zy \notin B_{xz^{-1}} \quad \text{for all } z \notin C_{x,y}.
\]  
Fix \((x,y) \in \bK^* \times \bK^*\). To construct \(C_{x,y}\), suppose that $z \in \bK^*$ and \(zy \in B_{x z^{-1}}\). By the definition of \(B_{xz^{-1}}\), this can occur in one of two ways:
\vspace{0.1cm}
\begin{enumerate}
    \item[(i)] $zy \in B_o$, or equivalently, $z \in B_o y^{-1}$. \vspace{0.2cm}
    
    \item[(ii)] There exists \(j \in J_{1/xz^{-1}}\) such that \(1/xz^{-1} \in S_j\) and  
    \[
    zy = y_{1/xz^{-1}}^{(j)}.
    \]  
    This leads to the equation  
    \[
    \rho(zy) = - \frac{\alpha_j + \beta_j \rho(1/xz^{-1})}{\gamma_j + \delta_j \rho(1/xz^{-1})} = - \frac{\alpha_j \rho(xz^{-1}) + \beta_j}{\gamma_j \rho(xz^{-1}) + \delta_j} = - \frac{\alpha_j \rho(x) + \beta_j \rho(z)}{\gamma_j \rho(x) + \delta_j \rho(z)}.
    \]
    Equivalently,  
    \[
    \delta_j \rho(y) \rho(z)^2 + (\gamma_j \rho(x) \rho(y) + \beta_j) \rho(z) + \alpha_j \rho(x) = 0.
    \]  
    Let \(D_{x,y}^{(j)}\) be the (possibly empty) set of all solutions \(z \in \bK^*\) to this quadratic equation for given \(j, x, y\). Clearly, we have \(|D_{x,y}^{(j)}| \leq 2\).  
\end{enumerate}
\vspace{0.1cm}
Now, define  
\[
C_{x,y} = B_o y^{-1} \cup \left( \bigcup_{j=1}^n D_{x,y}^{(j)} \right).
\]  
Since \(|D_{x,y}^{(j)}| \leq 2\), we have  
\[
|C_{x,y}| \leq |B_o| + 2n < \infty \quad \text{for all } (x,y) \in \bK^* \times \bK^*,
\]  
and, for all \(z \notin C_{x,y}\), we have \(zy \notin B_{xz^{-1}}\), completing the proof.
\end{proof}

\section{Existence of reduced models: Proof of Theorem \ref{Thm_ReducedModels}}
\label{AppendixB}

Let \(\bK\) be a countably infinite field of positive characteristic \(p = \Char(\bK)\), and let \((U, \rho_U)\) be a diagonalizable \(\bK^*\)-representation on a \(\bK\)-vector space \(U\), with weight spaces \(U_1, \ldots, U_r\) corresponding to weights \(m_1, \ldots, m_r\). \\

For each \(j = 1, \ldots, r\), write the weight \(m_j\) as \(m_j = p^{l_j} n_j\), where \(l_j \geq 0\) and \(n_j\) is not divisible by \(p\). Define the subfield \(\bK_j \subset \bK\) by
\[
\bK_j := \{ x^{p^{l_j}} \mid x \in \bK \}.
\]
Regard \(U_j\) as a \(\bK_j\)-vector space, and fix a Hamel basis \((e_{j,\alpha})_{\alpha \in I_j}\) for \(U_j\) over \(\bK_j\), where \(I_j\) is an appropriate index set. Let \(V_j\) be a \(\bK\)-vector space with a Hamel basis \((f_{j,\alpha})_{\alpha \in I_j}\), indexed by the same set. Define the map \(\Phi_j : V_j \to U_j\) by
\[
\Phi_j\left(\sum_{\alpha \in I_j} x_\alpha f_{j,\alpha}\right)
= \sum_{\alpha \in I_j} x_\alpha^{p^{l_j}} e_{j,\alpha}, \quad x_\alpha \in \bK.
\]
Since \(\bK\) has characteristic \(p\), the map \(\Phi_j\) is additive (but nonlinear unless \(l_j = 0\)). \\

To show that \(\Phi_j\) is an isomorphism, first note that it is injective: if \(\Phi_j(x) = 0\), then
\[
\sum_{\alpha \in I_j} x_\alpha^{p^{l_j}} e_{j,\alpha} = 0
\quad \Rightarrow \quad x_\alpha^{p^{l_j}} = 0 \text{ for all } \alpha
\quad \Rightarrow \quad x_\alpha = 0 \text{ for all } \alpha
\quad \Rightarrow \quad x = 0.
\]
To see that \(\Phi_j\) is also surjective, let \(y \in U_j\) and write
\[
y = \sum_{\alpha \in I_j} y_\alpha e_{j,\alpha}, \quad y_\alpha \in \bK_j.
\]
For each \(\alpha \in I_j\), there exists a unique \(x_\alpha \in \bK\) such that \(y_\alpha = x_\alpha^{p^{l_j}}\). Then
\[
x := \sum_{\alpha \in I_j} x_\alpha f_{j,\alpha} \in V_j \quad \text{ satisfies } \quad \Phi_j(x) = y,
\]
proving surjectivity. Now define
\[
V := \bigoplus_{j=1}^r V_j,
\]
and let \(\Phi : V \to U\) be the unique additive isomorphism such that \(\Phi|_{V_j} = \Phi_j\) for each \(j = 1, \ldots, r\). \\

Next, observe that for all \(a \in \bK^*\) and \(x \in V_j\),
\[
\Phi_j(a^{n_j} x) = \sum_{\alpha \in I_j} (a^{n_j} x_\alpha)^{p^{l_j}} e_{j,\alpha}
= a^{m_j} \sum_{\alpha \in I_j} x_\alpha^{p^{l_j}} e_{j,\alpha}
= a^{m_j} \Phi_j(x).
\]
This shows that \(\Phi\) intertwines \(\rho_U\) with the reduced representation \((V, \rho_V)\), defined by
\[
\rho_V(a) \cdot v_j := a^{n_j} v_j, \quad \text{for } v_j \in V_j, \quad a \in \bK^*.
\]
Hence, \((V, \rho_V)\) is a reduced model for \((U, \rho_U)\).


\begin{thebibliography}{99}
\bibitem{Bekka}
Bekka, B.
\emph{Infinite characters on $\GL_n(\bQ)$, on $\SL_n(\bZ)$, and on groups acting on trees.}
J. Funct. Anal. 277 (2019), no. 7, 2160--2178.
\bibitem{BLM}
Bergelson, V.; Leibman, A.; McCutcheon, R.
\emph{Polynomial Szemer\'edi theorems for countable modules over integral domains and finite fields.} 
J. Anal. Math. 95 (2005), 243--296.
\bibitem{BM}
Bergelson, V.; Moreira, J.
\emph{Ergodic theorem involving additive and multiplicative groups of a field and 
$\{x+y,xy\}$ patterns.} 
Ergodic Theory Dynam. Systems 37 (2017), no. 3, 673--692.
\bibitem{BS}
Bergelson, V.; Shapiro, D.B.
\emph{Multiplicative subgroups of finite index in a ring.}
Proc. Amer. Math. Soc. 116 (1992), no. 4, 885--896.
\bibitem{BB}
Bj\"orklund, M.; Bulinski, K.
\emph{Twisted patterns in large subsets of $\bZ^N$}. 
Comment. Math. Helv. 92 (2017), no. 3, 621--640.
\bibitem{BF1}
Bj\"orklund, M.; Fish, A.
\emph{Characteristic polynomial patterns in difference sets of matrices.} 
Bull. Lond. Math. Soc. 48 (2016), no. 2, 300--308.
\bibitem{BT}
Bogomolov, F; Tschinkel, Y.
\emph{Milnor $K_2$ and field homomorphisms.} 
Surveys in differential geometry. Vol. XIII. Geometry, analysis, and algebraic geometry: forty years of the Journal of Differential Geometry, 223--244.
Surv. Differ. Geom., 13
International Press, Somerville, MA, 2009
\bibitem{CLMS}
Cornelissen, G; Li, X.; Marcolli, M.; Smit, H.
\emph{Reconstructing global fields from dynamics in the abelianized Galois group.}
Selecta Math. (N.S.) 25 (2019), no. 2, Paper No. 24, 18 pp.
\bibitem{CdSLMS}
Cornelissen, G.; de Smit, B.; Li, X.; Marcolli, M.; Smit, H.
\emph{Characterization of global fields by Dirichlet L-series.} 
Res. Number Theory 5 (2019), no. 1, Paper No. 7, 15 pp.
\bibitem{CM}
Cornelissen, G.; Marcolli, M.
\emph{Quantum statistical mechanics, L-series and anabelian geometry I: Partition functions.} Trends in contemporary mathematics, 47--57.
\bibitem{Davies}
Davies, J.
\emph{Chromatic number of spacetime.} 
Acta Arith. 218 (2025), no. 1, 65--76.
\bibitem{D}
Deligne,  P.
\emph{La conjecture de Weil. II}, 
Inst. Hautes Etudes Sci. Publ.
Math. no. 52 (1980), 137--252.
\bibitem{EL}
Einsiedler, M.; Lindenstrauss, E.
\emph{Rigidity properties for commuting automorphisms on tori and solenoids.}
Ergodic Theory Dynam. Systems 42 (2022), no. 2, 691--736.
\bibitem{F}
Fish, A.
\emph{On product of difference sets for sets of positive density.} 
Proc. Amer. Math. Soc. 146 (2018), no. 8, 3449--3453.
\bibitem{G}
Glasner, E.
\emph{Ergodic theory via joinings.}
Math. Surveys Monogr., 101
American Mathematical Society, Providence, RI, 2003. xii+384 pp
\bibitem{HIS}
Hart, D.; Iosevich, A.; Solymosi, J.
\emph{Sum-product estimates in finite fields via Kloosterman sums.}
Int. Math. Res. Not. IMRN 2007, no. 5, Art. ID rnm007, 14 pp.
\bibitem{HeRo}
Hewitt, E.; Ross, K.
\emph{Abstract harmonic analysis. Vol. II: Structure and analysis for compact groups. Analysis on locally compact Abelian groups.}
Die Grundlehren der mathematischen Wissenschaften, Band 152
Springer-Verlag, New York-Berlin, 1970. ix+771 pp.
\bibitem{HR}
Howe, R.E.; Rosenberg, J.,
\emph{The unitary representation theory of $\GL(n)$ of an infinite discrete field}.
Israel J. Math. 67 (1989), no. 1, 67--81.
\bibitem{KoKr}
Kosheleva, O.; Kreinovich, V.
\emph{On chromatic numbers of space-times: open problems.}
Geombinatorics 19 (2009), no. 1, 14--17.
\bibitem{L}
Larick, P.G.
\emph{Results in polynomial recurrence for actions of fields}.
PhD-thesis, Ohio State University 1998 \url{http://rave.ohiolink.edu/etdc/view?acc_num=osu1487950153602457}.
\bibitem{Morris}
Morris, S. A.
\emph{Pontryagin duality and the structure of locally compact abelian groups}.
London Math. Soc. Lecture Note Ser., No. 29
Cambridge University Press, Cambridge-New York-Melbourne, 1977. viii+128 pp.
\bibitem{N}
Neukirch, J.
\emph{Kennzeichnung der p-adischen und der endlichen algebraischen Zahlk\"orper.}
Invent. Math. 6 (1969), 296--314.
\bibitem{P}
Pop, F.
\emph{On Grothendieck's conjecture of birational anabelian geometry.}
Ann. of Math. (2) 139 (1994), no. 1, 145--182.
\bibitem{U}
Uchida, K.
\emph{Isomorphisms of Galois groups}.
J. Math. Soc. Japan 28 (1976), no. 4, 617--620.
\bibitem{W}
Weil, A.
\emph{L'int\'egration dans les groupes topologiques et ses applications.} 
Actualit\'es Sci. Indust., No. 869 [Current Scientific and Industrial Topics]
Hermann \& Cie, Paris, 1940. 158 pp.
\end{thebibliography}
\end{document}